\documentclass[a4paper,10pt]{article}

\usepackage[a4paper, centering]{geometry}
\usepackage{float}

\usepackage{cancel}
\usepackage{enumitem}
\usepackage{amsmath}
\usepackage{amsfonts}
\usepackage{amsthm}
\usepackage{amssymb}
\usepackage{upref}
\usepackage{color}
\usepackage[dvipsnames]{xcolor}
\usepackage{graphicx}
\usepackage{hyperref}

\definecolor{darkgreen}{rgb}{0.0, 0.2, 0.13}
\definecolor{darkred}{rgb}{0.55, 0.0, 0.0}

\hypersetup{
	draft=false
	unicode=false,          
	pdftoolbar=true,        
	pdfmenubar=true,        
	pdffitwindow=true,      
	pdfsubject={},   
	pdfnewwindow=true,      
	pdfkeywords={}, 
	colorlinks=true,       
	linkcolor=darkred,          
	citecolor=darkgreen,        
	filecolor=magenta,      
	urlcolor=cyan           
}

\usepackage{array}
\usepackage[utf8]{inputenc}
\usepackage{algorithm,algorithmicx,algpseudocode}

\graphicspath{{./pics/}}

\theoremstyle{plain}
\newtheorem{thm}{Theorem}[section]

\newtheorem{prop}[thm]{Proposition}

\newtheorem{conj}[thm]{Conjecture}
\newtheorem{rem}[thm]{Remark}
\theoremstyle{definition}
\newcommand{\Per}{\operatorname{Per}}
\newcommand{\ds}{\displaystyle}
\newcommand{\di}{\operatorname{div}}
\newcommand{\bo}[1]{{\bf#1}}

\newcommand{\argmin}{\operatorname{argmin}}
\newcommand{\gconv}{\stackrel{\Gamma}{\longrightarrow}}
\usepackage{indentfirst} 

\newcommand{\Blue}[1]{{\color{blue} #1}}

\newcommand{\diam}{\operatorname{diam}}
\newcommand{\Lip}{\operatorname{Lip}}

\newcommand{\footremember}[2]{%
	\footnote{#2}
	\newcounter{#1}
	\setcounter{#1}{\value{footnote}}%
}

\title{Longest minimal length partitions}

\author{Beniamin Bogosel\footremember{cmap}{Centre de Mathématiques Appliquées, Ecole Polytechnique, UMR CNRS 7641, 91128 Palaiseau, France\newline email:
\texttt{beniamin.bogosel@polytechnique.edu}}  
and Edouard Oudet\footremember{grenoble}{Laboratoire Jean Kuntzmann, Univ. Grenoble Alpes, 
	Bâtiment IMAG, BP 53,
	38041 Grenoble, France\newline email:
\texttt{edouard.oudet@univ-grenoble-alpes.fr}}}

\begin{document}
	
	\maketitle
	
	\begin{abstract}
		This article provides numerical evidence that under volume constraint the ball is the set which maximizes the perimeter of the least-perimeter partition into cells with prescribed areas. We introduce a numerical maximization algorithm which performs multiple optimizations steps at each iteration to approximate minimal partitions. Using these partitions we compute perturbations of the domain which increase the minimal perimeter. The initialization of the optimal partitioning algorithm uses capacity-constrained Voronoi diagrams. A new algorithm is proposed to identify such diagrams, by computing the gradients of areas and perimeters for the Voronoi cells with respect to the Voronoi points.
	\end{abstract}

\section{Introduction}

In \cite{longest_shortest_fence} the authors answer a question raised by Polya in \cite{polya} and prove that among planar convex sets of given area the disk maximizes the length of the shortest area-bisecting curve. Denote by $\Omega \subset \Bbb{R}^d$ an open, connected region with Lipschitz boundary. Consider $c \in (0,1)$ and denote with $|\cdot |$ the usual Lebesgue measure (area in 2D, volume in 3D). Given $\Omega$ and $c$, define the \emph{shortest fence} set to be  
\begin{equation}
 SF(\Omega,c) = \argmin  \{ \Per_\Omega(\omega) : \omega \subset \Omega, |\omega| = c |\Omega|\}. 
 \label{eq:shortest_fence}
\end{equation}
In other words, $SF(\Omega,c)$ is one subset $\omega \subset \Omega$ which minimizes the relative perimeter $\Per_\Omega (\omega)$ when the measure $|\omega|$ is fixed to $c|\Omega|$. In the following, the relative perimeter of $SF(\Omega,c)$ is denoted by
 \begin{equation}
 I(\Omega,c) = \Per_\Omega (SF(\Omega,c))= \min\{\Per_\Omega(\omega) : \omega \subset \Omega,  |\omega| = c|\Omega|\}.
 \label{eq:length-shortest-fence}
 \end{equation}
 In the literature the mapping $c \mapsto I(\Omega,c)$ is sometimes called the \emph{isoperimetric profile} of the set $\Omega$. The paper \cite{longest_shortest_fence} cited above solves the problem of maximizing $I(\Omega,c)$ with respect to $\Omega$,
\begin{equation}
 \max_{|\Omega| = v_d} I(\Omega,c),
 \label{eq:longest_shortest_fence}
\end{equation}
in dimension two for $c=1/2$. In the following $v_d$ denotes the volume of the unit ball in $\Bbb{R}^d$. The choice of the volume constraint $|\Omega|=v_d$ does not reduce the generality of the problem, since changing this constant only rescales the solution via a homothety. Classical details regarding the existence of the sets $SF(\Omega,c)$ and the definition of the relative perimeter $\Per_\Omega(\cdot)$ are recalled in the next section.

\begin{figure}[!htp]
	\begin{center}
		\includegraphics[height=0.1\textwidth]{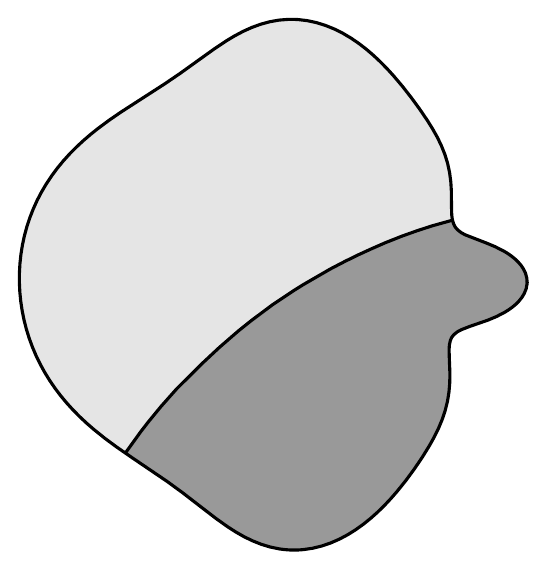}\quad
		\includegraphics[height=0.1\textwidth]{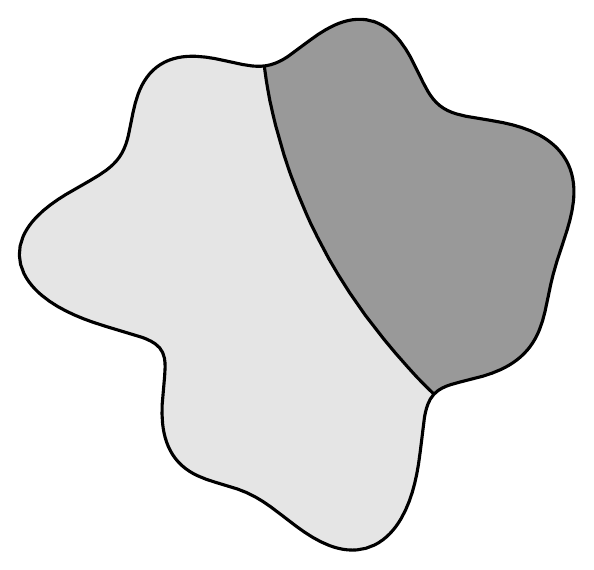}\quad
		\includegraphics[height=0.1\textwidth]{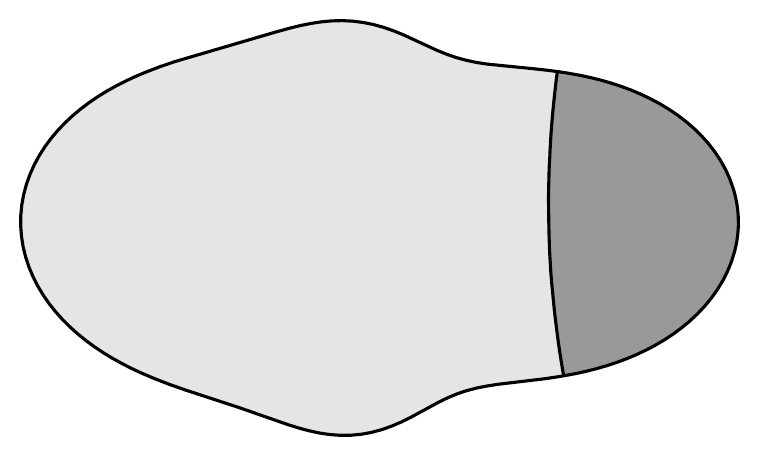}\quad
		\includegraphics[height=0.1\textwidth]{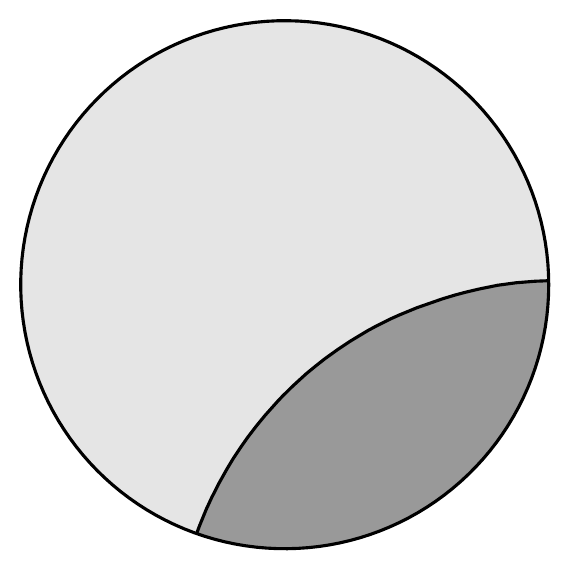}\quad
		\includegraphics[height=0.1\textwidth]{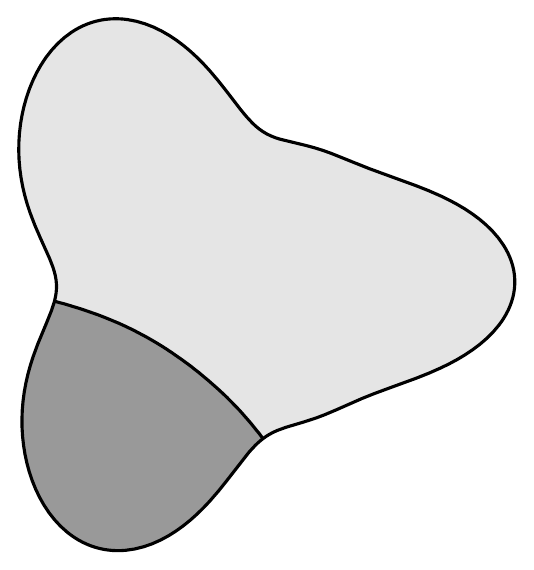}
	\end{center}
	\caption{Examples of minimizers of problem \eqref{eq:shortest_fence} for various shapes $\Omega$ and various constraints.}
	\label{fig:shortest_fences}
\end{figure}

This work was initiated by the note \cite{image_des_maths} published on the French CNRS website Images des mathématiques, where it is asked what happens to the solution of \eqref{eq:longest_shortest_fence} when the parameter $c$ varies in $(0,1)$. This conjecture is attributed to Wichiramala \cite{wichiramala} and the article \cite{morgan_blog} on F. Morgan's blog presents an extensive discussion regarding the history of the problem. The conjecture was partially solved in dimension two in the following works:
\begin{itemize}[topsep=0pt,noitemsep]
	\item In \cite{berry-et-al} the authors prove the conjecture in the plane for small fraction areas. The article contains many interesting results related to relative isoperimetric sets.
	\item In \cite{wang-wang} the authors prove the conjecture in the plane for domains symmetric with respect to both coordinate axes and perturbations of the unit disk.
\end{itemize}

Therefore the conjecture remains unanswered for large fraction areas, except for the case $c=1/2$. Moreover, other generalizations of this problem can be investigated. It is possible, for instance, to consider the analogue problem in the case of partitions of shortest total boundary measure. Given $\Omega \subset \Bbb{R}^d$ and $n>1$ consider $(\omega_1,...,\Blue{\omega_n})$ to be a partition of $\Omega$, in the sense that the union of $\omega_i$, $i=1,...,n$ is $\Omega$ and $\omega_i \cap \omega_j = \emptyset$.  Given a vector $\bo c=(c_1,...,c_n) \in \Bbb{R}^n$ with $\sum_{i=1}^n c_i=1$ consider the \emph{shortest partition} of $|\Omega|$ with \Blue{volume} constraints $\bo c$ to be
\begin{equation}
SP(\Omega,\bo c) = \argmin  \{ \sum_{i=1}^n \Per_\Omega(\omega_i) : (\omega_i) \text{ partition of } \Omega, |\omega_i| = c_i |\Omega|\}. 
\label{eq:shortest_partition}
\end{equation}
In the following we define the \emph{isoperimetric profile of a partition} given by the constraints $\bo c$ by
\begin{equation}
PI(\Omega,\bo c) = \min\{\sum_{i=1}^n \Per_\Omega(\omega_i) : (\omega_i)_{i=1}^n \text{ partition of }\Omega,  |\omega_i| = c_i|\Omega|\}.
\label{eq:length-shortest-partition}
\end{equation}
In other words $PI(\Omega,\bo c)$ is the minimal total relative perimeter of a partition with volume constraints given by $\bo c$. Now it is possible to formulate the problem
\begin{equation}
  \max_{|\Omega| = v_d} PI(\Omega, \bo c),
  \label{eq:longest_shortest_partition}
\end{equation} 
where the total relative perimeter of the shortest partition with constraints $\bo c$ is maximized when $\Omega$ has fixed volume. It is obvious that \eqref{eq:longest_shortest_fence} is a particular case of \eqref{eq:longest_shortest_partition} by considering $n=2$ and $\bo c = (c,1-c)$.

In this paper, problems \eqref{eq:longest_shortest_fence} and \eqref{eq:longest_shortest_partition} are investigated from both numerical and theoretical points of view. In order to approximate solutions of these problems multiple issues need to be addressed:
\begin{itemize}[noitemsep,topsep=0pt]
	\item Compute reliably a numerical approximation of the shortest partition $SP(\Omega,\bo c)$ once the domain $\Omega$ and the constraints vector $\bo c$ are given. It is important to avoid local minimizers at this stage, since the objective is to maximize the shortest perimeter. Any local minimizer may give a false candidate for the solutions of \eqref{eq:longest_shortest_fence} or \eqref{eq:longest_shortest_partition}. There are many works in the literature which deal with the investigation of minimal length partitions. In \cite{cox-partitions} Cox and Flikkema use the surface Evolver software to approximate minimal partitions. In this work the approach presented in \cite{oudet} is used, where the perimeter is approximated using a $\Gamma$-convergence result of Modica and Mortola \cite{modica-mortola}. This allows us to work with density functions rather than sets of finite perimeter and simplifies the handling of the partition condition. Moreover, working with densities directly allows changes in the topology of the partitions. 
	
	Given a domain $\Omega$, a mesh is constructed and finite elements are used in FreeFEM \cite{freefem} in order to approximate $SF(\Omega,c)$ or $SP(\Omega, \bo c)$. When dealing with partitions, in order to accelerate the convergence, an initialization based on Voronoi diagrams with prescribed areas is used. 
	
	\item Once the shortest partition $SP(\Omega, \bo c)$ is identified,  the bounding set $\Omega$ needs to be \emph{modified} in order to increase the  objective function $PI(\Omega,\bo c)$. In order to find a suitable \emph{ascent direction} classical results related to the shape derivative are used \cite{delfour-zolesio, henrot-pierre-english}.
	\item The family of star-shaped domains (which includes convex shapes) is parametrized using radial functions. Moreover, radial functions are discretized by considering truncations of the associated Fourier series. Using the shape derivative it is possible to compute the gradient of the objective function with respect to the discretization parameters. Once the gradient is known, an optimization algorithm is used in order to search for solutions of \eqref{eq:longest_shortest_fence} and \eqref{eq:longest_shortest_partition}. The choice of the optimization algorithm is also an important factor since the computation of $SP(\Omega,\bo c)$ is highly sensitive to local minima. Moreover, when changing $\Omega$ following a perturbation field found using a shape derivative argument, the configuration of the optimal partition might change. The chosen algorithm is a gradient flow with variable step size.  
\end{itemize}

Minimal length partitioning algorithms presented in \cite{oudet} or \cite{bogosel_oudet_manifolds} use random initializations. While this illustrates the flexibility of Modica-Mortola type algorithms and the ability of the algorithm to avoid many local minima, choosing random initializations leads to longer computation times required for the optimization algorithm. A classical idea is to use Voronoi diagrams as initializations. However, these Voronoi diagrams should consist of cells which verify the area constraints $|\omega_i|=c_i$. In the literature, the notion of \emph{capacity-constrained Voronoi diagrams} is employed and results in this direction can be found in \cite{capacity-constrained-discrete}, \cite{capacity-constrained-continuous} and \cite{xin2016}. In this work we propose a new way of computing capacity-constrained Voronoi diagrams by explicitly computing the gradients of the areas of the Voronoi cells with respect to variations in the Voronoi points. The gradient of the perimeters of the Voronoi cells is also computed, which allows the search of capacity-constrained Voronoi diagrams with minimal length. 

The numerical simulations give rise to the following conjectures:
\begin{itemize}[topsep=0pt,noitemsep]
	\item The result of the convex isoperimetric conjecture seems to generalize to every volume fraction $c \in (0,1)$ in dimensions two and three. 
	\item The same result seems to hold in the case of partitions. For $n>1$ and $\bo c = (c_i)_{i=1}^n \in \Bbb{R}^n$ with $\sum_{i=1}^n c_i=1$ arbitrary, the solution of \eqref{eq:longest_shortest_partition} is the disk in 2D and the ball in 3D. It is surprising that this result seems to hold even when the area constraints of the cells of the partition are not the same.
\end{itemize}

{\bf Outline and summary of results.} Section \ref{sec:theory} presents classical theoretical results regarding approximations of minimal perimeter partitions by $\Gamma$-convergence. 

Section \ref{sec:numerics-partitions} recalls basic aspects regarding the numerical computation of minimal length partitions. Section \ref{sec:voronoi} presents the computation of the gradients of the areas and perimeters of Voronoi cells and shows how to use prescribed-area Voronoi cells in order to construct initializations for our optimization algorithm. Section \ref{sec:sh-deriv} presents the computation of an ascent direction for the shape optimization algorithm using the notion of shape derivative. The choice of the discretization and the optimization algorithm for approximating solutions of problems \eqref{eq:longest_shortest_fence} and \eqref{eq:longest_shortest_partition} are presented in Section \ref{sec:optim-algo}. We underline that the maximization algorithm approximates solutions to a max-min problem, and the optimal partitioning algorithm presented in Section \ref{sec:numerics-partitions} is run at every iteration.

Finally, results of the optimization algorithm in dimensions two and three are presented in Section \ref{sec:results}. The numerical results suggest that the solution of problems \eqref{eq:longest_shortest_fence} and \eqref{eq:longest_shortest_partition} is the disk in dimension two and the ball in dimension three. A brief discussion of the optimality conditions is presented in Section \ref{sec:optim-cond}.

\section{Theoretical aspects}
\label{sec:theory}

\subsection{Minimal relative perimeter sets and partitions}
\label{sec:theory-part}

The appropriate framework to work with sets of finite relative perimeter in $\Omega$ is to consider the space of functions with bounded variation on $\Omega$ 
\[ BV(\Omega) = \{ u \in L^1(\Omega) : TV(u) <\infty\}\]
where
\[ TV(u) = \sup\left\{ \int_\Omega  u \di g : g \in C_c(\Omega),  \|g\|_\infty\leq 1\right\}.\]
As usual $C_c(\Omega)$ represents the space of $C^\infty$ functions defined on $\Omega$ with compact support in $\Omega$. 
Using the divergence theorem it is easy to observe that if $u$ is of class $C^1$ then 
\[ TV(u) = \int_\Omega |\nabla u|.\] 
If $\omega$ is a subset of $\Omega$ its generalized perimeter is defined by $\Per(\omega) = TV(\chi_\omega)$, where $\chi_\omega$ represents the characteristic function of $\omega$. All these aspects are classical and can be found, for example, in \cite{ambrosiofuscopallara,braides2}. 

The fact that problems \eqref{eq:shortest_fence} and \eqref{eq:shortest_partition} have solutions is classical and is a consequence of the fact that the generalized perimeter defined above is lower-semicontinuous for the $L^1$ convergence of characteristic functions. For more aspects related to solutions of these problems see \cite[Chapter 17] {maggi-book}. The book previously referenced also presents aspects related to the regularity of optimal partitions in Part Four. Aspects about optimal partitions in the smooth case are presented in \cite{morgan-bubbles} where qualitative properties of minimal  partitions in the plane and on surfaces are presented. 

Proving existence of solutions for problems \eqref{eq:longest_shortest_fence} and \eqref{eq:longest_shortest_partition} is more difficult since these are maximum problems and the perimeter is lower semicontinuous. We recall that problem \eqref{eq:longest_shortest_fence} was solved in  \cite{longest_shortest_fence} in the case $d=2$, $c=1/2$. In particular, existence was proved exploiting results in \cite{cianchi} which show that in this case the minimal relative perimeter sets are convex. In the following we prove that solutions exist in the class of convex domains for arbitrary area constraints. 

\begin{thm}
	Problem \eqref{eq:longest_shortest_fence} has solutions in the class of convex sets, i.e. given $c\in (0,1)$ there exist  convex sets $\Omega^*$ which maximize $I(\Omega,c)$ among convex sets with fixed volume $|\Omega|=v_d$.
	\label{thm:existence_fence}
\end{thm}

\emph{Proof:} We divide the proof into steps which allow us to apply classical methods in calculus of variations.

{\bf Step 1: Upper bounds.} In the following denote by $w(\Omega)$ the minimal $\mathcal H^{d-1}$ measure of the projection of $\Omega$ on a hyperplane (in dimension two this corresponds to the minimal width). For convex bodies the following reverse Loomis-Whitney inequality holds true:
\[ \min_{\{e_1,...,e_d\}} \prod_{i=1}^d \mathcal H^{d-1}(K|e_i^\perp) \leq \Lambda_d |K|^{d-1}, \]
where the minimum is taken over all orthonormal bases of $\Bbb{R}^d$ and $K|e_i^\perp$ represents the projection of $K$ onto a hyperplane orthogonal to $e_i$. In \cite{loomis_whitney} it is shown that there exists a constant $c_0$ such that $\Lambda_d \leq (c_0\sqrt{d})^d$. In particular, this shows that the minimal projection $w(\Omega)$ verifies $w(\Omega)^d \leq \Lambda_d |\Omega|^{d-1}$. As a direct consequence $w(\Omega)$ is bounded above in the class of convex sets $\Omega$ which satisfy $|\Omega|=v_d$.

 It is immediate to see that the quantity $w(\Omega)$ gives an upper bound for $I(\Omega,c)$. To justify this choose $e_0$ the direction for which $\mathcal H^{d-1}(\Omega|e_0^\perp)$ is minimal and slice $\Omega$ with a hyperplane orthogonal to $e_0$ which divides $\Omega$ into two regions $\omega$ and $\Omega\setminus \omega$ with volume $|\omega| = c$. The relative perimeter of the set $\omega$ in $\Omega$ is at most equal to $w(\Omega)$, the $\mathcal H^{d-1}$ measure of the projection. Therefore, we may conclude that in the class of convex sets with measure $|\Omega| = v_d$ the quantity $I(\Omega,c)$ is bounded from above, and the upper bound only depends on $d$ and $v_d$. This implies the existence of a maximizing sequence $(\Omega_h)_{h\geq 1}$ which verifies $I(\Omega_h,c) \leq I(\Omega_{h+1},c)$ and $I(\Omega_h,c) \to \sup_{|\Omega|=v_d} I(\Omega,c)$, where the supremum is taken in the class of convex sets.
 
{\bf Step 2: Compactness.} When dealing with a sequence of convex sets we may extract a subsequence converging in the Hausdorff distance provided the sets are uniformly bounded. For classical aspects related to the Hausdorff distance we refer to \cite[Chapter 2]{henrot-pierre-english}. Therefore, in the following we show that the diameters $\diam(\Omega_h)$ of convex sets $\Omega_h$ forming the maximizing sequence are uniformly bounded. 

First, let us note that since $(\Omega_h)$ is a maximizing sequence for $I(\Omega,c)$ there exists a positive constant $p_0>0$ such that $I(\Omega,c)>p_0$. Since $w(\Omega) \geq I(\Omega,c)$ we also have $w(\Omega_h)\geq p_0>0$ for $n\geq 1$. The results in \cite{Firey1965} show that the minimal perimeter projection, the diameter and the volume of a convex set $\Omega$ satisfy 
\[ w(\Omega) \diam (\Omega) \leq |\Omega|/d.\]
It is now immediate to see that $\diam(\Omega_h) \leq |\Omega_h|/(dw(\Omega))\leq v_d/(dp_0)$, and therefore the diameters of $(\Omega_h)$ are bounded. Without loss of generality we may assume that $(\Omega_h)$ are contained in a large enough ball. Applying the classical Blaschke selection theorem we find that there exists a maximizing sequence, denoted for simplicity by $(\Omega_h)$, such that $\Omega_h$ converges, with respect to the Hausdorff distance, to the convex set $\Omega$ . Moreover, the volume is continuous for the Hausdorff distance among bounded convex sets, so $\Omega$ also satisfies the volume constraint $|\Omega| = v_d$.

{\bf Step 3. Continuity.} The last step is to prove that $I(\Omega,c)$ is indeed equal to $\limsup_{n\to \infty} I(\Omega_h,c)$. This is a direct consequence of \cite[Theorem 4.1]{isoperimetric_convex}, which states that if $(\Omega_h)$ is a sequence of convex bodies in $\Bbb{R}^d$ and $\Omega_h \to \Omega$ in the Hausdorff distance then $I(\Omega_h,c) \to I(\Omega,c)$ for every $c \in [0,1]$. This finishes the proof as the limit $\Omega$ is indeed a maximizer for \eqref{eq:longest_shortest_fence}. \hfill $\square$

\begin{rem}
	Removing the convexity assumption is not straightforward. Nevertheless, using the regularity results regarding solutions of \eqref{eq:shortest_fence} it is possible that this result could be partially extended in the general case. There are multiple difficulties which follow the structure of the proof above:
	\begin{itemize}[noitemsep,topsep=0pt]
		\item Proving there exists an upper bound in \eqref{eq:longest_shortest_fence}.
		\item Proving that a maximizing sequence is bounded: long tails may not intersect the minimizing set in \eqref{eq:shortest_fence} therefore cutting them may increase $I(\Omega,c)$.
		\item Obtaining compactness results of a maximizing sequence: classically this should be possible when working in the class of sets of finite perimeter.
		\item Proving that the maximizing sequence converges to an actual maximizer. This would involve obtaining some continuity properties regarding the perimeter of a sequence of sets. This is not straightforward, as the perimeter is only lower-semicontinuous for the $L^1$ convergence of characteristic functions. Nevertheless, using the regularity of minimal relative perimeter sets might help obtain the desired results.
	\end{itemize}
\end{rem}

The case of partitions can be handled using a similar strategy in the class of convex sets. The missing ingredient is the convergence of the minimal perimeters of partitions, analogue to the results in \cite{isoperimetric_convex}. 

\begin{thm}
	Problem \eqref{eq:longest_shortest_partition} has solutions in the class of convex sets, i.e. given $\bo c = (c_i)_{i=1}^n\in \Bbb{R}^n, \sum_{i=1}^n c_i=1$ there exist  convex sets $\Omega^*$ which maximize $PI(\Omega,\bo c)$ among convex sets with fixed volume $|\Omega|=v_d$.
	\label{thm:existence_partitions}
\end{thm}

\emph{Proof:} As in the proof of Theorem \ref{thm:existence_fence} it is straightforward to give upper bounds for $PI(\Omega, \bo c)$ in terms of $w(\Omega)$ (the minimal $\mathcal H^{d-1}$ measure of the projection on a hyperplane). A maximizing sequence $(\Omega_h)$ would have a positive lower bound $0<p_0 \leq w(\Omega_h)$ for the sequence of minimal projections on hyperplanes. Therefore the diameters of $(\Omega_h)$ are bounded from above and we may assume that the convex sets $\Omega_h$ converge to a convex set $\Omega$ (with respect to the Hausdorff distance). The set $\Omega$ also verifies the volume constraint $|\Omega| = v_d$.

It only remains to prove the continuity of the minimal partition perimeters $PI(\Omega_h,\bo c)$ for the convergence with respect to the Hausdorff distance. In order to do this, the same tools as in the proof of Theorem 4.1 in \cite{isoperimetric_convex} can be used. 

\begin{enumerate}[noitemsep,topsep=0pt]
	\item {\bf Lower-semicontinuity.} Theorem 3.4 in \cite{isoperimetric_convex} shows that there exist bilipschitz maps $f_h : \Omega_h \to \Omega$ with Lipschitz constants $\Lip(f_h)$ converging to $1$, the Lipschitz constants of the inverse maps $\Lip(f_h^{-1})$ also converging to $1$.  The volumes and perimeters of the images of finite perimeter sets $E_h \subset \Omega_h$ have upper and lower bounds as follows:
	\begin{align*}
	 \frac{1}{\Lip(f_h^{-1})^d} |E_h| & \leq |f_h(E_h)| \leq \Lip(f_h)^d |E_h| \\ 
	 \frac{1}{\Lip(f_h^{-1})^{d-1}} \Per_{\Omega_h}(E_h) & \leq \Per_\Omega(f_h(E_h))  \leq \Lip(f_h)^{d-1} \Per_{\Omega_h}(E_h)
	\end{align*}
	
	Let $(\omega_h^i)_{i=1}^n$ be a minimal perimeter partition for $\Omega_h$ with constraint $\bo c \in \Bbb{R}^n$. Then $(f_h(\omega_h^i))$ is a partition of $\Omega$ with $\lim_{h \to \infty} |f_h(\omega_h^i)| = c_i|\Omega|$. Extracting a diagonal sequence, we may assume that $(\omega_h^i)_{i=1}^n$ converges  with respect to the Hausdorff distance to a partition $(\omega_i)_{i=1}^n$ of $\Omega$ as $h\to \infty$. Using the estimates above and the fact that the perimeter is lower semi-continuous with respect to the convergence of finite perimeter sets we have
	\begin{align*}
	 PI(\Omega,\bo c) & \leq \sum_{i=1}^n \Per_\Omega(\omega_i) \leq \liminf_{h\to \infty} \sum_{i=1}^n \Per_\Omega(f_i(\omega_h^i)) \\
	 & = \liminf_{h\to \infty} \sum_{i=1}^n \Per_{\Omega_h} (\omega_h^i) = \liminf_{h\to \infty} PI(\Omega_h,\bo c).
	 \end{align*}
	 \item {\bf Upper-semicontinuity.} It remains to prove that $PI(\Omega,\bo c) \geq \limsup_{h\to \infty} PI(\Omega_h,\bo c)$. Reasoning by contradiction, suppose that $PI(\Omega,\bo c)< \limsup_{h\to \infty} PI(\Omega_h,\bo c)$. Up to a subsequence we may assume that $PI(\Omega_h,\bo c)$ converges. Choose $(\omega_i)_{i=1}^n$ a minimal partition in $\Omega$ with constraints $|\omega_i| = c_i|\Omega|$. As in \cite{isoperimetric_convex} using these sets it is possible to construct better competitors on some $\Omega_h$ for large $h$ than the corresponding optimal partition. This leads to a contradiction. 
	 
	 Indeed, $(f_h^{-1}(\omega_h^i))_{i=1}^n$ forms a partition of $\Omega_h$, which may fail to satisfy the volume constraints. Optimality conditions imply that common boundaries of the sets in the partition are regular hypersurfaces. Therefore, it is possible to perturb these boundaries around regular points in order to attain the desired volume constraints. Moreover, for $h$ large enough this will produce partitions which verify 
	 \[ \sum_{i=1}^n \Per_{\Omega_h}(f_h^{-1}(\omega_h)) < PI(\Omega_h,\bo c),\]
	 contradicting the optimality of $PI(\Omega_h,\bo c)$.
\end{enumerate}

This concludes the proof of the existence of solutions for the given problem. \hfill $\square$

\begin{rem}
	Existence results obtained in this section may also be generalized to the case of manifolds, in particular when $\Omega$ is the boundary of a convex set in $\Bbb{R}^d$. There exist sets $\Omega$ which are surfaces of co-dimension $1$ that are boundaries of some convex set in $\Bbb{R}^d$ and have fixed $\mathcal H^{d-1}$ measure which maximize the minimal relative geodesic perimeter of a subset or partition with given $\mathcal H^{d-1}$ measure constraints.
\end{rem}

\subsection{Relaxation of the perimeter - Gamma convergence}
\label{sec:theory-relax}

A key point in our approach is to approximate minimal length partitions $SP(\Omega,\bo c)$. In order to avoid difficulties related to the treatment of the partition constraint it is convenient to represent each set in the partition $\omega_i$ as a density $u_i : \Omega \to [0,1]$. Then, the partition constraint can be simply expressed by the algebraic equality $\sum_{i=1}^n u_i = 1$ on $\Omega$. The next aspect is to approximate the perimeter of a set represented via its density function. A well known technique is to use a $\Gamma$-convergence relaxation for the perimeter inspired by a result of Modica and Mortola \cite{modica-mortola}. The main idea is to replace the perimeter with a functional that, when minimized, yields minimizers converging to those that minimize the perimeter. 

Let us briefly recall the concept of $\Gamma$-convergence and the property that motivates its use when dealing with numerical optimization. 

\begin{rem}
	Let $X$ be a metric space. For $\varepsilon>0$ consider the functionals $F_\varepsilon, F: X \to [0,+\infty]$. We say that $F_\varepsilon$ $\Gamma$-converges to $F$ and we denote $F_\varepsilon \gconv F$ if the following two properties hold:
	
	\begin{itemize}
		\item[(LI)] For every $x \in X$ and every $(x_\varepsilon) \subset X$ with $(x_\varepsilon) \to x$ we have
		\begin{equation} F(x) \leq \liminf_{\varepsilon \to 0} F_\varepsilon(x_\varepsilon)\label{LI} \end{equation}
		
		\item[(LS)] For every $x \in X$ there exists $(x_\varepsilon) \subset X$ such that $(x_\varepsilon) \to x$ and
		\begin{equation} F(x) \geq \limsup_{\varepsilon \to 0} F_\varepsilon(x_\varepsilon).\label{LS}\end{equation}
	\end{itemize}
	\label{def-g-conv}
\end{rem}
An important consequence is the following classical result concerning the convergence of minimizers of a sequence of functionals that $\Gamma$ converge.
\begin{prop}
	Suppose that $F_\varepsilon \gconv F$ and $x_\varepsilon$ minimizes $F_\varepsilon$ on $X$. Then every limit point of $(x_\varepsilon)$ is a minimizer for $F$ on $X$.
	\label{prop:gconv}
\end{prop}
Therefore, in practice, in order to approximate the minimizers of $F$ it is possible to search for minimizers of $F_\varepsilon$, for $\varepsilon$ small enough. 

Let us now state the two theoretical results that are used in this work concerning the $\Gamma$-convergence relaxation of the perimeter and of the total perimeter of a partition, with integral constraints on the densities. The first result is the classical Modica-Mortola theorem \cite{modica-mortola}. Various proofs can be found in \cite{gammaconvalberti,braides2,buttazzogconv}. In the following $\Omega$ is a bounded, Lipschitz open set. Consider a double well potential $W: \Bbb{R} \to [0,\infty)$ which verifies the following assumptions: $W$ is of class $C^1$, $W(z)=0$ if and only if $z \in \{0,1\}$ and $W$ has exactly three critical points. For such a double well potential $W$ described previously, denote $\gamma =2\int_0^1 \sqrt{W(s)}ds$. In the following $c\in [0,1]$ represents the fraction used for the volume constraint.

\begin{thm}[Modica-Mortola] Define $F_\varepsilon,F : L^1(\Omega) \to [0,+\infty]$ by
	\[ F_\varepsilon(u) = \begin{cases}\displaystyle\int_\Omega \left( \varepsilon |\nabla u|^2 +\frac{1}{\varepsilon} W(u) \right) & u \in H^1(\Omega), \int_\Omega u = c|\Omega| \\ +\infty & \text{otherwise} \end{cases} \]
	and
	\[ F(u)             = \begin{cases} \gamma\Per_\Omega(\{u=1\}) & u \in BV(\Omega; \{0,1\}), \int_\Omega u = c|\Omega| \\ +\infty & \text{otherwise}\end{cases}.\]
	Then $F_\varepsilon \stackrel{\Gamma}{\longrightarrow} F$
	in the $L^1(\Omega)$ topology.
	\label{thm:modica-mortola}
\end{thm}

In \cite{oudet} this result was generalized to the case of partitions and was used to compute approximations for $SP(\Omega,\bo c)$. For $\bo c \in \Bbb{R}^n$ with $\sum_{i=1}^n c_i=1$, in order to simplify notations, let us denote by $X(\Omega,\bo c)$ the space of admissible densities which verify the integral constraints and the algebraic non-overlapping constraint
\[ X(\Omega, \bo c) = \{ \bo u = (u_i)_{i=1}^n \in L^1(\Omega)^n : \int_\Omega u_i = c_i|\Omega|, \sum_{i=1}^n u_i = 1\}.\]
The $\Gamma$-convergence result in the case of partitions is recalled below.

\begin{thm}
	Define $G_\varepsilon , G: L^1(\Omega) \to [0,+\infty]$ by
	\[ G_\varepsilon (\bo u) = \begin{cases}
	\ds \sum_{i=1}^n\int_\Omega \left( \varepsilon |\nabla u_i|^2 +\frac{1}{\varepsilon} W(u_i) \right) & \text{ if } \bo u \in (H^1(\Omega))^n \cap X(\Omega,\bo c)\\
	+\infty & \text{ otherwise}
	\end{cases}
	\]
	\[ G(\bo u) = \begin{cases}
	\gamma\sum_{i=1}^n\Per_\Omega(\{u_i = 1\}) & \text{ if } \bo u \in (BV(\Omega,\{0,1\}))^n  \cap X(\Omega,\bo c)\\
	+\infty & \text{ otherwise}
	\end{cases}
	\]
	Then $G_\varepsilon \gconv G$ in the $(L^1(\Omega))^n$ topology.
	\label{thm:gconv-part}
\end{thm}
A proof of this result can be found in \cite{oudet}. In the numerical simulations the double well potential is chosen to be $W(s) = s^2(1-s)^2$ which gives the factor $\gamma = 1/3$ in the results shown above.
 
\begin{rem}
	It can be seen that $SF(\Omega,c)$ corresponds to a density that is a minimizer of $F$ in Theorem \ref{thm:modica-mortola}. Moreover, $SP(\Omega,\bo c)$ corresponds to a family of densities which minimizes $G$ in Theorem \ref{thm:gconv-part}. Using the result recalled in Proposition \ref{prop:gconv} it is possible to approximate these minimizers by those of $F_\varepsilon$ and $G_\varepsilon  $, respectively, for $\varepsilon$ small enough. From a numerical point of view, dealing with the minimization of $F_\varepsilon$ and $G_\varepsilon  $ is easier since the variable densities are $H^1$ regular. 
\end{rem}

\begin{rem}
	The structures of minimizers of $F_\varepsilon$ was widely studied in the literature as can be seen in the papers \cite{gurtin-matano}, \cite{luckhaus-modica}, \cite{sternberg}.
	It can immediately be seen that, assuming $W$ is at least of class $C^1$, minimizers $u$ of $F_\varepsilon$ verify an optimality condition of the form
	\begin{equation} \int_\Omega \left( 2\varepsilon \nabla u \cdot \nabla \varphi +\frac{1}{\varepsilon} W'(u)\varphi + \mu \varphi\right) = 0, \text{ for every }\varphi \in H^1(\Omega)
	\label{eq:opt-cond-MM}
    \end{equation}
	where $\mu \in \Bbb{R}$ is a Lagrange multiplier for the volume constraint. Classical regularity theory results that can be found in \cite{gilbarg-trudinger} allow us to employ a bootstrap argument and conclude that $u$ is of class $C^\infty$ in the interior of $\Omega$ and $u$ has the regularity of $\Omega$ up to the boundary. For example, for smooth domains $\Omega$ the optimizer $u$ is also smooth up to the boundary of $\Omega$. Moreover, it can be proved that the minimizer $u$ takes values in $[0,1]$. In the case when $\Omega$ is convex results found in \cite[Chapter 3]{grisvard} show that solutions of the above problem are in $H^2(\Omega)$.
	
	The same type of results hold for minimizers of $G_\varepsilon$ in the case of partitions, with eventual singularities at junction points between three or more phases in the partition. Nevertheless, the contact between the optimal partition and the boundary $\partial \Omega$ has the desired regularity.
	\label{rem:smoothness-MM}
\end{rem}

\begin{rem}
	The results in \cite{luckhaus-modica} show that the Lagrange multiplier $\mu$ for the volume constraint in \eqref{eq:opt-cond-MM} has a geometric interpretation. Given a volume fraction $c \in (0,1)$, 
	as $\varepsilon\to 0$, the Lagrange multiplicator $\mu$ converges to $\gamma$ times the mean curvature of the shortest fence set $SF(\Omega,c)$, where $\gamma$ was defined before Theorem \ref{thm:modica-mortola}. Recall that this minimal set $SF(\Omega,c)$, being optimal for the relative perimeter under a volume constraint has constant mean curvature inside $\Omega$. Moreover, as shown in \cite{sternberg}, taking $\varphi=1$ as a test function in \eqref{eq:opt-cond-MM} gives an explicit formula for the Lagrange multiplier
	\begin{equation}
	 \mu = -\frac{1}{\varepsilon|\Omega|}\int_\Omega W'(u).
	 \label{eq:formula-lagrange-mult}
	\end{equation}
	\label{rem:lag-mult}
\end{rem}

Since in the numerical section we deal with the minimization of $F_\varepsilon, G_\varepsilon  $ for fixed $\Omega$, we briefly recall existence results related to these problems. In the following we suppose that the double well potential $W$ is Lipschitz continuous on $\Bbb{R}$. This is not restrictive since minimizers of $F_\varepsilon, G_\varepsilon $ are densities which take values in $[0,1]$, which means that values of $W$ far away from this interval do not matter in the analysis. 

\begin{thm}
 (i)	Problems 
	\[ \min_{ u \in L^1(\Omega)} F_\varepsilon(u) \text{  and  } \min_{\bo u \in L^1(\Omega)^n} G_\varepsilon (\bo u) \]
	admit solutions for $\Omega$ a Lipschitz domain with finite volume. In the following we denote by $I_\varepsilon(\Omega,c)$ and $PI_\varepsilon (\Omega, \bo c)$ the optimal values obtained when minimizing $F_\varepsilon$ and $G_\varepsilon$, respectively.	
	
 (ii)   Given $c \in (0,1)$ and $\bo c = (c_i) \in \Bbb{R}^n, \sum_{i=1}^n c_i=1$ problems 
 \[ \max_{|\Omega|=v_d} I_\varepsilon (\Omega,c) \text{ and } \max_{|\Omega|=v_d} PI_\varepsilon (\Omega, \bo c)\]
 admit solutions in the class of convex sets.	
 \label{thm:existence-density}
\end{thm}

\emph{Proof:} The proof of (i) is classical. Note that the constraints on the density functions are embedded in the definition of the functionals $F_\varepsilon, G_\varepsilon $ to be minimized. We give the ideas for $G_\varepsilon $ as $F_\varepsilon$ is just a particular case. The existence proof goes as follows:
\begin{itemize}[noitemsep,topsep=0pt]
	\item The functional $ G_\varepsilon $ is obviously bounded from below by zero. Moreover, truncating the density functions $(u_i)$ to take values in $[0,1]$ does not increase the value of $G_\varepsilon$. This allows us from now on to assume that the densities have values in this interval.
	\item Minimizing sequences exist and they are bounded in $H^1(\Omega)^n$, which allows us to extract a subsequence weakly converging in $H^1$. The constraints are stable under the $L^2$ convergence. Moreover, the lower-semicontinuity of the $H^1$ norm and Fatou's lemma allow us to see that any weak $H^1$-limit point of the minimizing sequence is a minimizer. 
\end{itemize}

The proof of (ii) follows the same lines as the proofs of Theorems \ref{thm:existence_fence}, \ref{thm:existence_partitions}. As in the proof of these theorems, we start by noticing that the minimal $\mathcal H^{d-1}$ measure $w(\Omega)$ of the projection of $\Omega$ on a hyperplane is bounded from above. We detail the proof for $I_\varepsilon $, while the proof in the case of partitions follows the same path. In order to emphasize the dependence of $F_\varepsilon$ on $\Omega$ we write $F_\varepsilon(u) = F_\varepsilon(\Omega,u)$.

{\bf Upper bound.} Choose $e_0$ the direction for which $\mathcal H^{d-1}(\Omega| e_0^\perp)$ is minimal and equal to $w(\Omega)$. Given a hyperplane $\zeta$ orthogonal to $e_0$ consider the function $u_\varepsilon = \varphi_\varepsilon(d(x))$, where $d(x)$ is the signed distance to the hyperplane $\zeta$ (choosing an orientation) and $\psi_\varepsilon, \varphi_\varepsilon$ are given by
\[ \psi_\varepsilon(t) = \int_0^t \frac{\varepsilon}{\sqrt{\varepsilon+W(s)}}ds,\ \ \varphi_\varepsilon(t) = \begin{cases}
0 & t\leq 0\\
\psi_\varepsilon^{-1}(t) & 0 \leq t \leq \psi_\varepsilon(1) \\
1 & t \geq \psi_\varepsilon(1).
\end{cases}\]
This type of construction is standard when proving the limsup part of the $\Gamma$-convergence proof for the Modica-Mortola type results in Theorems \ref{thm:modica-mortola}, \ref{thm:gconv-part} (see for example \cite{modica_gradient_theory}). The coarea formula and the fact that $|\nabla d(x)| = 1$ allows us to write
\[ \int_\Omega u_\varepsilon = \int_{\Bbb{R}} \int_{\{d(x)=t\}} \varphi_\varepsilon(t)d\mathcal H^{d-1}dt.\] 
The definition of $\varphi_\varepsilon$ and a continuity argument allow us to deduce that there is a position of the hyperplane for which the constraint $\int_\Omega u_\varepsilon = c|\Omega|$ is verified. 

Using the coarea formula to evaluate $F_\varepsilon(\Omega,u_\varepsilon)$ we obtain
\begin{align*}
F_\varepsilon(\Omega,u_\varepsilon) & = \int_\Omega \left( \varepsilon |\varphi_\varepsilon'(d(x))|^2 +\frac{1}{\varepsilon} W(\varphi_\varepsilon(d(x))) \right)\\
& = \int_{\Bbb{R}} \int_{\{d(x)=t\}}\left( \varepsilon|\varphi_\varepsilon'(t)|^2+\frac{1}{\varepsilon}W(\varphi_\varepsilon(t))\right)d\mathcal H^{d-1}dt\\
& \leq  w(\Omega) \int_{0}^{\psi_\varepsilon(1)} \left( \varepsilon|\varphi_\varepsilon'(t)|^2+\frac{1}{\varepsilon}W(\varphi_\varepsilon(t))\right) dt.
\end{align*}
The last inequality comes from the fact that $\{d(x)=t\}$ is a slice of $\Omega$ orthogonal to $e_0$ and its $\mathcal H^{d-1}$ measure is at most $w(\Omega)$. Moreover, we can restrict the bounds in the one dimensional integral to $0$ and $\psi_\varepsilon(1)$ since for $t$ not in this interval the integrand is zero. A simple computation gives
\[ \varphi_\varepsilon'(t) = \frac{1}{\psi_\varepsilon'(\psi_\varepsilon^{-1}(t))}=\frac{1}{\varepsilon} \sqrt{\varepsilon+W(\varphi_\varepsilon)}.\]
Thus we obtain
\[ I_\varepsilon(\Omega,c) \leq F_\varepsilon(\Omega,u_\varepsilon)\leq \frac{2w(\Omega)}{\varepsilon} \int_0^{\psi_\varepsilon(1)} (\varepsilon+W(\varphi_\varepsilon))dt = 2w(\Omega) \int_0^1 \sqrt{\varepsilon+W(s)}ds,\]
where the last equality comes from the change of variables $s= \varphi_\varepsilon(t)$. This quantity depends only on $W$ and $w(\Omega)$ and  is bounded from above independently of $\Omega$.

{\bf Compactness.} The same argument used in the proof of Theorem \ref{thm:existence_fence} can be applied in order to conclude that there exists a maximizing sequence $(\Omega_h)$ converging in the Hausdorff distance to a convex set $\Omega$ with non-empty interior and volume $|\Omega|=v_d$. Moreover, we may assume that there exists a bounded open set $D$ such that $(\Omega_h)_{h\geq 1}, \Omega$ are contained in $D$.

Following the ideas in \cite[Chapter 2]{henrot-pierre-english} we may assume that $(\Omega_h)$ and $\Omega$ satisfy an $\varepsilon$-cone condition, or equivalently that they are Lipschitz regular with a uniform Lipschitz constant. In this case, the convergence with respect to the Hausdorff distance implies that $|\Omega_h\setminus \Omega| + |\Omega\setminus \Omega_h| \to 0$.

{\bf Continuity.} It now remains to prove that $I_\varepsilon (\Omega_h,c) \to I_\varepsilon (\Omega,c)$ as $h \to \infty$. Let us note first that since $\Omega_h$ is a maximizing sequence we have $I_\varepsilon (\Omega,c) \leq \lim_{h \to \infty} I_\varepsilon (\Omega_h,c)$. 
Consider a minimizer $u \in H^1(\Omega)$ such that $F_\varepsilon(\Omega,u) = I_\varepsilon (\Omega,c)$. 

Since $(\Omega_h)$ and $\Omega$ have a uniform Lipschitz constant $L$ (as recalled above), using the extension theorems recalled in \cite[Theorem 3.4] {extension_theorems}, there exists an extension $\tilde u \in W^{1,p}(D)$ of $u$ which verifies $\|\tilde u\|_{W^{1,p}(D)} \leq \text{Const}(L) \|u\|_{W^{1,p}(\Omega_h)}$. Together with the results recalled in Remark \ref{rem:smoothness-MM} we find that $\varepsilon |\nabla \tilde u|^2+\frac{1}{\varepsilon} W(\tilde u) \in L^\infty(D)$. Combining this with the fact that $|\Omega\setminus \Omega_h|+|\Omega_h\setminus \Omega|\to 0$ implies that
\[ F_\varepsilon(\Omega_h,\tilde u) \to F_\varepsilon(\Omega, u).\]
We cannot conclude yet, since $\tilde u$ may not satisfy the integral constraints on $\Omega_h$. 

 In order to fix this, let $x_0$ be a point in the interior of $\Omega$. For $h$ large enough there exists a ball $B_\delta$ of radius $\delta>0$ such that $B_\delta \subset \Omega \cap \Omega_h$. Denote by $d_\delta$ the function which is equal to the distance to $\partial B_\delta$ inside $B_\delta$ and zero outside. We use this function to construct functions $u_h = \tilde u + x_h d_\delta$, for $x_h \in \Bbb{R}$, which verify the integral constraints $\int_{\Omega_h} u_h = c|\Omega_h|$. Since
\[ \left | \int_\Omega \tilde u -\int_{\Omega_h} u\right| \leq \left|\int_{\Omega\setminus \Omega_h} \tilde u\right|+\left|\int_{\Omega_h \setminus \Omega} u\right| = O(|\Omega\setminus \Omega_h|+|\Omega_h\setminus \Omega|) \to 0, \]
we necessarily have $x_h\to 0$. This immediately shows that 
\[ |F_\varepsilon(\Omega_h,u_h) - F_\varepsilon(\Omega_h,\tilde u)| \to 0\]
as $h \to \infty$.

Since $I_\varepsilon (\Omega_h,c) \leq F_\varepsilon(\Omega,u_h)$ we find that $\limsup_{h\to \infty} I_\varepsilon (\Omega_h,c) \leq F_\varepsilon(\Omega,u) = I_\varepsilon (\Omega,c)$. This concludes the proof of the existence result. The case of partitions can be handled in a similar manner with the additional difficulty that the area constraints and sum constraints need to be handled simultaneously. This can be achieved by modifying the candidate densities in a finite family of balls. \hfill $\square$

\section{Numerical modeling}

\subsection{Numerical framework for approximating minimal perimeter partitions}
\label{sec:numerics-partitions}

In this section the numerical minimization of $F_\varepsilon$ and $ G_\varepsilon $ is discussed. Since $\Omega$ is a general domain in this work, we choose to work with finite element discretizations. Given $\mathcal T_h$ a triangulation of $\Omega$, denote by $(x_j)_{j=1}^N$ the set of the nodes.  Working with $\Bbb P_1$ Lagrange finite elements, a piecewise affine function $u$ defined on the mesh $\mathcal T_h$ is written $\sum_{j=1}^N u_j \phi_j$. As usual, $\phi_j$ are the piece-wise linear functions on each triangle, characterized by $\phi_j(x_k) = \delta_{jk}$. For a $\Bbb P_1$ finite element function, the values $u_j$ are given by $u(x_j)$ and we denote $\bo u = (u_j) = (u(x_j))\in \Bbb{R}^N$. With these notations, it is classical to introduce the mass matrix $M$ and the rigidity matrix $K$ defined by
\[ M = \left( \int_{\mathcal T_h} \phi_i \phi_j\right)_{1\leq i,j \leq N} \text{  and   } K = \left( \int_{\mathcal T_h} \nabla \phi_i \cdot \nabla \phi_j\right)_{1\leq i,j \leq N}\]
As an immediate consequence of the linearity of the decompositions $u = \sum_{j=1}^N u_j \phi_j$, $v = \sum_{j=1}^N v_j \phi_j$ we have that 
\[ \int_{\mathcal T_h} u v = \bo u^T M \bo v \text{ and } \int_{\mathcal T_h} \nabla u \cdot \nabla v = \bo u^T K \bo v.\]
This immediately shows that the functionals $F_\varepsilon$ and $ G_\varepsilon $ can be expressed in terms of the mass and rigidity matrices $M$ and $K$ using the expression
\begin{equation}
  \int_{\mathcal T_h} \left(\varepsilon |\nabla u|^2 + \frac{1}{\varepsilon} u^2(1-u)^2\right) = \varepsilon \bo u^T K \bo u + \frac{1}{\varepsilon} \bo v^T M \bo v =: \mathcal F(\bo u)
 \label{eq:discrete-MM}
\end{equation}
where $\bo v = (u_j(1-u_j))_{j=1}^N$. The gradient of this expression w.r.t. $\bo u$ can be computed and is given by
\begin{equation}
 \nabla \mathcal F(\bo u) = 2\varepsilon K \bo u + \frac{2}{\varepsilon} M \bo v \odot (1-2\bo u),
 \label{eq:grad-discrete-MM}
\end{equation}
where $\odot$ denotes pointwise multiplication of two vectors: $\bo u \odot \bo v = (u_jv_j)_{j=1}^N$. 

It is obvious that with \eqref{eq:discrete-MM} and \eqref{eq:grad-discrete-MM} it is possible to implement a gradient-based optimization algorithm in order to minimize $F_\varepsilon$ and $G_\varepsilon$. The software FreeFEM \cite{freefem} is used for constructing the finite element framework and the algorithm LBFGS from the package Nlopt \cite{nlopt} is used for the minimization of \eqref{eq:discrete-MM}. We address the question of handling the constraints in the next section.

\subsection{Area constraints and projections}
\label{sec:constraints}

The area or volume constraint can be expressed with the aid of the vector $\bo m =  M \bo e$ with $\bo e = (1,1,...,1) \in \Bbb{R}^N$. Indeed, with this notation, for a finite element function $u$ we have $\int_{\mathcal T_h} u = \bo m\cdot \bo u$. 

{\bf Projection for one phase.} Let us start with the projection of one function onto the integral constraint. Given a $\Bbb{P}_1$ finite element function $u$ and its values $\bo u$ at the nodes we search a function $u_0 $ with values at nodes $\bo u_0= \bo u + \alpha \bo m$ verifying the constraint $\bo m\cdot  \bo u_0 = c$ by solving
\[ (\bo u + \alpha \bo m) \cdot \bo m = c,\]
which leads to $\alpha = (c- \bo u \cdot \bo m) / (\bo m \cdot \bo m)$.

An alternative way of handling the constraint during the optimization process is to project the initial vector on the constraint and project the gradient onto the hyperplane $\bo x \cdot \bo m = 0$  at each iteration. This can simply be done by using $c=0$ in the relation above. Such a modification of the gradient allows us to use efficient \emph{black-box} optimization toolboxes, since quasi-Newton algorithms like LBFGS will perform updates based on a number of gradients stored in memory. If these gradients verify $\bo x \cdot \bo m=0$, the integral constraint will be preserved throughout the optimization process. 

{\bf Projection for multiple phases.} In the case of partitions projections on the integral constraints were already proposed in \cite{oudet} (when using finite differences) and in \cite{bogosel_oudet_manifolds} (when using finite elements). A drawback of using orthogonal projections parallel to the vector $\bo m$ is the fact that the vector $\bo u$ is modified almost everywhere in the domain $\Omega$, which also includes the regions where it is $0$ or $1$. As observed in \cite{many_cells}, this can cause resulting optimal densities to be non-zero at interfaces between two cells and at triple points. The solution proposed was to use instead projections parallel to $\sqrt{2W(u_i)}$. It is possible to note that, using this method, the functions $u_i$ are mainly modified only at the interface of transition between $0$ or $1$. 

Let us now describe the construction of the projection algorithm on the constraints
\[ \int_\Omega u_i = c_i|\Omega|,\ \ \sum_{i=1}^n u_i = 1.\]
with the compatibility condition $\sum c_i =1$. Consider  $\lambda \in H^1(\Omega)$ and $(\mu_i) \in \Bbb{R}^n$ and perform the transformation
\[ u_i+\lambda \sqrt{2W(u_i)} + \mu_i \sqrt{2W(u_i)}\]
in order to satisfy the constraints
\begin{equation}
 \int_\Omega u_i + \int_\Omega\lambda \sqrt{2W(u_i)}+ \mu_i \int_\Omega\sqrt{2W(u_i)}=c_i|\Omega|, i=1,...,n
 \label{eq:area-constr-proj}
\end{equation}
and 
\begin{equation}
   \sum_{i=1}^n u_i +\lambda \sum_{i=1}^n \sqrt{2W(u_i)}  + \sum_{i=1}^n \mu_i  \sqrt{2W(u_i)} = 1
   \label{eq:sum-constr-proj}
\end{equation}
It is easy to note that:
\begin{itemize}[noitemsep,topsep=0pt]
	\item in view of \eqref{eq:sum-constr-proj}, given $\mu_i$ we can find $\lambda$:
	\[ \lambda = \frac{1-\sum_{i=1}^n u_i -\sum_{i=1}^n \mu_i \sqrt{2W(u_i)}}{\sum_{i=1}^n \sqrt{2W(u_i)}}.\]
	\item in view of \eqref{eq:area-constr-proj}, given $\lambda$ we can find $\mu_i$ using the relations above.
\end{itemize}
In the following we introduce the quantities $\bar \lambda_i = \int_\Omega \lambda \sqrt{2W(u_i)}$. Again, in view of \eqref{eq:area-constr-proj}, if $\bar \lambda_i$ are known, then $\mu_i$ are known and so is $\lambda$. In order to obtain a system for $\bar \lambda_i$, let us note that
\[ \mu_i = \frac{c_i|\Omega|-\int_D u_i - \bar \lambda_i}{\int_D \sqrt{2W(u_i)}}.\]
With this in mind we get
\begin{align*}
\bar\lambda_i & = \int_\Omega\lambda \sqrt{2W(u_i)}  = \int_\Omega \left(\frac{1-\sum_{j=1}^n u_j -\sum_{j=1}^n \mu_j \sqrt{2W(u_j)}}{\sum_{j=1}^n \sqrt{2W(u_j)}} \right)\sqrt{2W(u_i)}\\
&=\int_\Omega \left(\frac{1-\sum_{j=1}^n u_j -\sum_{j=1}^n  \left(\frac{c_j|\Omega|-\int_\Omega  u_j - \bar \lambda_j}{\int_D \sqrt{2W(u_j)}}\right) \sqrt{2W(u_j)}}{\sum_{j=1}^n \sqrt{2W(u_j)}} \right)\sqrt{2W(u_i)}
\end{align*} 
In order to further simplify the above expression, let's make the following notations
\begin{itemize}[noitemsep,topsep=0pt]
	\item $\sqrt{2W(u_i)} = w_i$
	\item $1-\sum_{j=1}^n u_i = E$
	\item $c_j|\Omega|-\int_D u_j = F_j$, $j=1,...,n$
\end{itemize}
This gives
\begin{align*}
\bar\lambda_i & =\int_\Omega \left(\frac{E -\sum_{j=1}^n  \left(\frac{F_j - \bar \lambda_j}{\int_\Omega w_j}\right) w_j}{\sum_{j=1}^n w_j} \right)w_i\\
& = \int_\Omega \frac{E-\sum_{j=1}^n F_jw_j /\int_\Omega w_j }{\sum_{j=1}^n w_j} w_i+ \sum_{j=1}^n \int_\Omega  \frac{w_i w_j/\int_\Omega w_j}{\sum_{j=1}^n w_j} \bar \lambda_j
\end{align*} 
which can be written in the form $(I-A)\bar \lambda = b$
with $\bar \lambda = (\bar \lambda_1,...,\bar \lambda_n)$ and
\[ A = \left(\int_\Omega  \frac{w_i w_j/\int_\Omega w_j}{\sum_{j=1}^n w_j}\right)_{i,j=1,...,n}, \ \ b =\left(\int_\Omega \frac{E-\sum_{j=1}^n F_jw_j /\int_\Omega w_j }{\sum_{j=1}^n w_j} w_i\right)_{j=1,...,n}.\]
One may note that the above system $(I-A)\bar \lambda = b$ is singular since the sum on the columns of $A$ is equal to $1$ and therefore
\[ (I-A^T)\bo e = 0,\]
where $\bo e = (1,...,1) \in \Bbb{R}^n$. This is due to the fact that one of the constraints is redundant, in view of the compatibility condition. In practice we simply discard one unknown and set it to zero.

As noted previously, the same procedure can be applied to the gradients $g_i$ associated to each $u_i$ in order to satisfy at every iteration
\[ \int_\Omega g_i = 0,\ \ \sum_{i=1}^n g_i = 0.\]
This allows us to preserve the constraints when using a black-box LBFGS optimizer when initial parameters satisfy the integral and sum constraints.

\subsection{Initializations for 2D partitions - Voronoi diagrams}
\label{sec:voronoi}

The optimization algorithm for approximating $I_\varepsilon (\Omega,c)$ and $PI_\varepsilon (\Omega,\bo c)$ is ready to be implemented, following the ideas shown in the previous sections. There is, however, the choice of the initialization which is non-trivial and which has an impact on the performance of the optimization algorithm. It was already noted in \cite{oudet} and \cite{bogosel_oudet_manifolds} that starting from random initializations is possible, but some additional work needs to be done in order to avoid constant phases, which are encountered at some local minimizers. Keeping in mind that the optimal partition problem needs to be solved multiple times during the optimization algorithm, we propose below a different initialization strategy, based on Voronoi diagrams. The use of Voronoi diagrams for generating initializations is a rather natural idea when dealing with partitions and was already mentioned in \cite{cox-partitions}. In this section $\Omega$ is assumed to be a polygon.

Using random Voronoi diagrams is not very helpful, since area constraints are not verified in general. This led us to consider Voronoi diagrams which verify the area constraints, which in the literature are called \emph{capacity-constrained} Voronoi diagrams. Algorithms for computing such diagrams were proposed in \cite{capacity-constrained-discrete} for the discrete case and in \cite{capacity-constrained-continuous} for the continuous case. In the continuous case the method employed in \cite{capacity-constrained-continuous} was to optimize the position of one Voronoi point at a time using the gradient-free Nelder-Mead method. In \cite{xin2016} the authors propose efficient ways of generating such diagrams, but for weighted Voronoi diagrams only. In the following we propose an alternative method for constructing capacity-constrained Voronoi diagrams by computing the sensitivity of the areas of the Voronoi cells with respect to the position of the points generating the respective Voronoi diagram. Since we are also interested in minimizing the perimeter, the computation of the sensitivity of the perimeter of Voronoi cells is also described. 

{\bf Terminology related to Voronoi diagrams.} Given a set of points $p_1,...,p_n \in \Bbb{R}^2$ (called Voronoi points) the associated Voronoi diagram consists of $n$ Voronoi cells $V_1,...,V_n$ defined for $i=1,...,n$ by
\[ V_i = \{ x\in \Bbb{R}^2 : |x-p_i|\leq |x-p_j|, j = 1,...,n, j \neq i\}.\]
Each Voronoi cell $V_i$ is a polygonal region (possibly unbounded). The vertices of $V_i$ are simply called \emph{vertices} in the following (please observe the difference between the \emph{Voronoi points} and the \emph{Voronoi vertices}). The edges of $V_i$ are called \emph{ridges}, some of which can be unbounded. Each ridge connects two Voronoi vertices (possibly at infinity, for unbounded ridges) called \emph{ridge vertices}. Moreover, each ridge separates two of the initial points, called \emph{ridge points}. All structure information of a Voronoi diagram associated to a set of points can be recovered as an output to some freely available software like \texttt{scipy.spatial.Voronoi}. The Voronoi diagrams are not restricted to a bounded domain. It is possible, however, to consider restrictions of a Voronoi diagram to a bounded set $\Omega$ by simply intersectiong the regions $V_i$ with the set $\Omega$. In our implementation the intersection of polygons is handled using the \texttt{Shapely} Python package for computational geometry. 

In the following, given the points $p_i, i=1,...,n$, we consider the Voronoi regions restricted to a finite domain $\Omega$ re-defined by $V_i = V_i \cap \Omega$. Note that in some cases, some $V_i$ may be void if $\Omega$ does not contain the associated point $p_i$. We explain below how to compute the gradients of the areas and perimeters of $V_i$ with respect to positions of the points $p_i$. 

{\bf Gradient of the areas of the Voronoi cells.} The derivative of a functional that can be represented as an integral over the Voronoi cell $V_i$ with respect to the Voronoi points can be computed if the normal displacement of the cell is known. This fact was recalled in \cite{Du1999} and \cite{Iri} and is classical in the shape derivative theory. However, since the functionals considered there were sums over all Voronoi cells $V_i$, the contributions coming from the variations of the boundary cancelled themselves and only the variation of the integrand mattered.

This is no longer the case in our situation. The area of the voronoi cell $V_i$ is $A_i = \int_{V_i} 1 dx$ and its directional derivative when perturbing a point $p_j$ in direction $d$ is given by the integral on the boundary of the normal variation of  $V_i$: $A_i'(d)= \int_{\partial V_i} \theta.n$, where $\theta$ is the infintesimal displacement of the boundary of $V_i$ when moving $p_j$ in the direction $d_i$. More explicitly, if $V_i(t)$ is the Voronoi cell for $p_j+td$ then $\theta = \lim_{t\to 0} \frac{v(t)}{t}$, where the vector field $v(t)$ is defined by $V_i+v(t)=V_i(t)$ on the boundary of $V_i$.

\begin{figure}
	\centering
		\includegraphics[width=0.2\textwidth]{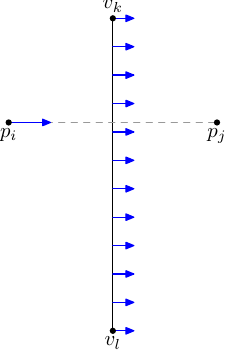}\quad
		\includegraphics[width=0.2\textwidth]{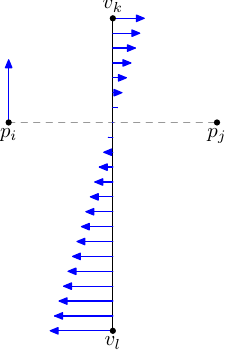}\quad
		\includegraphics[width=0.2\textwidth]{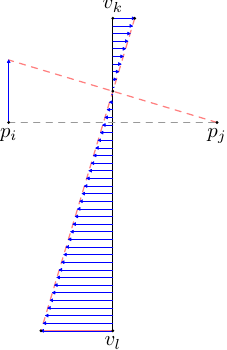}
\caption{Normal perturbation of the Voronoi ridge when moving one of the point in the normal and tangent directions.}
\label{fig:AreaPerturb}
\end{figure}

  For a given ridge $v_kv_l$ with associated ridge points $p_i,p_j$, we perturb the point $p_i\mapsto p_i+\delta$ and investigate the derivative of the normal perturbation of $v_kv_l$ as $\delta \to 0$. The two main perturbations are the following:
\begin{itemize}[noitemsep,topsep=0pt]
	\item {\bf $\delta$ is collinear with $p_ip_j$}: in this case the perturbation induced on the ridge is just $\delta/2$. (see Figure \ref{fig:AreaPerturb} (left)). The associated infinitesimal normal perturbation  is constant equal to $1/2$. 
	\item {\bf $\delta$ is orthogonal to $p_ip_j$}: in this case the infinitesimal perturbation induced on the ridge is a rotation around the intersection $m_{ij}$ of $p_ip_j$ and $v_kv_l$. The associated infinitesimal normal perturbation varies linearly on $v_kv_l$ from $-|v_l-m_{ij}|/|p_i-p_j|$ to $|v_k-m_{ij}|/|p_i-p_j|$ (the signs vary with respect to the orientation of the orthogonal perturbation). (see Figure \ref{fig:AreaPerturb} (middle)). In order to prove this it is enough to consider the normal perturbation $v(t)$ of the ridge $v_kv_l$ illustrated in Figure \ref{fig:AreaPerturb} (right) and take the limit $v(t)/t$ as $t \to 0$.  
\end{itemize}

For a general perturbation $\delta$ of $p_i$ we denote by $\bo n$ the normal vector to $v_kv_l$ pointing outwards to $V_i$ and by $\bo t$ the unit vector collinear with $\overrightarrow{v_lv_k}$. Furthermore, consider the notations for the normal and tangential contributions (computed as one dimensional integrals on $v_kv_l$ of the infinitesimal perturbations described above):
\begin{equation} z_n = (\delta \cdot\bo n)/2 |v_k-v_l|,\  z_t = (\delta \cdot \bo t) \frac{1}{2|p_i-p_j|} (|v_k-m_{ij}|^2-|v_l-m_{ij}|^2).
\label{eq:area-normal-tangential}
\end{equation}
By symmetry, these contributions will be similar, but with changed signs when perturbing $p_j$ with $\delta$.
The contributions to the gradients of the areas of the cells $V_i$ and $V_j$ when perturbing $p_i$ or $p_j$ are described in the table below:
\begin{equation}
	\begin{tabular}{|c|c|c|}
		\hline
		 & $V_i$ & $V_j$  \\ \hline 
		 $p_i$ & $z_n+z_t$ & $-z_n-z_t$ \\ \hline 
		 $p_j$ & $z_n-z_t$ & $-z_n+z_t$  \\ 
		\hline 
	\end{tabular}
	\label{eq:grad-contribution}
\end{equation}
The algorithm for computing the gradient of the areas of the cells simply iterates over all the Voronoi ridges that intersect $\Omega$ and for each ridge adds the contributions described in \eqref{eq:grad-contribution}.

Algorithm \ref{alg:Voronoi-grad-areas} describes the computation of the gradient of the areas of the cells. The coordinates of the $n$ input points are given in the vector $\bo x \in \Bbb{R}^{2n}$ and the output is the real matrix $M$ of size $2n\times n$ containing as columns the gradients of the areas of the $n$ cells with respect to the $2n$ coordinates. 
 \begin{algorithm}
	\caption{Compute gradients of areas of Voronoi cells}
	\label{alg:Voronoi-grad-areas}
	\begin{algorithmic}[1]
		\Require $\bo x = (x_1,y_1,...,x_n,y_n)$, coordinates of points $p_1,...,p_n$, bounding polygon $\Omega$
		\State Initialize $M=0$ (of size $2n\times n$)
		\State Compute the Voronoi diagram associated to the points $(p_i)_{i=1}^n$ and the intersections of the Voronoi cells with the polygon $\Omega$.
		\State Set \texttt{Voronoi\_ridges} as the set of Voronoi ridges that intersect the bounding polygon $\Omega$.
		\For{\texttt{r} in \texttt{Voronoi\_ridges}}
		    \State For the Voronoi ridge \texttt{r} Find the associated Voronoi points $p_i$ and $p_j$ and the Voronoi vertices $v_k,v_l$.
		    \State Set $\delta = (1,0)$ and compute the contributions $z_n, z_t$ as above.
		    \State Perform the updates using \eqref{eq:grad-contribution}: 
		    \[ M_{2i-1,i} \gets M_{2i-1,i}+z_n+z_t,\ \  M_{2i-1,j} \gets M_{2i-1,j}-z_n-z_t,   \] 
		    \[ M_{2j-1,i} \gets M_{2j-1,i}+z_n-z_t,\ \  M_{2j-1,j} \gets M_{2j-1,j}-z_n+z_t.   \] 
		    \State Set $\delta = (0,1)$ and compute the contributions $z_n, z_t$ as above.
		    \State Perform the updates using \eqref{eq:grad-contribution}: 
		    \[ M_{2i,i} \gets M_{2i,i}+z_n+z_t,\ \  M_{2i,j} \gets M_{2i,j}-z_n-z_t,   \] 
		    \[ M_{2j,i} \gets M_{2j,i}+z_n-z_t,\ \  M_{2j,j} \gets M_{2j,j}-z_n+z_t.   \] 
		\EndFor
		
		\Return $M$
	\end{algorithmic}
\end{algorithm}

\begin{figure}
	\centering
	\includegraphics[height=0.26\textwidth]{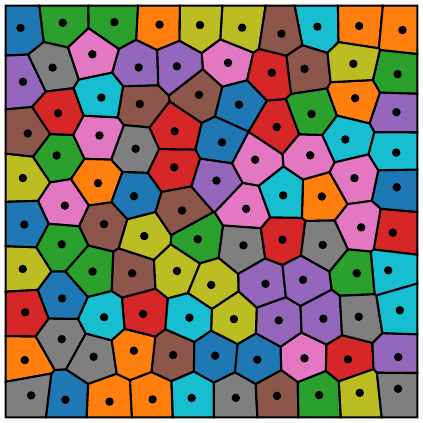} \quad
	\includegraphics[height=0.26\textwidth]{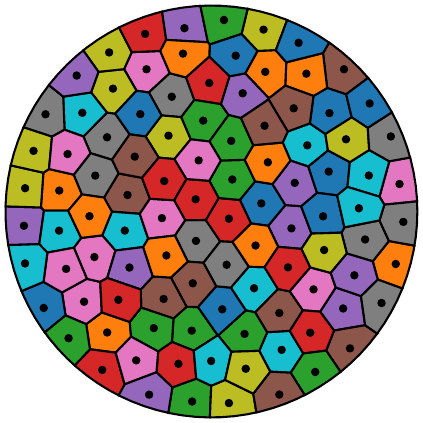} \quad
	\includegraphics[height=0.26\textwidth]{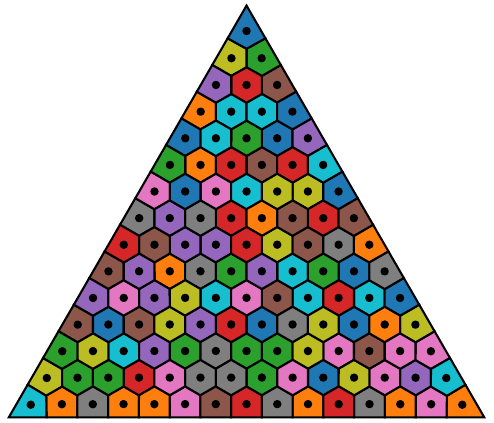} 
	
	\includegraphics[height=0.26\textwidth]{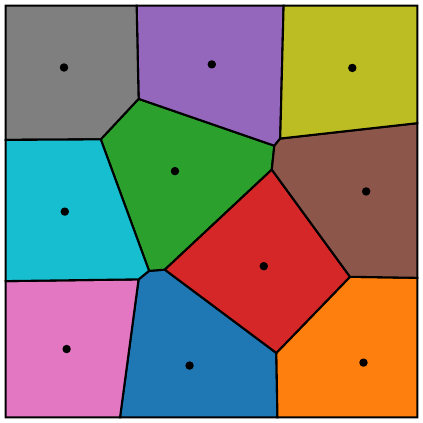} \quad
	\includegraphics[height=0.26\textwidth]{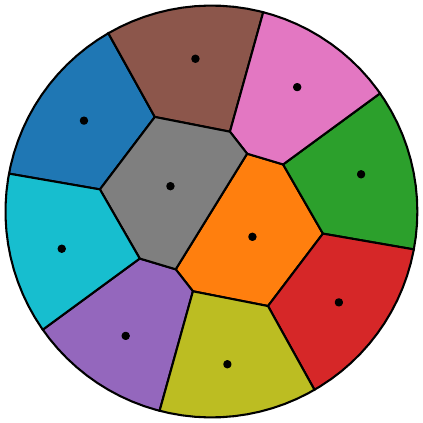} \quad
	\includegraphics[height=0.26\textwidth]{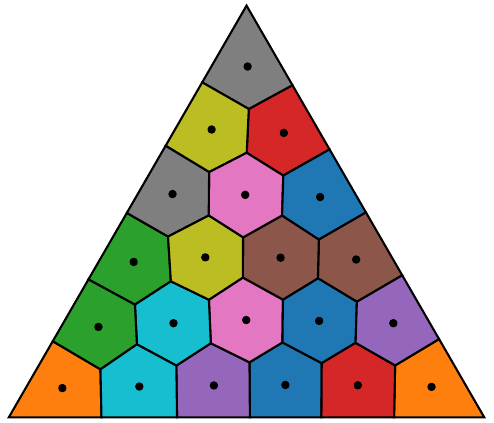} 
	\caption{(top) Voronoi diagrams with more than $100$ cells with equal areas obtained when minimizing \eqref{eq:capacity-constrained-min}. (bottom) Voronoi diagram obtained when minimizing \eqref{eq:Lloyd-functional} under capacity constraints.}
	\label{fig:capacity-constrained-Voronoi}
\end{figure}

The explicit formulas for the gradients of the areas allow us to easily find capacity-constrained Voronoi diagrams as results of an optimization algorithm. For given constraints $|V_i|=c_i|\Omega|$ with $\sum_{i=1}^n c_i=1$, it is enough to minimize the functional
\begin{equation}
 (p_1,...,p_n) \mapsto \sum_{i=1}^n (\text{Area}(V_i)-c_i)^2.
 \label{eq:capacity-constrained-min}
\end{equation}
In order to obtain more regular structures it is also possible to minimize the energy
\begin{equation}
 (p_1,...,p_n) \mapsto \sum_{i=1}^n \int_{V_i} |x-p_i|^2
 \label{eq:Lloyd-functional}
\end{equation}
under the capacity constraints $|V_i| = c_i$. The energy \eqref{eq:Lloyd-functional} is employed for characterizing centroidal Voronoi diagrams where each Voronoi point $p_i$ coincides with the centroid of the cell $V_i$. In particular, Centroidal Voronoi diagrams are critical points for \eqref{eq:Lloyd-functional}. See \cite{xin2016} for more details regarding this functional. Examples in this sense are shown in Figure \ref{fig:capacity-constrained-Voronoi}. The constrained minimization is done using the MMA algorithm \cite{Svanberg2002} from the NLOPT library \cite{nlopt}. Note that all constraints are coded as inequality constraints in this algorithm: $|V_i|\leq c_i$. Since $V_i$ form a partition of $\Omega$ it is immediate to see that if the sets satisfy the inequality constraints, they, in fact, also satisfy the equality constraints $|V_i|=c_i$.

\begin{rem}
	It is also possible to generalize the gradient formulas when a density is involved, when dealing with quantities of the type $\int_{V_i}\rho$, where $\rho \in L^1(\Omega)$ is a given density. The shape derivative of $\int_{V_i} \rho$ is $\int_{\partial V_i} \rho \theta.n$, where $\theta$ is the perturbation of the boundary $\partial V_i$. The boundary perturbations are obviously the same, but the computations in \eqref{eq:area-normal-tangential} are no longer explicit, and a one-dimensional numerical integration needs to be performed for each Voronoi ridge.
\end{rem}

\begin{figure}[!htp]
	\centering
	\includegraphics[width=0.2\textwidth]{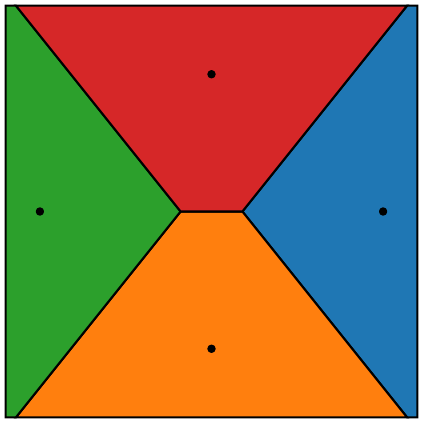} \quad
	\includegraphics[width=0.2\textwidth]{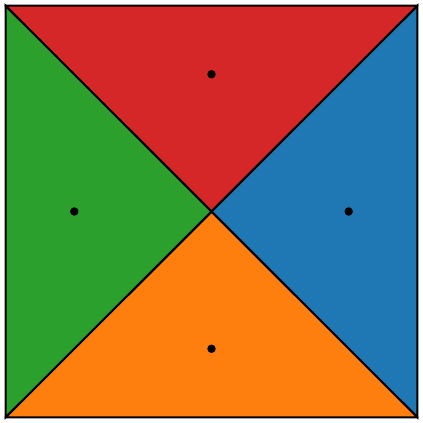} \quad
	\includegraphics[width=0.2\textwidth]{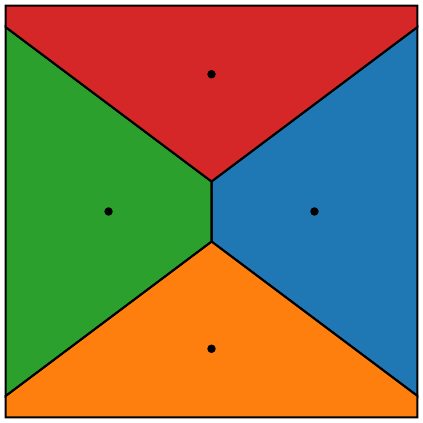} \quad
	\includegraphics[height=0.2\textwidth]{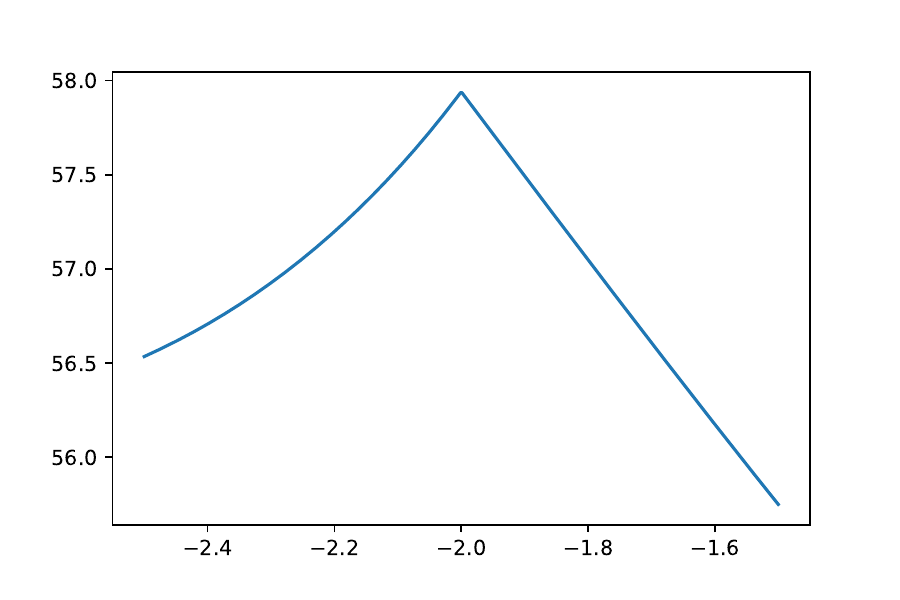} \quad
	\caption{Variation of the perimeter corresponding to a four point singularity in a square. The Voronoi points are $(-t,0),(0,-2),(t,0),(0,2)$ for $t \in [1.5,2.5]$. From left to right: starting configuration, singular configuration, final configuration, graphical representation of the total perimeter.}
	\label{fig:non-smooth-perim}
\end{figure}

{\bf Gradient of the perimeter of the Voronoi cells.} We saw that in order to compute the gradient of the areas of the Voronoi cells, the normal displacement of the Voronoi ridges needed to be understood, when moving the Voronoi points. On the other hand, the variation of the perimeter of a Voronoi region depends on the tangential perturbation of the Voronoi ridges. In order to understand this perturbation one needs to see how the Voronoi vertices move when perturbing the Voronoi points. Moreover, it can be observed that when two Voronoi vertices merge, i.e. a Voronoi ridge collapses, the total perimeter of the cells is not smooth. This behavior is illustrated by an example shown in Figure \ref{fig:non-smooth-perim}.

\begin{figure}[H]
  \centering
  \includegraphics[width=0.3\textwidth]{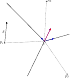}\quad 
  \includegraphics[width=0.3\textwidth]{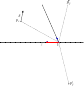}
  \caption{(left) Perturbation of the circumcenter when moving one point and projections on the Voronoi ridges. (right) Computing the perturbation of a boundary point by transforming it into a circumcenter.}
  \label{fig:circumcenter-perturb}
\end{figure}

 Therefore, we suppose in the following that each Voronoi vertex is in contact with at most three Voronoi ridges. Moreover, the definition of the Voronoi cells allows us to conclude that in this situation each Voronoi vertex is the circumcenter of the triangle determined by the three points associated to the neighboring Voronoi regions. This allows us to transform perturbations of the Voronoi points into perturbations of the Voronoi vertices, by looking at the following well known formulas for the circumcenter of a triangle with vertices $(A_x,A_y),(B_x,B_y),(C_x,C_y)$:
\begin{align}
 O_x & = \frac{1}{D}[(A_x^2+A_y^2)(B_y-C_y)+(B_x^2+B_y^2)(C_y-A_y)+(C_x^2+C_y^2)(A_y-B_y)] \notag \\
 O_y & = \frac{1}{D}[(A_x^2+A_y^2)(C_x-B_x)+(B_x^2+B_y^2)(A_x-C_x)+(C_x^2+C_y^2)(B_x-A_x)]
\label{eq:circumcenter}
\end{align}
where $D = 2[A_x(B_y-C_y)+B_x(C_y-A_y)+C_x(A_y-B_y)]$. The formulas above are well defined as long as the three points $A,B,C$ are not colinear. Moreover, it is immediate to see that in this case the circumcenter varies smoothly with respect to the coordinates of the vertices of the triangle. The infinitesimal perturbation of the circumcenter when moving $(A_x,A_y)$ can be computed by simply differentiating the above formulas w.r.t. $A_x$ and $A_y$. Once the derivative of the circumcenter is known, in order to find the gradient of the prerimeter it is enough to project this derivative on all the Voronoi ridges going through the respective circumcenter and add the contribution to the gradient of the perimeter of each cell with respect to the corresponding coordinates. See Figure \ref{fig:circumcenter-perturb} for more details.

Variations induced by the Voronoi nodes are enough to compute the gradient of perimeters of Voronoi cells that do not intersect the boundary of the bounding polygon. For the boundary cells, it is necessary to describe the perturbation of intersections between Voronoi ridges and the bounding polygon. Fortunately, this can also be described using variations of circumcenters for some particular triangles.

Indeed, let $v_kv_l$ be a Voronoi ridge intersecting a side $\ell$ of the bounding polygon $\Omega$ at the point $q$ and let $p_i,p_j$ be the associated Voronoi points. Consider now $p_j'$ the reflection of the point $p_j$ with respect to the line supporting $\ell$. Then obviously $q$ is the circumcenter of the triangle $p_ip_jp_j'$ and the variation of $q$ with respect to perturbations of $p_i$ can be found using the same procedure as above. See Figure \ref{fig:circumcenter-perturb} for more details. The algorithm for computing the gradients for the perimeters of the Voronoi cells is presented in Algorithm \ref{alg:Voronoi-grad-perims}, assuming that every Voronoi vertex is a circumcenter of exactly one triangle determined by the Voronoi points.

\begin{algorithm}
	\caption{Compute gradients of perimeters of Voronoi cells}
	\label{alg:Voronoi-grad-perims}
	\begin{algorithmic}[1]
		\Require $\bo x = (x_1,y_1,...,x_n,y_n)$, coordinates of points $p_1,...,p_n$, bounding polygon $\Omega$
		\State Initialize $M=0$ (of size $2n\times n$)
		\State Compute the Voronoi diagram associated to the points $(p_i)_{i=1}^n$ and the intersections of the Voronoi cells with the polygon $\Omega$.
		\State Set \texttt{Voronoi\_vertices} as the set of Voronoi ridges that intersect the bounding polygon $\Omega$.
		\For{\texttt{v} in \texttt{Voronoi\_vertices}}
		\State Let $p_i,p_j,p_k$ be the three Voronoi points which are associated to ridges going through \texttt{v}. 
		\State Compute the derivative $\vec d$ of the circumcenter of $p_ip_jp_k$ when moving $p_i$ in the direction $\delta = (1,0)$. See Figure \ref{fig:circumcenter-perturb}.
		\State For all ridges \texttt{r} going through \texttt{v} project $\vec d$ on \texttt{r} and add this to the gradient w.r.t. the $x$ coordinate of the perimeter of the cells $\{V_1,V_2\}$ neighbors to the ridge \texttt{r} (determined by the ridge points associated to the ridge \texttt{r}): these are elements $M_{2i-1,V_1}, M_{2i-1,V_2}$ in matrix $M$.
		\State Repeat the above with $\delta = (0,1)$ in order to get the gradients with respect to the $y$-coordinates. 
		\State Do the same instructions for $p_j$ and $p_k$.
		\EndFor
		
		\noindent Set \texttt{Voronoi\_ridges} as the set of Voronoi ridges that intersect the boundary polygon $\Omega$.
		\For{\texttt{r} in \texttt{Voronoi\_ridges}} 
		\State Denote by $p_i, p_j$ the associated ridge points and by $\ell$ the edge of the boudnary polygon $\Omega$ cut by \texttt{r}
		\State Let $p_j'$ be the reflection of $p_j$ with respect to $\ell$.
		\State For $\delta = (1,0)$ compute the derivative $\vec d$ of the circumcenter of $p_ip_jp_j'$ when moving $p_i$ in the direction $\delta$. See Figure \ref{fig:circumcenter-perturb}.
		\State Project $\vec d$ on the ridge \texttt{r} and add this projection to the gradient of the cells $i$ and $j$ w.r.t. the $x$ coordinate: $M_{2i-1,i}$ and $M_{2i-1,j}$ in matrix $M$.
		\State Project $\vec d$ on $\ell$ and add this to the gradient of the cells $i$ and $j$ (with the proper sign).
		\State Repeat the above with $\delta=(0,1)$ in order to get the gradients with respect to the $y$ coordinates.
		\State Do the same instructions for $p_j$.
		\EndFor
		
		\Return $M$
	\end{algorithmic}
\end{algorithm}

\begin{figure}[H]
	\centering
	\includegraphics[height=0.2\textwidth]{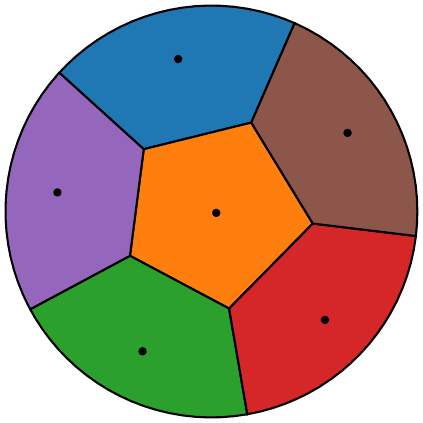}\quad
	\includegraphics[height=0.2\textwidth]{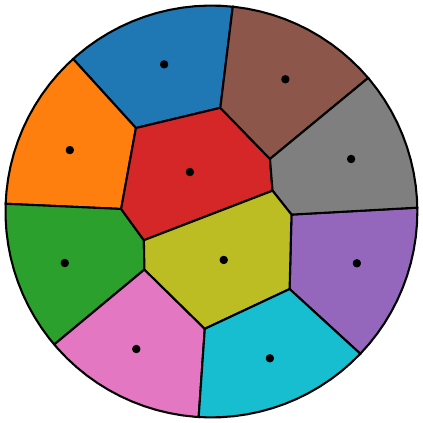}\quad
	\includegraphics[height=0.2\textwidth]{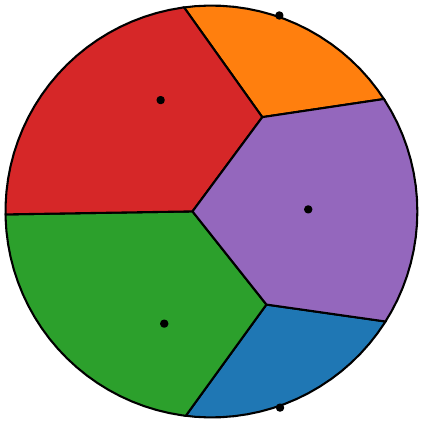}\quad
	\includegraphics[height=0.2\textwidth]{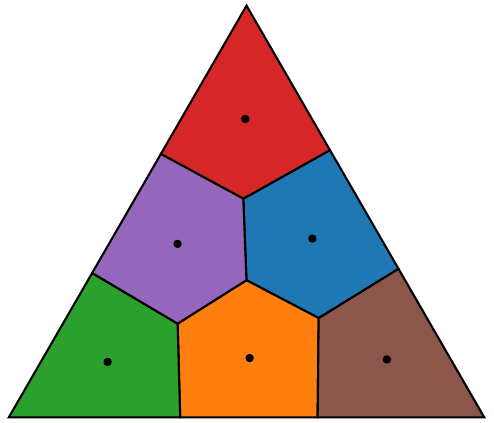}
	\caption{Initializations obtained when minimizing the perimeter of Voronoi cells under area constraints. The areas are equal, except the third case where two cells have areas three times smaller than the other three.}
	\label{fig:examples-initializations}
\end{figure}

Using the gradients for the area and perimeters of Voronoi cells it is possible to perform a constrained minimization of the perimeter under area constraint starting from random Voronoi initializations. The optimization is performed with Nlopt \cite{nlopt} optimization toolbox in Python using the MMA \cite{Svanberg2002} algorithm. Some examples of initializations obtained are shown in Figure \ref{fig:examples-initializations}. The initial Voronoi points are chosen randomly inside the polygon $\Omega$. In order to accelerate the convergence of the optimization algorithm a few iterations of Lloyd's algorithm are performed before starting the optimization process. Recall that Lloyd's algorithm consists in replacing the Voronoi points by the centroids of the respective cells iteratively (see for example \cite{xin2016} for more details). In order to deal with local minima multiple optimizations (typically $10$) are performed for every polygon $\Omega$ and the one with the partition having the least perimeter is retained as a valid initialization. Note that the algorithm gives similar topologies with the best known ones shown in \cite{cox-partitions} for the case of equal areas and in \cite{cox-nonequal} for the case of cells with two different areas. 

{\bf Initialization of a partition.} Having at our disposal the gradients of areas and perimeters of Voronoi cells, we are now ready to propose initialization algorithms for optimal partitioning algorithm. In practice we use one of the options below:
\begin{enumerate}[noitemsep,topsep=0pt]
	\item Compute minimizers of \eqref{eq:capacity-constrained-min} starting from random Voronoi points $(p_i)$. Repeat the procedure a number of times and keep the configuration having the smallest total perimeter. This works well when the areas of the cells are the same.
	\item Optimize the total perimeter of the Voronoi cells under capacity constraints starting from random Voronoi points $(p_i)$. Repeat the procedure a number of times and keep the configuration having the smallest total perimeter. This approach gives good results when the areas of the cells are different and the optimization process is more difficult, since more local minima are present.
	\item When $n\leq 4$ random initializations work very well.
	\item In dimension three random initializations were used for $n\leq 4$ and random Voronoi initializations were used for $n\geq 5$. 
\end{enumerate}

\subsection{Shape derivative}
\label{sec:sh-deriv}

In order find perturbations of the domain $\Omega$ that decrease the value of a functional $J(\Omega)$ the concept of shape derivative is used. For a Lipschitz domain $\Omega$, the functional $J$ is said to be shape differentiable at $\Omega$ if there exists a linear form $\theta \mapsto J'(\Omega)(\theta)$ such that for every vector field $\theta \in W^{1,\infty}(\Bbb R^d,\Bbb R^d)$ we have
\[ J((I+\theta)(\Omega)) = J(\Omega)+J'(\Omega)(\theta)+o(\|\theta\|_{W^{1,\infty}}),\]
where $I$ denotes the identity mapping. Classical results from \cite[Section 2.31]{sokolowski-zolesio} and \cite[Section 5.2]{henrot-pierre-english} show that for a function $f_\Omega \in H^1(\Omega)$ that varies smoothly with respect to perturbations of $\Omega$, the functional $J(\Omega) = \int_\Omega f_\Omega$ is shape differentiable with 
\begin{equation}
  J'(\Omega)(\theta) = \int_\Omega f'_\Omega(\theta)+ \int_{\partial \Omega} f_\Omega \theta.n.
  \label{eq:sh-deriv-integral}
\end{equation}
 A nice overview of the basic notions regarding shape derivatives, together with the associated references is given in the paper \cite{sh-deriv-gauss}.

We are interested in finding the shape derivatives of $I_\varepsilon(\Omega,c)$ and $PI_\varepsilon(\Omega,\bo c)$ which are minimal values obtained through constrained optimization. In the following we perform these computations assuming that the corresponding shape derivatives exist. Methods used in the following are inspired from the following works:
\begin{itemize}[noitemsep,topsep=0pt]
	\item In \cite[Chapter 10, Sections 2.3, 5.4]{delfour-zolesio} the shape derivation of a minimum problem and, respectively, a saddle point is described.
	\item In \cite[Capter 3]{fiacco} the derivative of the minimal value problem of a constrained problem where the objective and the constraints depend on a parameter is given in the finite dimensional case.
	\item In \cite[Chapter 4]{bonnans} the differentiability of the minimal value given by a constrained parametric problem is considered in Banach spaces.
\end{itemize}
In our case we are in the framework of the derivation of a saddle point. However, we were not able apply the first result cited above to prove rigorously that the shape derivative exists. In the following we use a formal approach in order to identify the formula for the shape derivatives that are of interest for us.

{\bf The case of one phase.} Consider $u_\Omega$ which minimizes $F_\varepsilon(u)$ from Theorem \ref{thm:modica-mortola} and suppose that $u_\Omega$ is unique. We assume that $u_\Omega$ varies smoothly with respect to perturbations of the boundary of $\Omega$ and that its shape derivative $u_\Omega'(\theta)$ exists and belongs to $H^1(\Omega)$. Remark \ref{rem:smoothness-MM} underlines the fact that $u_\Omega$ is $C^\infty$ in the interior of $\Omega$ and has the regularity of $\Omega$ up to the boundary. In the following we suppose that $\Omega$ convex (or, in general, at least of class $C^2$), which implies that $u_\Omega$ is indeed in $H^2(\Omega)$ and the gradient $\nabla u_\Omega$ has a well defined trace on $\partial \Omega$. The function $I_\varepsilon(\Omega,c) = F_\varepsilon(u_\Omega)$ has the same structure as in \eqref{eq:sh-deriv-integral}. Using the classical chain rule for the shape derivatives (see for example \cite[Lemma 26]{sh-deriv-gauss}) we obtain that
\begin{align}
 I_\varepsilon (\Omega,c)'(\theta) =  (F_\varepsilon(u_\Omega))'(\theta)&= \int_{\partial \Omega} \left(\varepsilon |\nabla u_\Omega|^2 +  \frac{1}{\varepsilon}W(u_\Omega)\right) \theta.n \notag \\  &+ \int_\Omega \left(2\varepsilon \nabla u_\Omega\cdot \nabla u'_\Omega(\theta)+\frac{1}{\varepsilon}W'(u_\Omega)u'_\Omega(\theta)\right) 
 \label{eq:sh-deriv-init}
\end{align}
where $u_\Omega'(\theta)$ is the shape derivative of $u_\Omega$ with respect to $\theta$. 

Recall (see Remark \ref{rem:smoothness-MM}) that minimizing $F_\varepsilon(u)$ under the constraint $\int_\Omega u = c|\Omega|$ implies the existence of a Lagrange multiplier $\mu \in \Bbb{R}$ such that 
\[ \int_\Omega \left(2\varepsilon \nabla u_\Omega \cdot \nabla \phi + \frac{1}{\varepsilon}W'(u_\Omega)\phi+\mu \phi\right) = 0,\]
for every $\phi \in H^1(\Omega)$. Taking $\phi = u'_\Omega(\theta)$ in the previous equation gives 
\begin{equation} \int_\Omega \left(2\varepsilon \nabla u_\Omega \cdot \nabla u_\Omega'(\theta) + \frac{1}{\varepsilon}W'(u_\Omega)u'_\Omega(\theta)\right)=-\mu \int_\Omega u'_\Omega(\theta)
\label{eq:remove-euler-deriv}
\end{equation}

Recall that $u_\Omega$ also verifies the constraint $\int_\Omega u_\Omega -c|\Omega|=0$. Differentiating this with respect to the shape $\Omega$ gives
\begin{equation}
\int_\Omega u'_\Omega(\theta)+\int_{\partial \Omega} (u_\Omega-c)\theta.n = 0.
\label{eq:sh-deriv-constraint}
\end{equation}
Combining \eqref{eq:sh-deriv-init}, \eqref{eq:remove-euler-deriv} and \eqref{eq:sh-deriv-constraint} gives 
\begin{equation}
I_\varepsilon(\Omega,c)'(\theta) = \int_{\partial \Omega} \left(
\varepsilon |\nabla u_\Omega|^2 +  \frac{1}{\varepsilon}W(u_\Omega)+\mu(u_\Omega -c)
\right)\theta.n 
\label{eq:sh-deriv-single}
\end{equation}
The previous formula corresponds to the shape derivative of the Lagrangian 
\[ \mathcal L(\Omega,u,\nu) = \int_\Omega \left( \varepsilon|\nabla u|^2+\frac{1}{\varepsilon}W(u)+\nu(u-c)\right)\]
with respect to the shape when $u = u_\Omega$ and $\nu = \mu$. This is in accord with similar results in \cite[Chapter 10, Section 5.4]{delfour-zolesio}, \cite[Chapter 4]{bonnans} and \cite[Chapter 3]{fiacco}.

Note that the formula \eqref{eq:sh-deriv-single} does give a valuable and reasonable assumption on the perturbation of the boundary $\Omega$ that increases the length of the set $SF(\Omega,c)$ having minimal perimeter. Notice that the term $\left(\varepsilon |\nabla u_\Omega|^2 + \frac{1}{\varepsilon}W(u_\Omega)\right)$ is non-zero (and strictly positive) only in the neighborhood of the contact points of the minimal relative perimeter set with the boundary $\partial \Omega$. Moving the boundary outwards at these points with a small enough step size will increase the minimal perimeter $I_\varepsilon (\Omega,c)$. 

On the other hand the term containing $\mu(u_\Omega-c)$ models the movement of the relative isoperimetric set when the boundary of $\Omega$ is perturbed away from the contact points. Indeed, recall that $|SF(\Omega,c)|=c|\Omega|$ and the Lagrange multiplier $\mu$ is proportional (as $\varepsilon \to 0$) to the curvature of the isoperimetric set $SF(\Omega,c)$, as shown in Remark \ref{rem:lag-mult}. Therefore, when perturbing $\Omega$ away from the contact points with $SF(\Omega,c)$ the inner boundary of $SF(\Omega,c)$ is \emph{pulled} in one direction or another (corresponding to the change in the volume) and the corresponding variation of the minimal perimeter is given by the mean curvature of the isoperimetric set.

\begin{figure}
	\centering
	\includegraphics[height=0.29\textwidth]{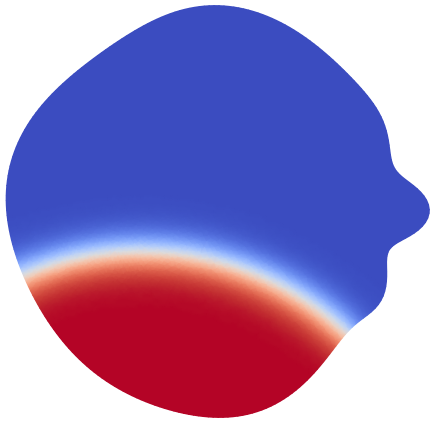}\quad
	\includegraphics[height=0.3\textwidth]{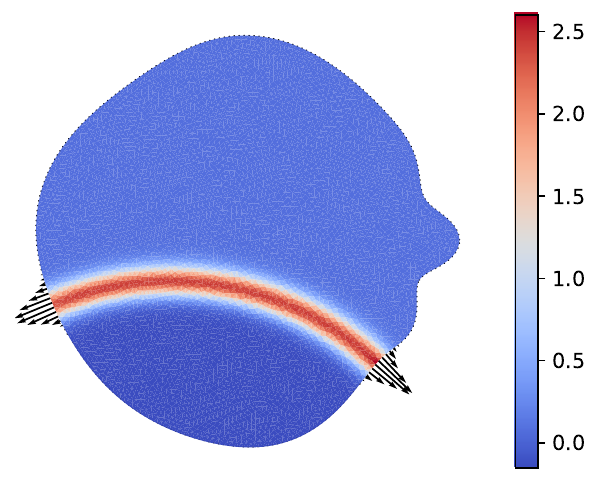}
	\caption{(left) Minimization of the Modica-Mortola functional $F_\varepsilon$ with integral constraint $0.3|\Omega|$. (right) Representation of the shape gradient and the normal perturbation producing an ascent direction for $SF(\Omega,0.3)$.}
	\label{fig:exampleOnePhase}
\end{figure}

In the case the solution $u_\Omega$ is not unique, we cannot assume that $u_\Omega$ varies smoothly with $\Omega$. Indeed, perturbing the boundary of $\Omega$ with the normal velocity given by $\left(\varepsilon |\nabla u_\Omega|^2 + \frac{1}{\varepsilon}W(u_\Omega)\right)+\mu (u_\Omega-c)$ may drastically change the topology of the minimal set. Nevertheless, perturbing the boundary of $\Omega$ with this normal velocity will eventually increase the value of $I_\varepsilon (\Omega,c)$.
In Figure \ref{fig:exampleOnePhase} the numerical approximation of $u_\Omega$ is shown together with the values of $\left(\varepsilon |\nabla u_\Omega|^2 + \frac{1}{\varepsilon}W(u_\Omega)\right)+\mu(u_\Omega-c)$. Perturbing $\Omega$ in the normal direction as shown will increase the minimal value, provided the solution $u_\Omega$ is unique. It can be noted that the term $\left(\varepsilon |\nabla u_\Omega|^2 + \frac{1}{\varepsilon}W(u_\Omega)\right)$ is dominant in the shape derivative, as the second term is of order $1/\varepsilon$. Recall that the results given in \cite{luckhaus-modica} show that the the Lagrange multiplier $\mu$ as a bounded limit as $\varepsilon\to 0$, proportional to the constant mean curvature of the limiting minimal interface.

\begin{rem}
	Multiple approaches may exist in order to make the above computation a rigorous proof of the shape derivative formula and we describe a few below. The technical difficulties involved do not allow us to present such a complete proof.
	\begin{itemize}[topsep=0pt,noitemsep]
		\item In \cite[Chapter 9, Theorem 5.1]{delfour-zolesio} a method that computes the directional derivative of a saddle point is described. However, it is not clear if this result applies to our case. On the other hand, one may note that this method would give the same result as the formal method described previously.
		\item The computation shown above strongly depends on the existence of the shape derivative $u_\Omega'(\theta)$. The existence of this derivative is not obvious even under the assumption that $u_\Omega$ is a unique minimizer. It may be possible to apply techniques similar to those in \cite[Section 5.7]{henrot-pierre-english} or \cite[Section 2.29]{sokolowski-zolesio} in order to deduce that the unique minimizer $u_\Omega$ is differentiable with respect to the shape $\Omega$. 
	\end{itemize}
\end{rem}

\begin{rem}
	The hypothesis regarding the uniqueness of $u_\Omega$ for the shape derivative to exist is similar to the hypothesis needed when differentiating the eigenvalue of an operator with respect to the shape. When dealing with multiple eigenvalues, the shape derivative does not exist, but directional derivatives are available. For more details see \cite[Chapter 5]{henrot-pierre-english}. It is possible that such theoretical results could be obtained in our case, but this goes outside the scope of this article.
\end{rem}

{\bf The case of partitions.} Consider $\bo u_\Omega=(u_\Omega^i)_{i=1}^n$ which minimizes $G_\varepsilon $ in Theorem \ref{thm:gconv-part}. As in the previous paragraphs, we suppose that $\bo u_\Omega$ is unique and varies smoothly with respect to perturbations in $\Omega$. Differentiating with respect to the domain using the formula \eqref{eq:sh-deriv-integral} we get 
\begin{align}
PI_\varepsilon (\Omega,\bo c)'_\Omega(\theta) =  ( G_\varepsilon (\bo u_\Omega))_\Omega'(\theta)&= \int_{\partial \Omega} \sum_{i=1}^n \left(\varepsilon |\nabla u_\Omega^i|^2 +  \frac{1}{\varepsilon}W(u_\Omega^i)\right) \theta.n\notag  \\   &+ \sum_{i=1}^n \int_\Omega \left(2\varepsilon \nabla u_\Omega^i \cdot \nabla (u_\Omega^i)'(\theta) +\frac{1}{\varepsilon}W'(u_\Omega^i)(u_\Omega^i)'(\theta)\right)
\label{eq:sh-deriv-part-init}
\end{align}
Since $\bo u_\Omega$ minimizes $G_\varepsilon (\bo u)$ under the constraints $\int_{\Omega} u_i = c_i|\Omega|$ and $\sum_{i=1}^n u_i=1$ there exist Lagrange multipliers $\mu_i\in \Bbb{R}$, $i=1,...,n$ and $\lambda \in L^1(\Omega)$ such that 
\begin{equation}
\int_\Omega \sum_{i=1}^n\left( 2\varepsilon \nabla u_\Omega^i \cdot \nabla \phi_i + \frac{1}{\varepsilon}W'(u_\Omega^i)\phi_i+\mu_i \phi_i+\lambda \phi_i\right) = 0 \text{ for every } \phi_i \in H^1(\Omega),\ i=1,...,n.
\label{eq:lagmult-partitions}
\end{equation}
Since the sum constraint $\sum_{i=1}^n u_i=1$ implies that one of the area constraints is redundant we may note that the Lagrange multipliers are not uniquely defined. Adding a constant to $\lambda$ and subtracting the same constant from each $\mu_i$ gives another set of valid multipliers. Therefore, it is not restrictive to assume that the multiplier $\lambda$ verifies $\int_\Omega \lambda = 0$.

Replacing $\phi_i = (u_\Omega^i)'(\theta)$ in \eqref{eq:lagmult-partitions} and using \eqref{eq:sh-deriv-part-init} we obtain
\begin{align*}
PI_\varepsilon (\Omega,\bo c)'_\Omega(\theta) &= \int_{\partial \Omega} \sum_{i=1}^n \left(\varepsilon |\nabla u_\Omega^i|^2 +  \frac{1}{\varepsilon}W(u_\Omega^i)\right) \theta.n\notag  \\   &-\sum_{i=1}^n \mu_i \int_\Omega (u_\Omega^i)'(\theta)  -\int_\Omega \lambda \sum_{i=1}^n (u_\Omega^i)'(\theta)
\end{align*}

On the other hand, differentiating the constraint $\sum_{i=1}^n u_\Omega^i = 1$ gives $\sum_{i=1}^n (u_\Omega^i)'(\theta)=0$. Moreover, differentiating the constraints $\int_\Omega u_\Omega^i - c_i|\Omega| = 0$ we obtain 
\[ \int_\Omega (u_\Omega^i)'(\theta) = -\int_{\partial \Omega} u_\Omega^i\theta.n  +c_i\int_{\partial \Omega} \theta.n.\]
Therefore we obtain
\begin{align}
PI_\varepsilon (\Omega,\bo c)'_\Omega(\theta) &= \int_{\partial \Omega} \sum_{i=1}^n \left(\varepsilon |\nabla u_\Omega^i|^2 +  \frac{1}{\varepsilon}W(u_\Omega^i)\right) \theta.n\notag  \\   &+\sum_{i=1}^n \mu_i \int_{\partial \Omega} (u_\Omega^i-c_i).
\theta.n
\label{eq:sh-deriv-part}
\end{align}
The Lagrange multipliers $\mu_i$ can be found by using $\phi_i = \delta_{ij}$ in \eqref{eq:lagmult-partitions}, which gives $\mu_i = -1/(\varepsilon|\Omega|)\int_\Omega W'(u_\Omega^i)$.

\begin{figure}
	\centering
	\includegraphics[height=0.29\textwidth]{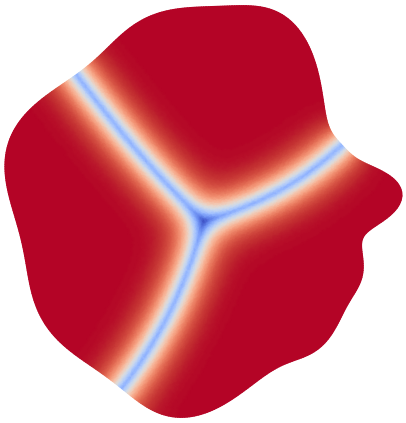}\quad
	\includegraphics[height=0.3\textwidth]{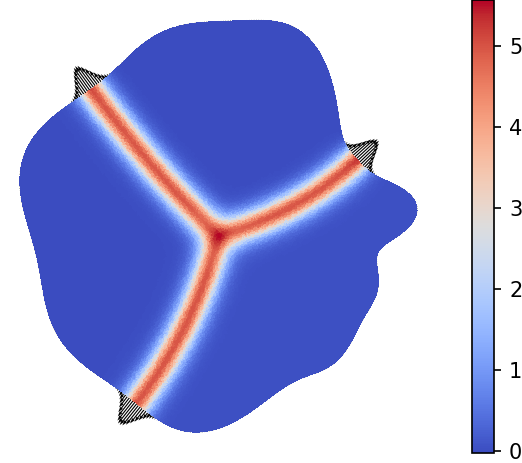}
	\caption{(left) Minimization of the minimal partition functional $ G_\varepsilon $ for three cells with equal areas. (right) Representation of the shape gradient and the normal perturbation producing an ascent direction for $SF(\Omega,\bo c)$.}
	\label{fig:examplePartition}
\end{figure}

As discussed above, in the case $\bo u_\Omega$ is not unique, perturbing the boundary of $\Omega$ with normal velocity given by $\sum_{i=1}^n\left(\varepsilon |\nabla u_\Omega^i|^2 + \frac{1}{\varepsilon}W(u_\Omega^i)+\mu_i(u_\Omega^i-c_i)\right)$ will eventually increase the value of $PI_\varepsilon (\Omega,\bo c)$ and will reduce the multiplicity of the family of optimal partitions $\bo u_\Omega$. In Figure \ref{fig:examplePartition} the numerical approximation of $\bo u_\Omega$ is shown together with the value of $\sum_{i=1}^n\left(\varepsilon |\nabla u_\Omega^i|^2 + \frac{1}{\varepsilon}W(u_\Omega^i)+\mu_i(u_\Omega^i-c_i)\right)$. Perturbing $\Omega$ in the normal direction as shown will increase the length of the minimal partition, provided the solution $\bo u_\Omega$ is unique.

\begin{figure}
	\centering 
    \includegraphics[width=0.4\textwidth]{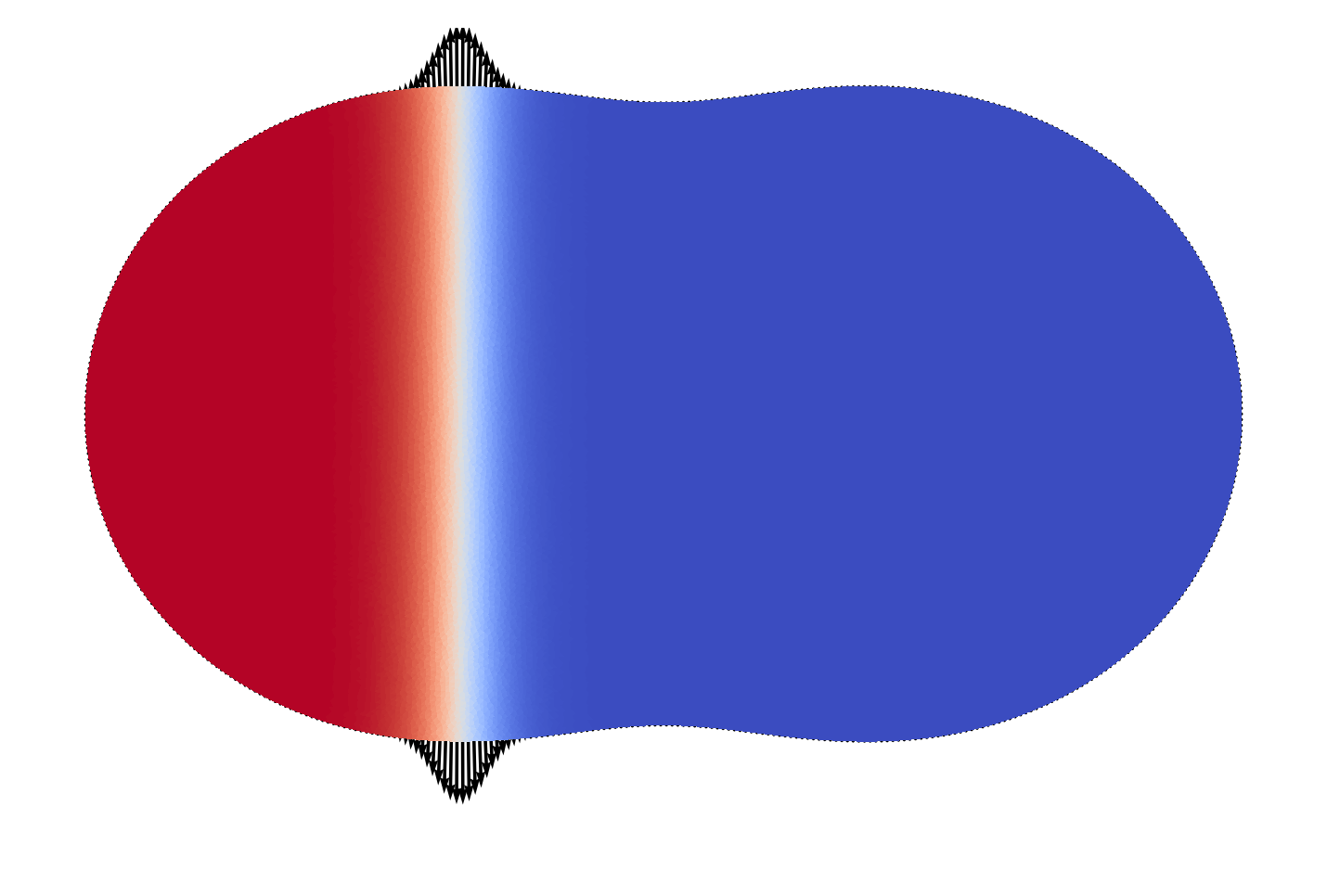}\quad
    \includegraphics[width=0.4\textwidth]{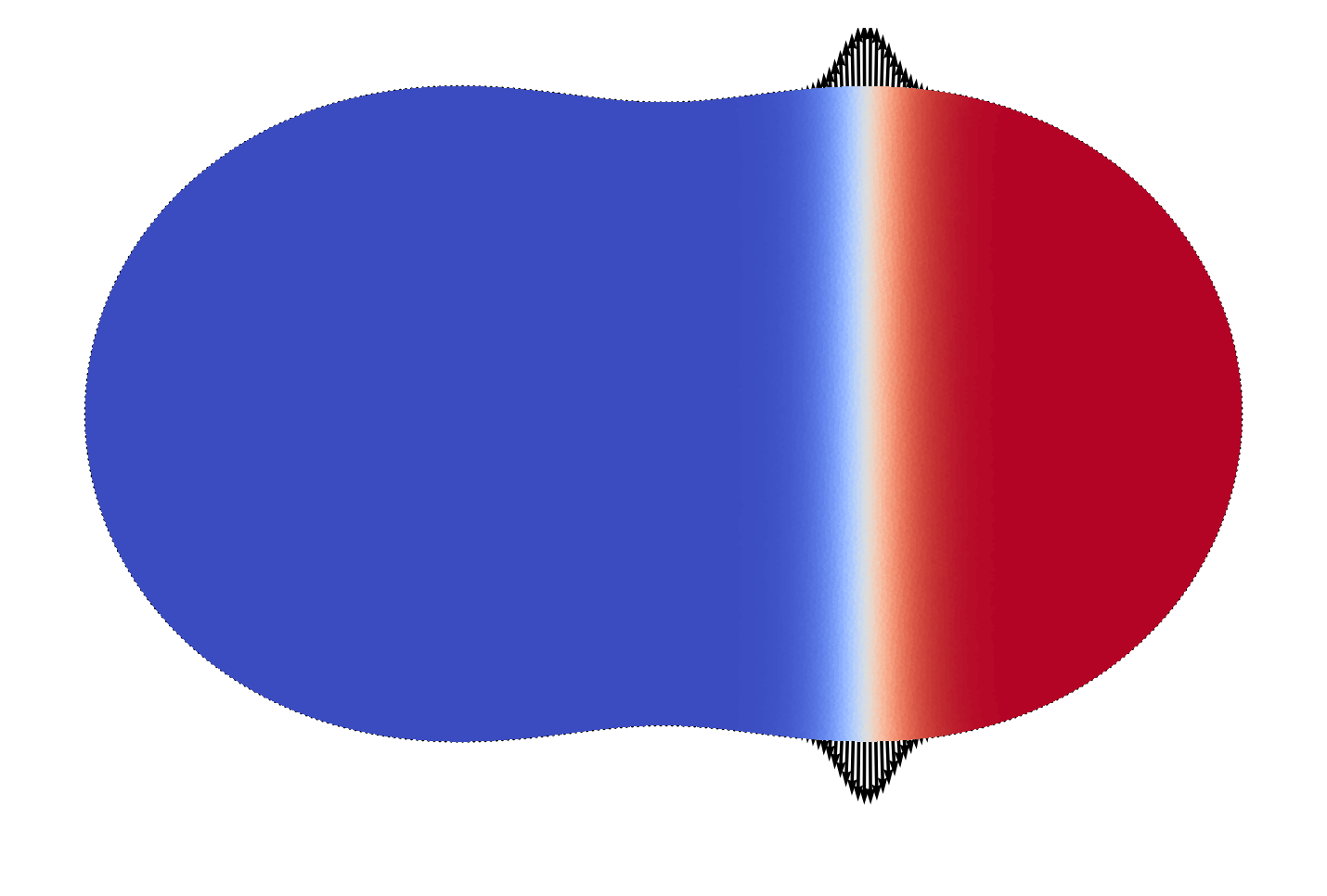}
    \caption{Symmetric domain with two minimizers $I_\varepsilon (\Omega,c)$. The optimal densities obtained numerically are represented together with the perturbation fields obtained from \eqref{eq:sh-deriv-single}}
    \label{fig:example_non_differentiability}
\end{figure}

{\bf An example of non differentiability - multiplicity greater than one.} It is not difficult to imagine domains $\Omega$ for which there are multiple minimizers for $I_\varepsilon (\Omega,c)$, given $c>0$. It is enough to consider $c$ small enough and a symmetric domain like in Figure \ref{fig:example_non_differentiability}. Moreover, there exists a value of $c$ for which the Lagrange multiplier vanishes: $\mu = -1/(\varepsilon|\Omega|)\int_\Omega W'(u_\Omega)=0$. It is enough to remember that $\mu$ is proportional to the curvature of the minimal interface as $\varepsilon\to 0$. Let us show that in such a case the functional $I_\varepsilon (\Omega,c)$ does not admit a shape derivative. Indeed, suppose that $\Omega$ is symmetric like in Figure \ref{fig:example_non_differentiability}. Denote $u_1,u_2$ the two solutions and $\theta_1,\theta_2$ two vector fields such that $\theta_i.n = \varepsilon |\nabla u_i|^2+\frac{1}{\varepsilon} u_i^2(1-u_i)^2$ (also illustrated in Figure \ref{fig:example_non_differentiability}). 

Suppose that $J(\Omega) := I_\varepsilon (\Omega,c)$ is differentiable at $\Omega$. Then for $\{i,j\} = \{1,2\}$ it is clear that $F_\varepsilon((I+t\theta_i)(\Omega),u_j)=F_\varepsilon(\Omega,u_j)$, i.e. the minimal value of $F_\varepsilon$ does not change when modifying $\Omega$ with only one of the vector fields $\theta_i$. This would imply that $J'(\Omega)(\theta_1) = J'(\Omega)(\theta_2) = 0$ and by linearity $J'(\Omega)(\theta_1+\theta_2) = 0$. However, this last equality is clearly false, since if we modify $\Omega$ with the combined vector field $\theta_1+\theta_2$ we clearly have 
\[ J((I+(\theta_1+\theta_2)(\Omega)) = J(\Omega)+\int_{\partial \Omega} \left(\varepsilon |\nabla u_1|^2 + \frac{1}{\varepsilon}u_1^2(1-u_1)^2\right) \theta_1.n + o(\|\theta_1\|_{W^{1,\infty}}). \]
In order to deduce the above equality it is enough to work with one half $\Omega_i$ of the domain $\Omega$. We use the fact that this half can be extended to a non-symmetric domain for which $u_i$ is a unique minimizer and apply the formula for the shape derivative found previously. 

Therefore, we arrive at a contradiction, showing that when multiple minimizers of $F_\varepsilon(\Omega,u)$ exist, the functional $I_\varepsilon (\Omega,c)$ is not shape differentiable. The same kind of argument can be applied for $PI_\varepsilon (\Omega,\bo c)$ in the case where optimal partitions are not unique.

\subsection{Radial parametrization and optimization algorithm}
\label{sec:optim-algo}

The results of \cite{longest_shortest_fence} and existence results obtained in Section \ref{sec:theory-part} are restricted to convex domains $\Omega$. We therefore choose to search for domains maximizing $SF(\Omega, c)$ and $SP(\Omega,\bo c)$ in the larger class of star-shaped domains which includes the class of convex sets. These domains can be parametrized using an associated radial function in dimensions two and three. Furthermore, a spectral decomposition of the radial function with a Finite number of Fourier coefficients is used in order to work with a finite, but sufficiently large number of parameters in the computations.

{\bf Planar domains.} In dimension two, the radial function $\rho : [0,2\pi] \to \Bbb{R}_+$ is discretized using $2N+1$ Fourier coefficients
\[ \rho(t) = a_0+\sum_{k=1}^N (a_k \cos(kt) + b_k \sin(kt)).\]
Consider a shape functional $J(\Omega)$ whose shape derivative is expressed by $J'(\Omega)(\theta) = \int_{\partial \Omega} \mathcal G \theta.n$. Using the discretization above, given $\bo v = (a_0,a_1,...,a_N,b_1,...,b_N)$ that defines $\Omega$ via the radial function $\rho$, a finite dimensional function is obtained $j(\bo v) = J(\Omega)$. It is classical to compute the gradient of $j$ using the shape derivative, by choosing the appropriate boundary perturbation for each Fourier coefficient. Using the notation $\bo r = x/|x|$ we have $\bo r.n = \rho/\sqrt{\rho^2+(\rho')^2}$. Therefore, denoting $v_n=\rho/\sqrt{\rho^2+(\rho')^2}$, we obtain 
\begin{equation} \frac{\partial j}{\partial a_k} = \int_{\partial \Omega} \mathcal G \cos(k t) v_n \text{ and }\frac{\partial j}{\partial b_k} = \int_{\partial \Omega} \mathcal G \sin(kt) v_n. 
\label{eq:grad-2D}
\end{equation}

{\bf Domains in $\Bbb{R}^3$.} In dimension three we choose to parametrize the unit sphere using $(\phi,\psi) \in [-\pi,\pi]\times [0,2\pi] \mapsto (\cos \psi \cos \phi,\sin \psi \cos \phi,\sin \phi)$. Next, we are interested in parametrizing radial functions $\rho :  [-\pi,\pi]\times [0,2\pi]$ which are constant for $\phi \in \{-\pi,\pi\}$. This is needed in order to be able to create 3D meshes in FreeFEM \cite{freefem} by deforming two dimensional meshes. One way of attaining this objective is to use two dimensional Fourier parametrizations which contain only sines for the $\phi$ coordinate, together with an affine function in $\phi$ in order to allow different values at the extremities $\phi \in \{-\pi,\pi\}$:
\[ \rho(\phi,\psi) = a\phi+b + \sum_{k=1}^N\sum_{l=1}^M \left(c_{k,l}\sin(2 k \phi)\cos(l\psi) +d_{k,l}\sin(2 k \phi)\sin(l\psi)\right).\]
As in dimension two, it is straightforward to infer the gradient of the discretized functional with respect to each one of the parameters. A simple computation yields $v_n=\bo r.n = \rho/\sqrt{\rho^2+(\rho'_\theta)^2/\cos^2\phi + (\rho'_\phi)^2}$:
\[ \frac{\partial j}{\partial a} = \int_{\partial \Omega} \mathcal G \phi v_n,\ \frac{\partial j}{\partial b} = \int_{\partial \Omega} \mathcal G v_n ,\]
\begin{equation} \frac{\partial j}{\partial c_{k,l}} = \int_{\partial \Omega} \mathcal G \sin(2k\phi)\cos (l\psi) v_n ,  \  \frac{\partial j}{\partial d_{k,l}} = \int_{\partial \Omega} \mathcal G \sin(2k\phi)\sin (l\psi)v_n
\label{eq:grad-3D}
\end{equation}

{\bf Optimization algorithm.} Given the discretization and the gradients expressed above it is straightforward to implement a gradient descent algorithm. The delicate issue is the fact that at each iteration, the objective function and its gradient are computed as a result of a minimization algorithm. If a local minimal perimeter set/partition is found instead of the global one, this might give a wrong ascent direction. Therefore, we choose to work with a \emph{gradient flow} type algorithm, which consists in advancing at each iteration in the direction given by the gradient of the functional with a prescribed step, regardless of the fact that the objective function increases or decreases. In this way, even the optimization algorithms solved at one of the iterations yields a local minimum, the global optimization algorithm may still correct itself at subsequent iterations. The area constraint is imposed by a projection algorithm: the next iterate is rescaled to have the desired area via a homothety. The precise description is given in Algorithm \ref{alg:global-max}.

\begin{algorithm}
	\caption{Global maximization algorithm}
	\label{alg:global-max}
	\begin{algorithmic}[1]
		\Require Initial Fourier coefficients, area constraints ($c$ for one phase, a vector $\bo c$ for the partitions), the number of iterations \texttt{Niter}, $\varepsilon$, initial step $\alpha$, the number of iterations \texttt{Nmod} after which the step is halved
		\For{\texttt{i} in \{1,2,...,\texttt{Niter}\}}
		\State Construct the mesh of $\Omega$ from the Fourier coefficients $\bo v$: the size of triangles/tetrahedra should be at most $\varepsilon/2$.
		\State Approximate $I_\varepsilon (\Omega,c)$ (or $PI_\varepsilon (\Omega,\bo c)$ in the case of partitions)
		\State Compute the gradient $\nabla j(\bo v)$: use \eqref{eq:grad-2D} or \eqref{eq:grad-3D} with $\mathcal G$ given by \eqref{eq:sh-deriv-single} (or \eqref{eq:sh-deriv-part} for the partitions case)
		\State Advance in the direction of the gradient in order to increase the value of $j(\bo v)$:
		\[ \bo v \gets \bo v + \alpha \nabla j(\bo v).\]
		\State Project on the area/volume constraint of $\Omega$ using a homothety
		\State If $i \mod \texttt{Nmod} \equiv 0$ decrease the step: $\alpha \gets \alpha/2$.
		\EndFor
		
		\Return the final set of Fourier coefficients $\bo v$
	\end{algorithmic}
\end{algorithm}

\begin{figure}
	\centering
	\includegraphics[height=0.3\textwidth]{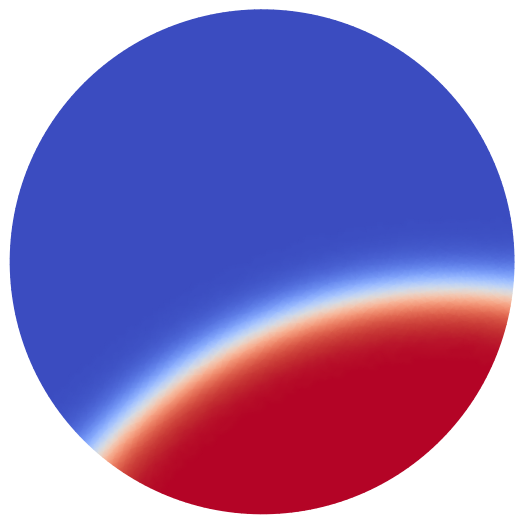}\quad
	\includegraphics[height=0.3\textwidth]{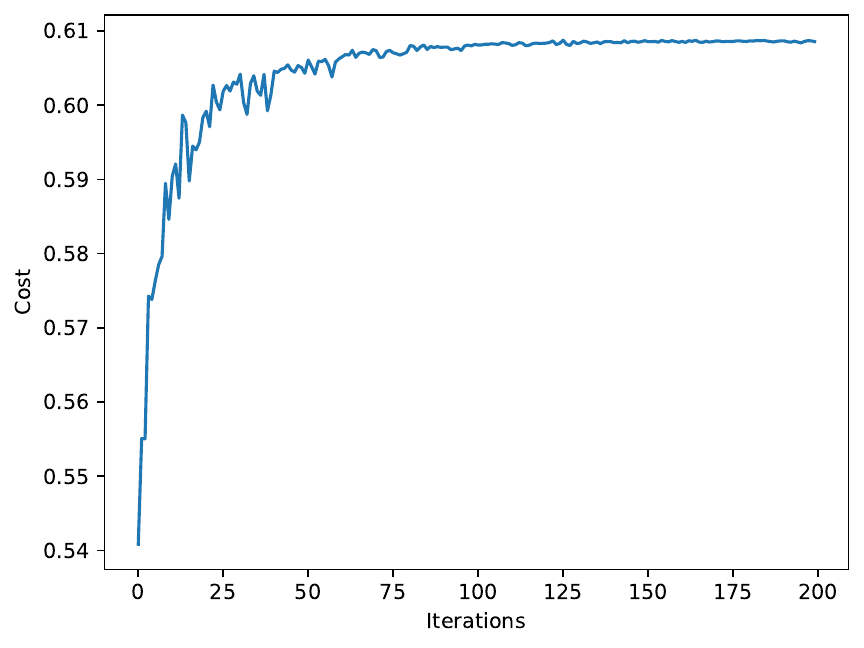}
	\caption{Maximization of $I(\Omega,0.3)$ in dimension two together with the evolution of the cost function.}
	\label{fig:example-opt-single}
\end{figure}

\begin{figure}
	\centering
	\includegraphics[height=0.3\textwidth]{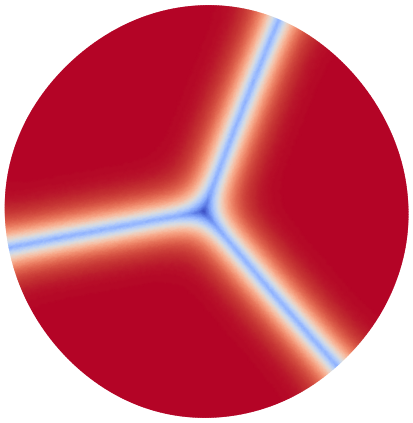}\quad
	\includegraphics[height=0.3\textwidth]{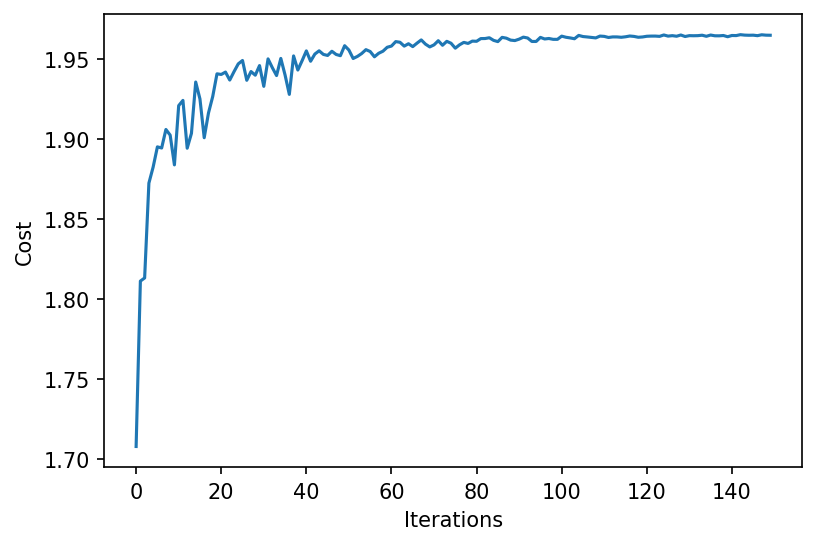}
	\caption{Maximization of $PI(\Omega,(1/3,1/3,1/3))$ in dimension two together with the evolution of the cost function.}
	\label{fig:example-opt-part}
\end{figure}

\begin{figure}
	\centering
	\begin{tabular}{cccccc}
	\includegraphics[height=0.14\textwidth]{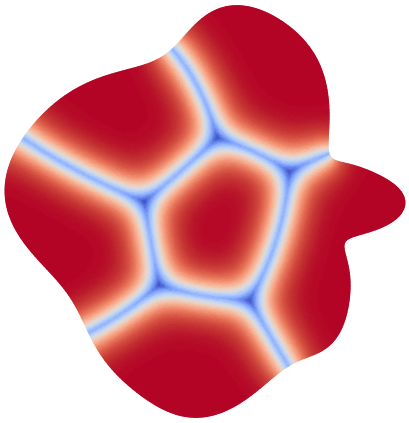}&
	\includegraphics[height=0.14\textwidth]{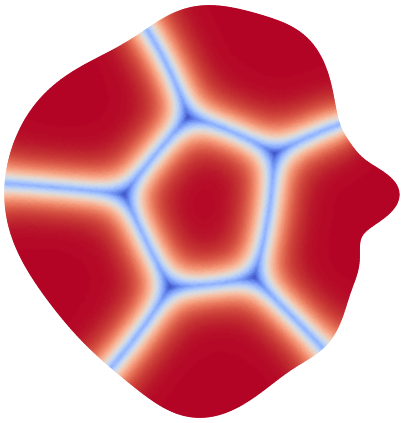}&
	\includegraphics[height=0.14\textwidth]{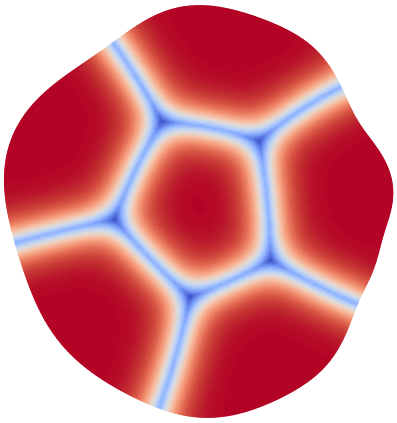}&
	\includegraphics[height=0.14\textwidth]{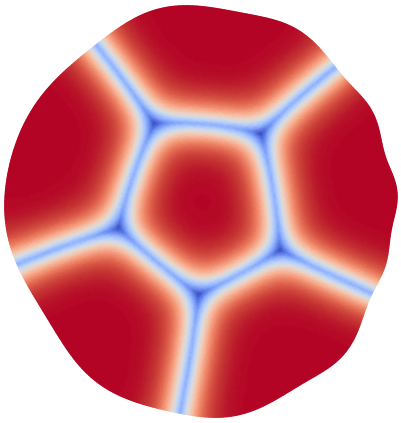}&
	\includegraphics[height=0.14\textwidth]{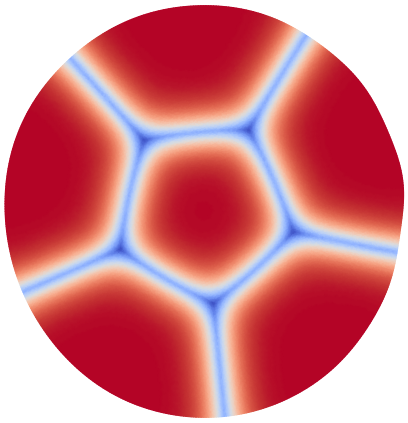}&
	\includegraphics[height=0.14\textwidth]{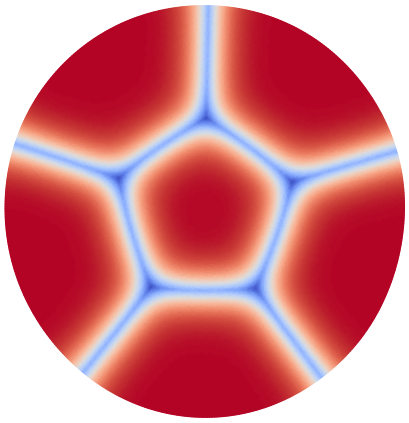} \\
	Iter $1$: $3.305$ &
	Iter $6$: $3.383$ &
	Iter $13$: $3.384$ & 
	Iter $20$: $3.422$ & 
	Iter $70$: $3.441$ & 
	Iter $150$: $3.450$
		\end{tabular}
	\caption{Illustration of the gradient flow algorithm in dimension two for $n=6$: the numerical optimal partition and its associated cost are represented for a couple of iterations.}
	\label{fig:gflow}
\end{figure}

\begin{figure}
	\centering
	\begin{tabular}{cccccc}
		\includegraphics[height=0.14\textwidth]{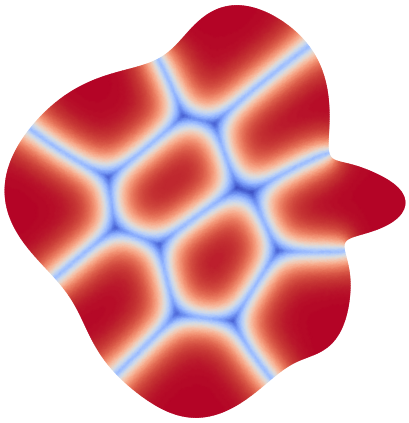}&
		\includegraphics[height=0.14\textwidth]{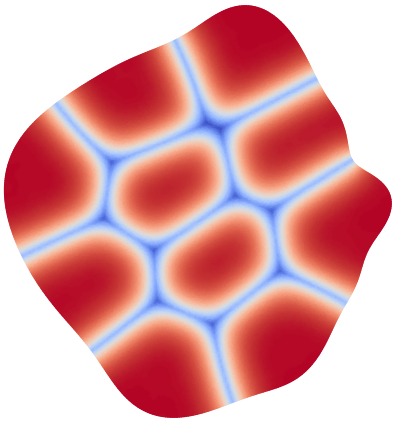}&
		\includegraphics[height=0.14\textwidth]{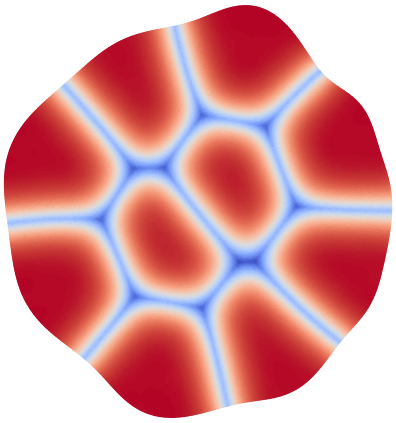}&
		\includegraphics[height=0.14\textwidth]{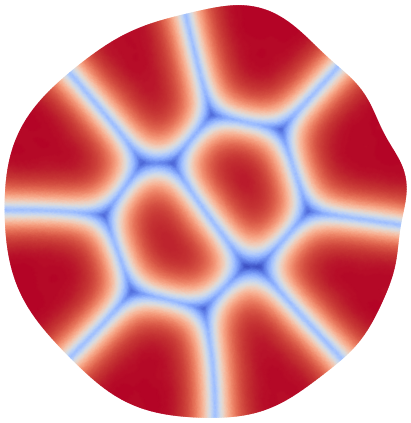}&
		\includegraphics[height=0.14\textwidth]{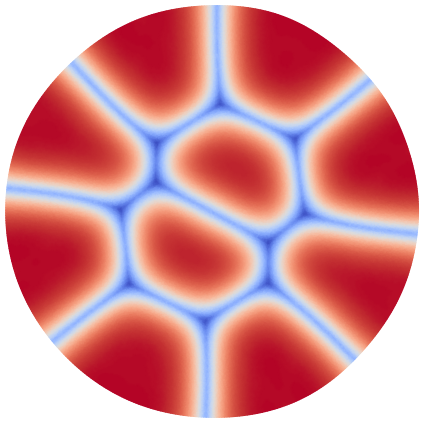}&
		\includegraphics[height=0.14\textwidth]{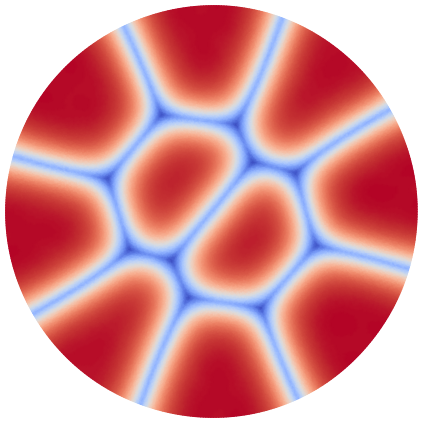} \\
		Iter $1$: $4.655$ &
		Iter $5$: $4.738$ &
		Iter $15$: $4.784$ & 
		Iter $20$: $4.861$ & 
		Iter $70$: $4.895$ & 
		Iter $150$: $4.902$
	\end{tabular}
	\caption{Illustration of the gradient flow algorithm in dimension two for $n=10$: the numerical optimal partition and its associated cost are represented for a couple of iterations.}
	\label{fig:gflow10}
\end{figure}

An example of a result obtained with this algorithm for the maximization of $I_\varepsilon (\Omega,c)$ is shown in Figure \ref{fig:example-opt-single} together with the graph of the objective function. It can be seen that the objective function increases and stabilizes as the size of the step decreases. Oscillations in the curve describing the cost have two main causes: first, the optimization algorithm at the current iteration might yield a local minimum instead of the global one and secondly, the size of the step may be too big. An example for the case of partitions is shown in Figure \ref{fig:example-opt-part}. Multiple instances of the gradient flow maximization algorithm are represented in Figure \ref{fig:gflow} for $n=6$ and in Figure \ref{fig:gflow10} for $n=10$.

{\bf Numerical aspects.} When minimizing $F_\varepsilon$ and $G_\varepsilon$ it is classical to consider meshes with elements that have size smaller than $\varepsilon$. This is due to the fact that the phase transition from $0$ to $1$ typically takes place in a region of width proportional to $\varepsilon$ and the mesh needs to be fine enough to capture this. In dimension two we consider $\varepsilon=0.05$ giving rise to meshes having around $23k$ nodes. 

In dimension three using $\varepsilon=0.1$ gives meshes of about $25k$ nodes. When dealing with more cells in dimension three we start with $\varepsilon=0.07$ and we interpolate and re-optimize the result on a finer mesh corresponding to $\varepsilon=0.04$. This gives meshes with around $35k$ nodes. For postprocessing and plotting purposes, the final mesh is further refined using MMG3D \cite{MMG} such that more tetrahedra are present where phases change quickly. The final partition is interpolated and re-optimized (with $\varepsilon=0.025$) on this fine mesh (with around $270k$ nodes) before plotting.

{\bf Code.} The finite element software used for the optimization algorithm described in Section \ref{sec:numerics-partitions} is FreeFEM \cite{freefem}, which provides an interface to the LBFGS optimizer from Nlopt \cite{nlopt}. 

The partition initialization via Voronoi diagrams is coded in Python, where optimization algorithms from \texttt{Scipy.optimize} and \texttt{Nlopt} are used for unconstrained and, respectively, constrained optimizations. Codes and examples are provided in the following Github repository: \href{https://github.com/bbogo/LongestShortestPartitions/tree/main/GradientVoronoi}{\nolinkurl{https://github.com/bbogo/LongestShortestPartitions/tree/main/GradientVoronoi}}.

The visualization is done with Python using Matplotlib in dimension two and Mayavi \cite{mayavi} in dimension three. The graphical representation of partitions is done by extracting surface meshes of an iso-level for each cell in the optimal partition using FreeFEM \cite{freefem} and MMG3D \cite{MMG}. These surface meshes are then plotted with Mayavi \cite{mayavi}.

Some codes used for obtaining the results illustrated in the paper can be found in the Github repository: \href{https://github.com/bbogo/LongestShortestPartitions/tree/main/FreeFEMcodes}{\nolinkurl{https://github.com/bbogo/LongestShortestPartitions/tree/main/FreeFEMcodes}}.

\section{Results}
\label{sec:results}

\begin{figure}
	\centering
	\includegraphics[width=0.3\textwidth]{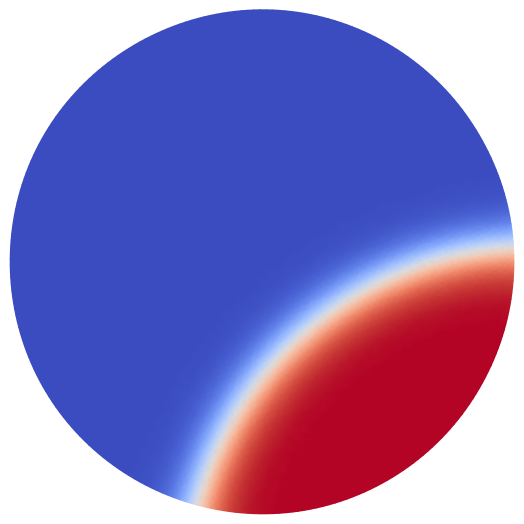}\quad
	\includegraphics[width=0.3\textwidth]{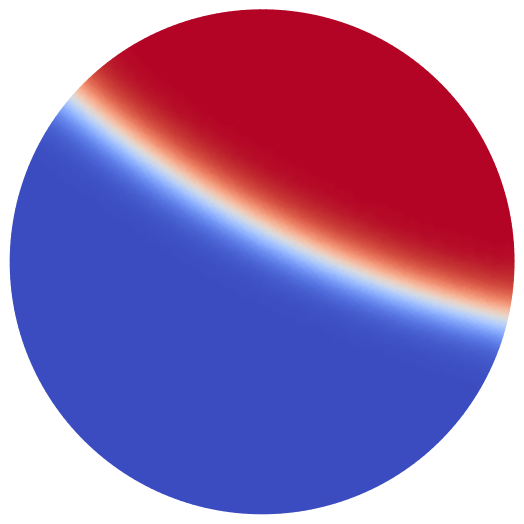}\quad
	\includegraphics[width=0.3\textwidth]{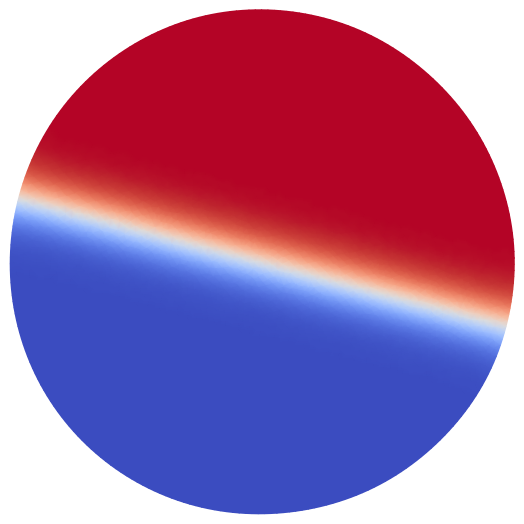}

	\includegraphics[width=0.3\textwidth]{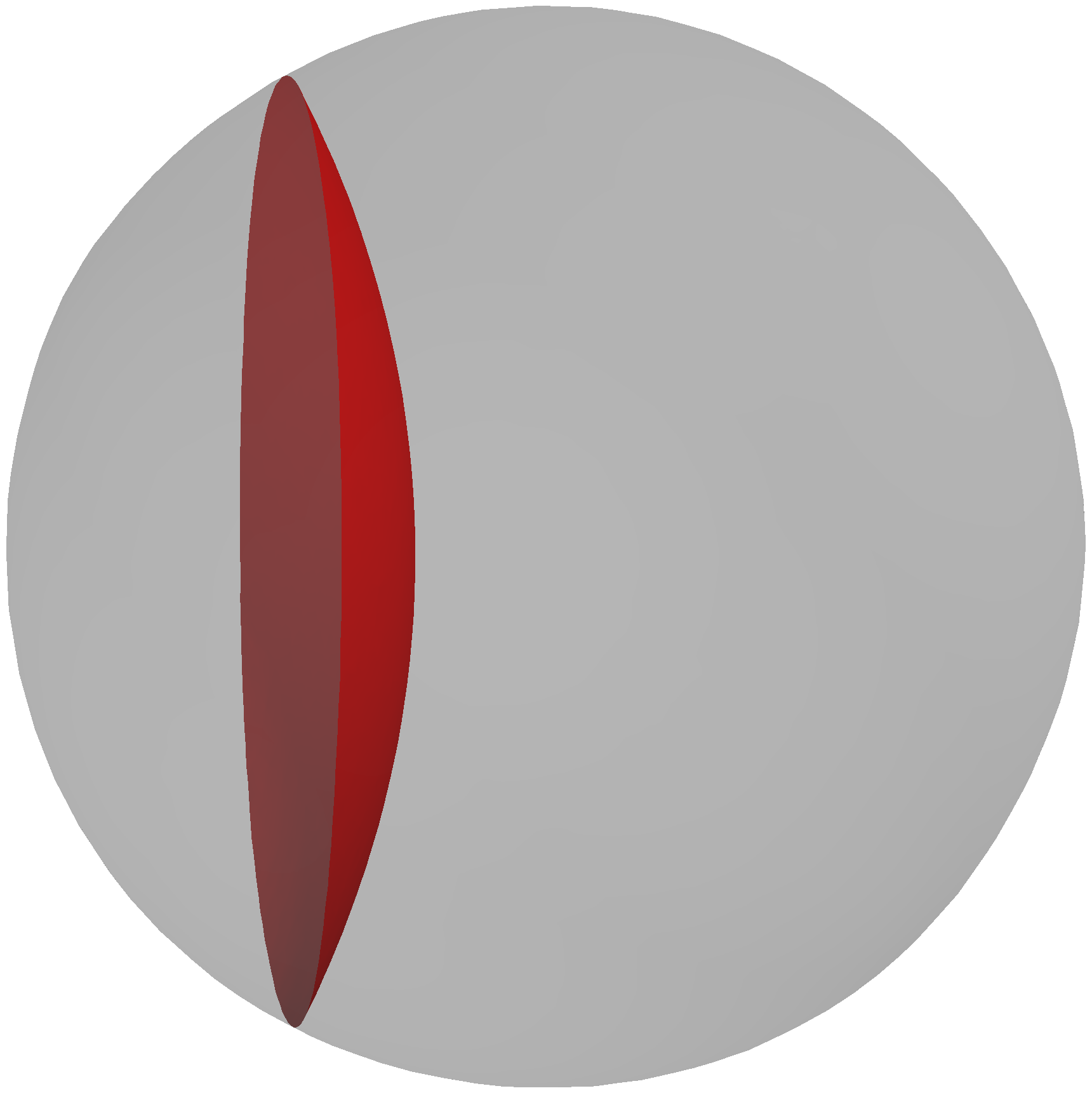}\quad
	\includegraphics[width=0.3\textwidth]{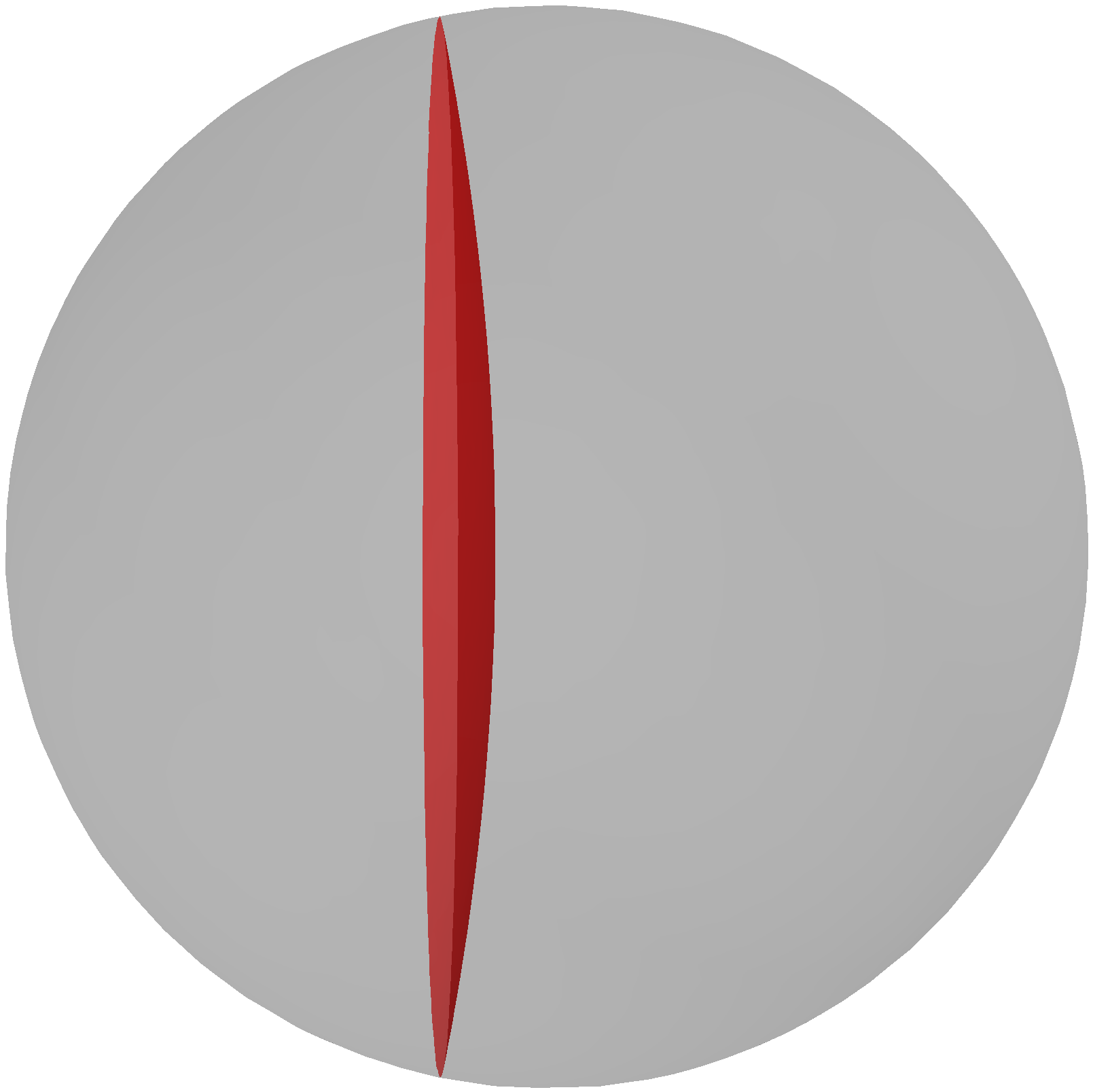}\quad
	\includegraphics[width=0.3\textwidth]{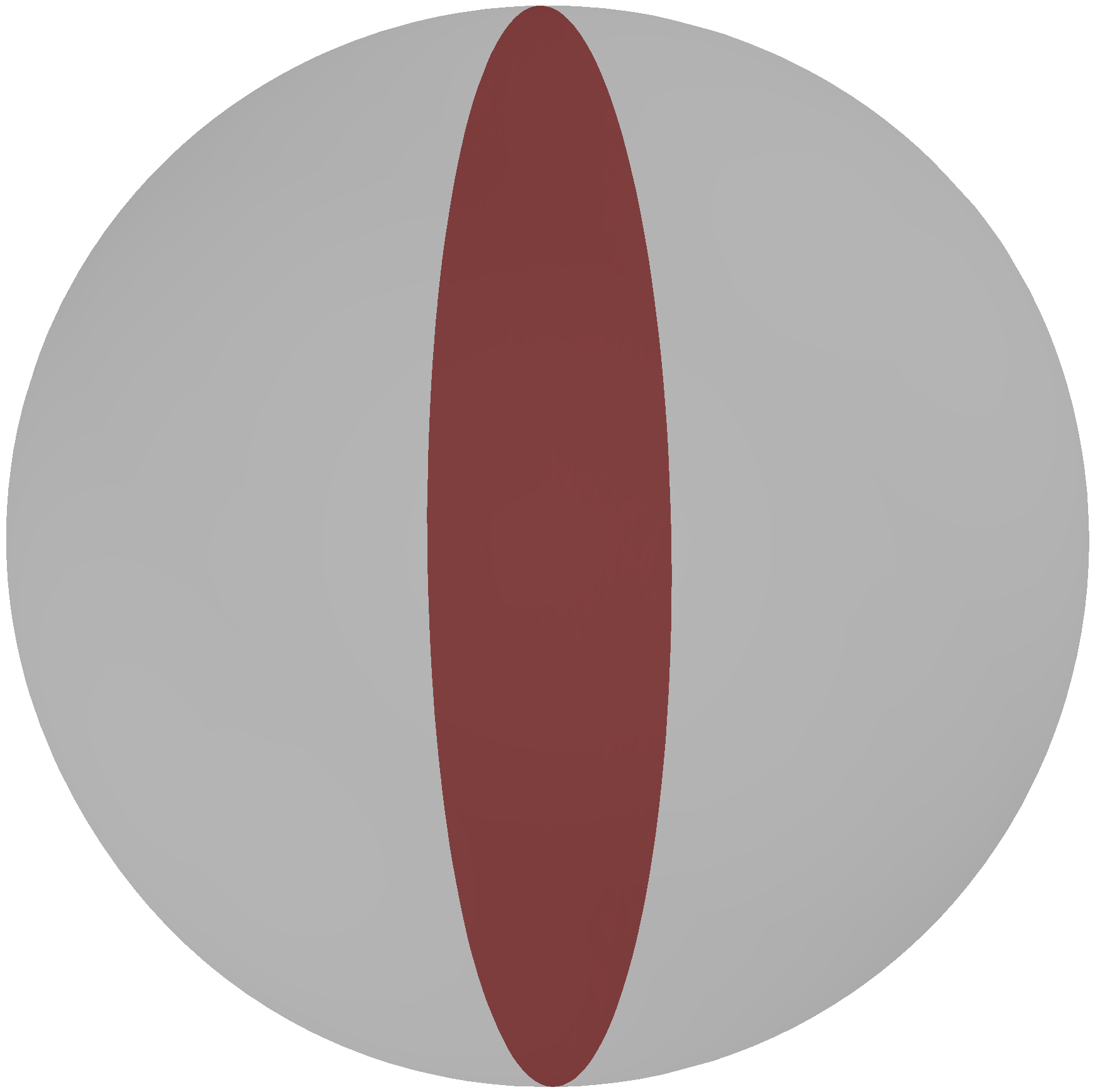}
	\caption{Maximization of the minimal relative perimeter in 2D and 3D with volume constraints $c \in \{0.25,0.4,0.5\}$. The optimal set $\Omega$ (the disk/ball) together with the set obtained numerically when minimizing the relative perimeter for the given volume fraction.}
	\label{fig:results-one-phase}
\end{figure}

In this section we use the algorithm described previously in order to study problems \eqref{eq:longest_shortest_fence} and \eqref{eq:longest_shortest_partition}. Results from \cite{longest_shortest_fence} show that problem \eqref{eq:longest_shortest_fence} is solved by the disk in dimension two for $c=1/2$. We perform simulations for various values of $c<1/2$ (note that considering $c$ or $1-c$ for the constraint gives the same result) and the numerical result is always the disk in dimension two. In dimension three the same phenomenon occurs: for various values of the volume fraction $c$ the shape which maximized the relative minimal perimeter of a subset with volume $c|\Omega|$ is the ball. Some examples are shown in Figure \ref{fig:results-one-phase}.

\begin{figure}
	\centering 
	\includegraphics[width=0.3\textwidth]{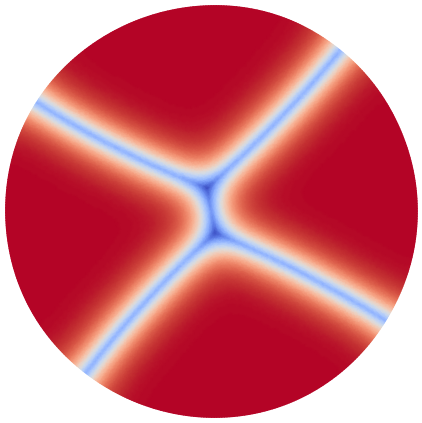}\quad
	\includegraphics[width=0.3\textwidth]{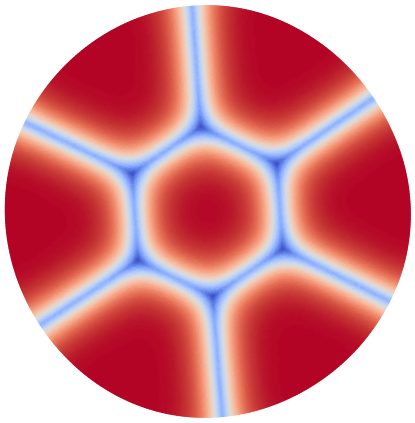}\quad
	\includegraphics[width=0.3\textwidth]{GFlow150_10}\quad
	\caption{Maximization of the length of the minimal perimeter partition into equal areas for $n \in \{4,7,10\}$}
	\label{fig:part-2D-equal}
\end{figure}

\begin{figure}
	\centering
	\includegraphics[width=0.2\textwidth]{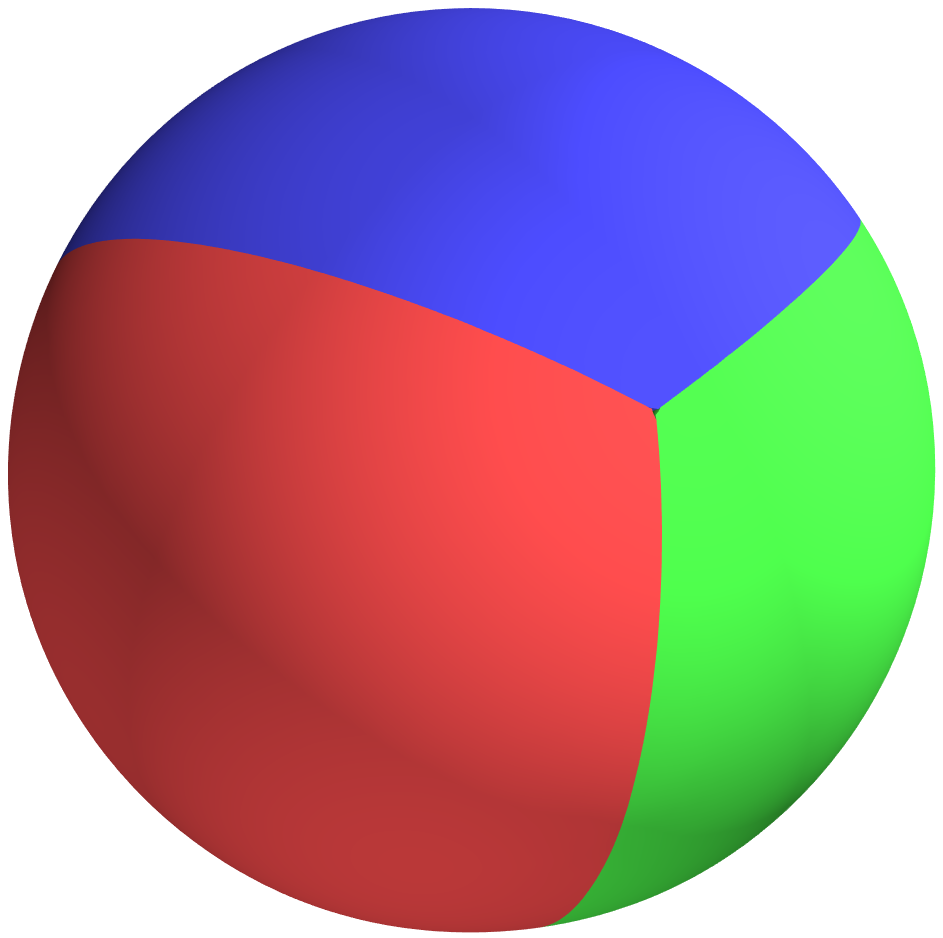}\quad
	\includegraphics[width=0.2\textwidth]{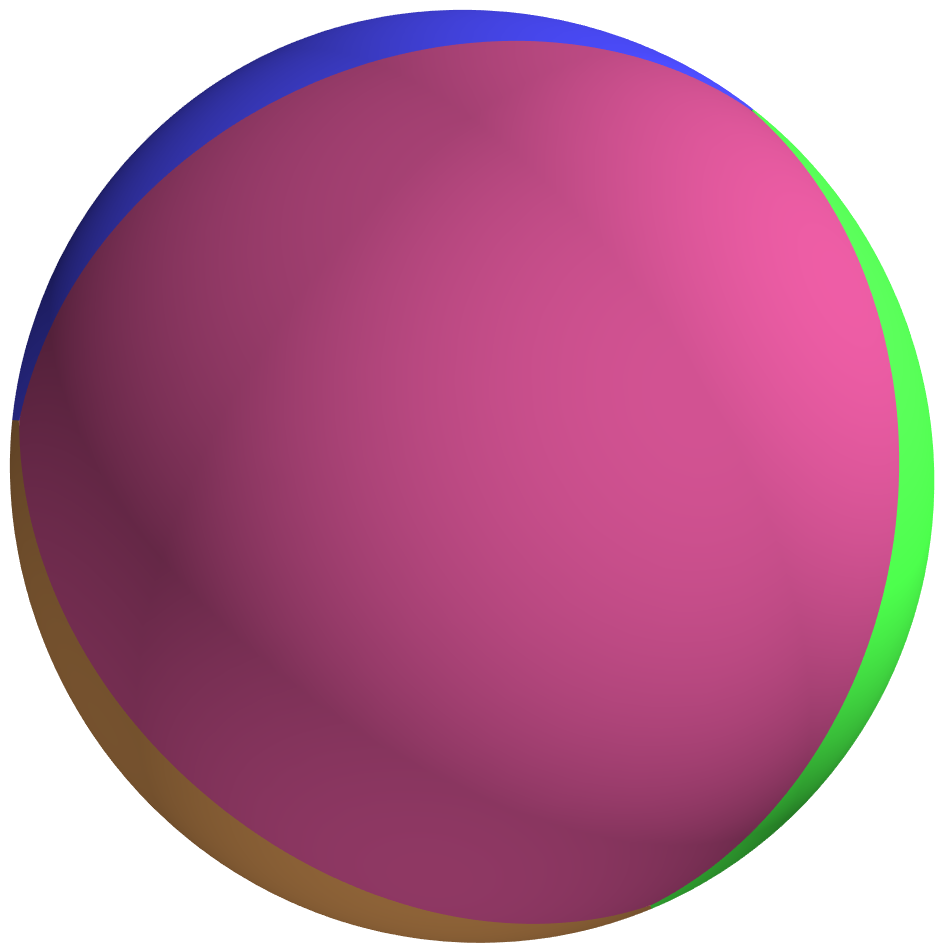}\qquad
	\includegraphics[width=0.2\textwidth]{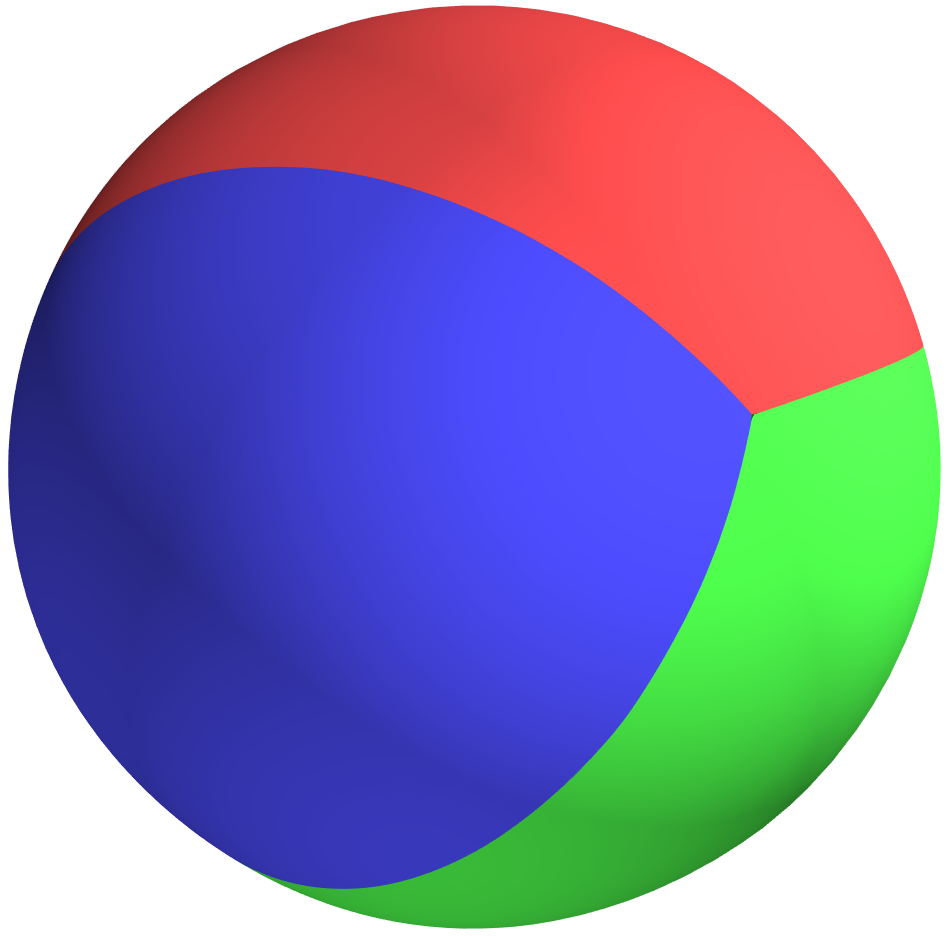}\quad
	\includegraphics[width=0.2\textwidth]{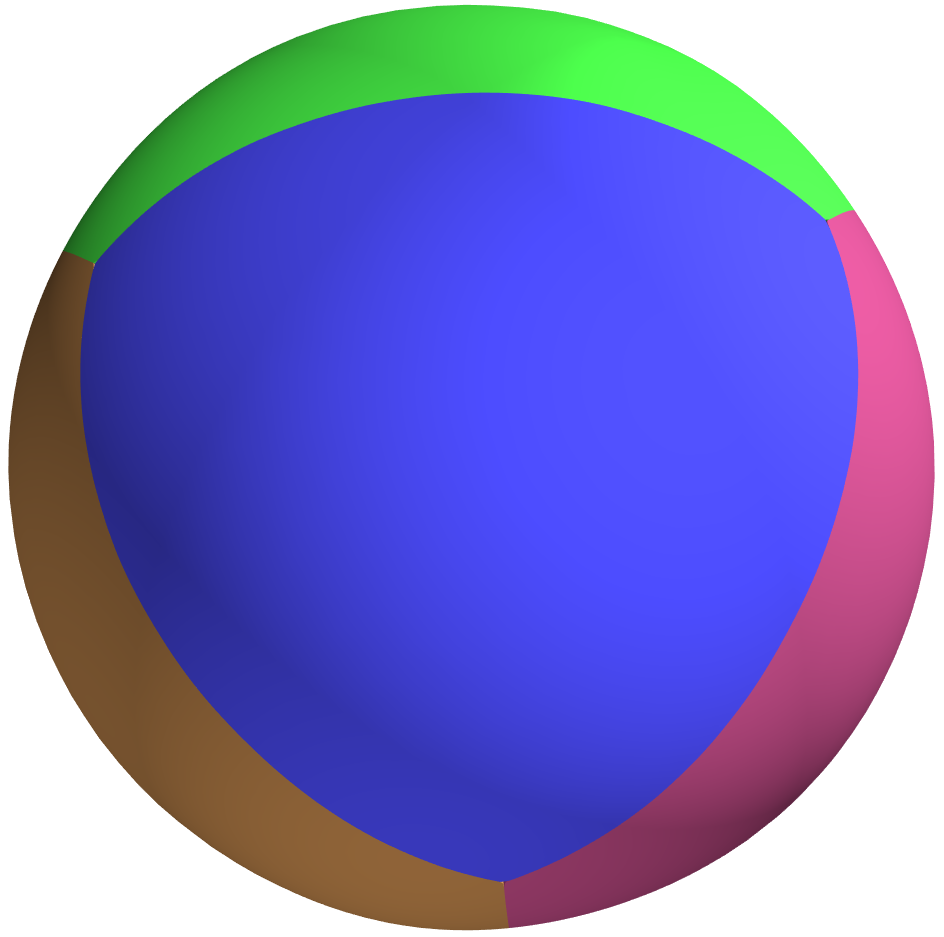}\qquad
	\caption{(left) Maximization of the length of the minimal perimeter partition into equal areas for $n \in \{3,4\}$. (right) Results obtained when the area constraints are not the same: $n=3$: ratios $1:2:2$, $n=4$: ratios $1:2:2:2$.}
	\label{fig:part-3D-equal}
\end{figure}

\begin{figure}
	\centering
	\includegraphics[width=0.2\textwidth]{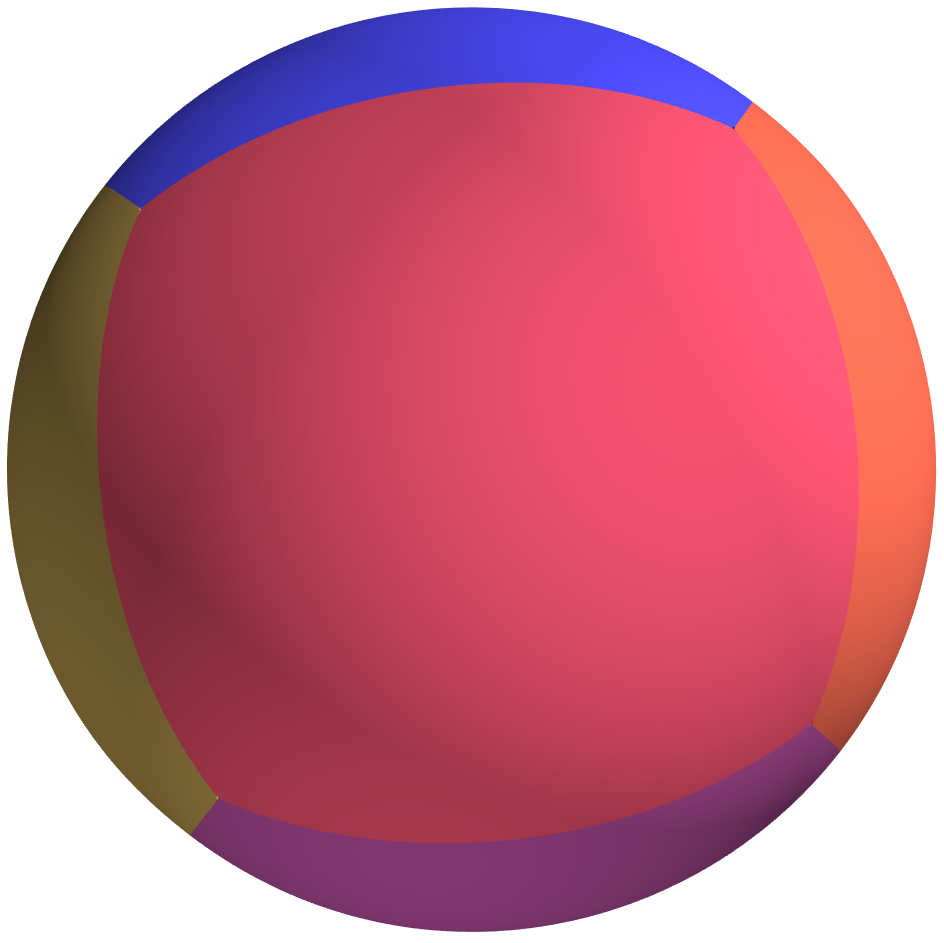}\quad
	\includegraphics[width=0.2\textwidth]{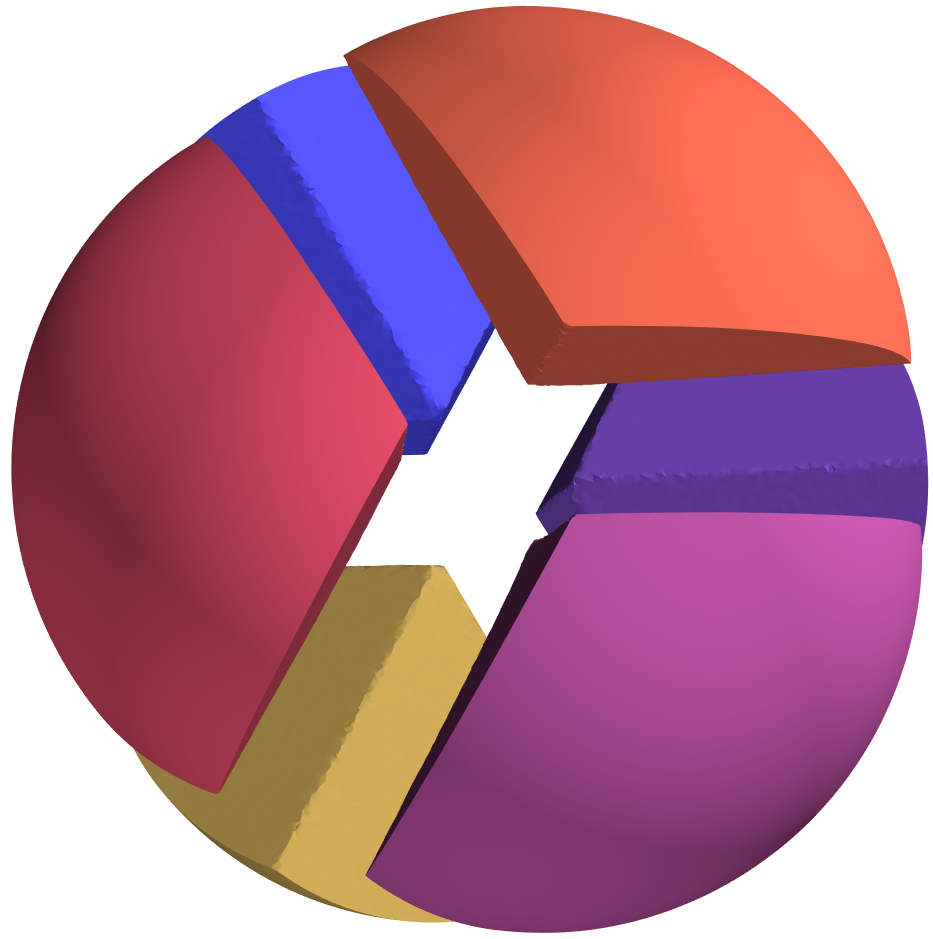}\quad
	\includegraphics[width=0.2\textwidth]{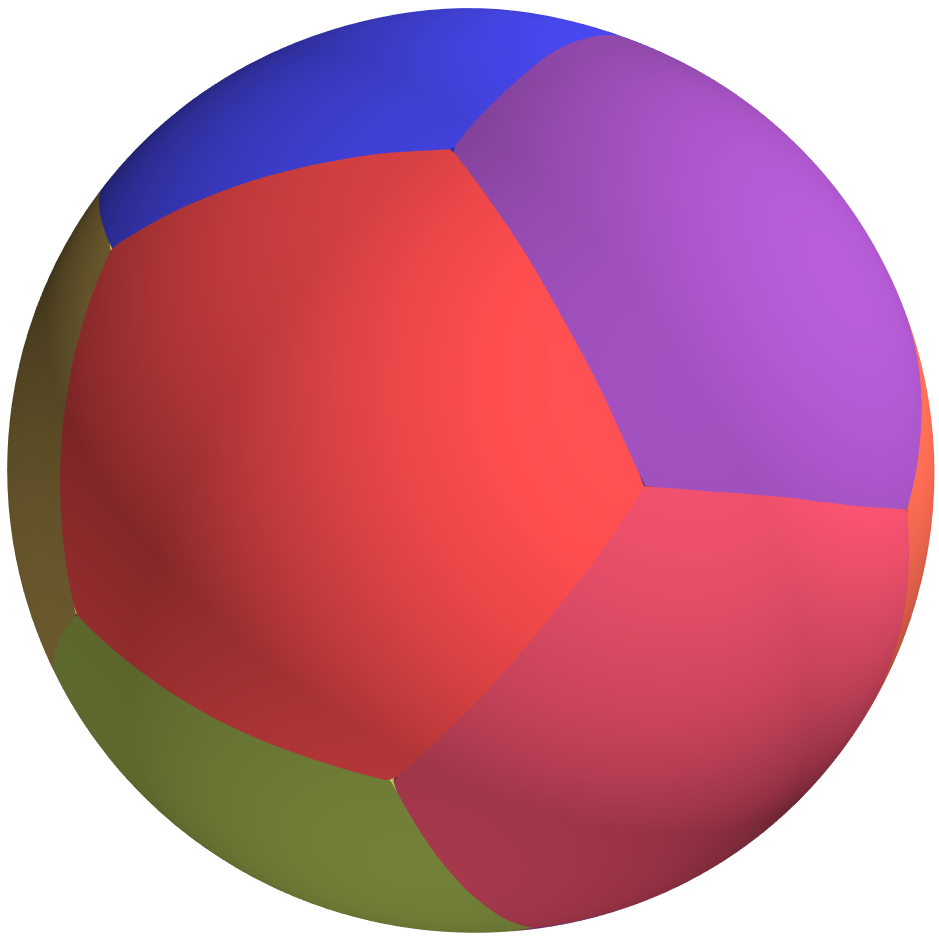}\quad
	\includegraphics[width=0.2\textwidth]{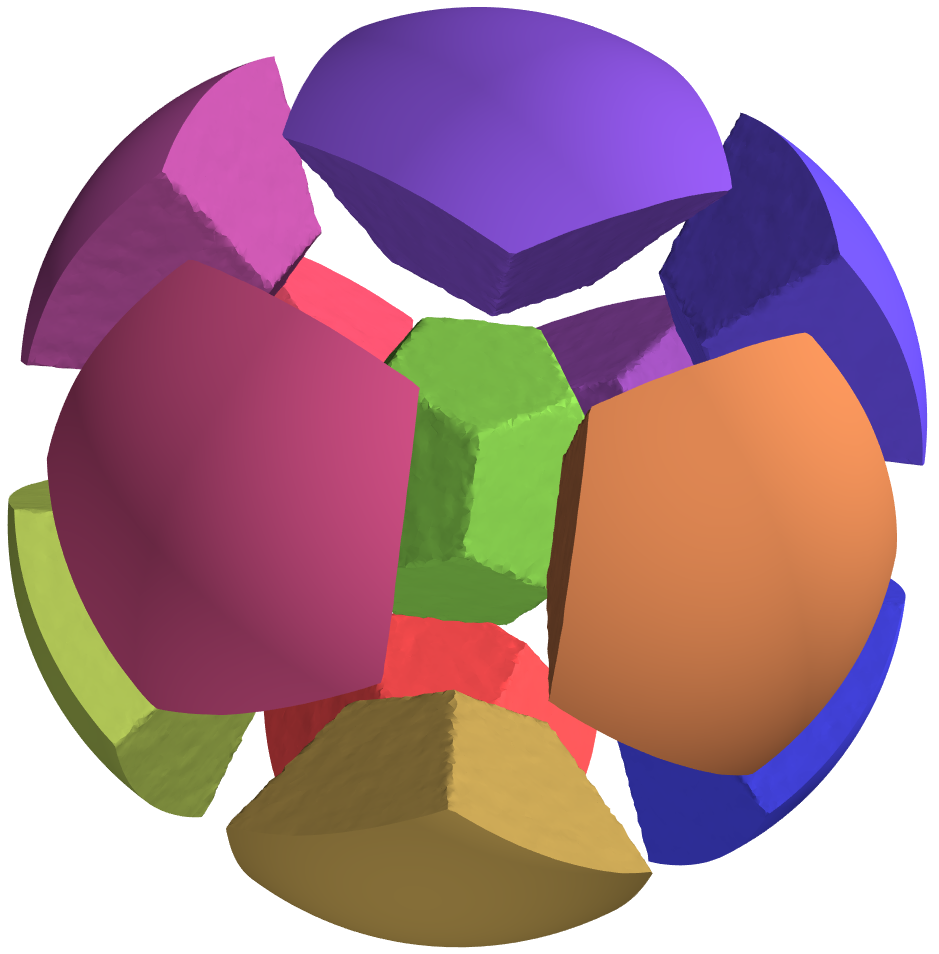}
	\caption{Maximization of the length of the minimal perimeter partition into equal areas for $n \in \{6,13\}$. An expanded view of the optimal partition is also illustrated for each case.}
	\label{fig:part-3D-large}
\end{figure}

Surprisingly, the case of partitions shows similar results. When considering equal area constraints the set with fixed area maximizing the length of the minimal partition is still the disk (see Figure \ref{fig:part-2D-equal} for some examples). In dimension three for $n\in \{3,4,6,13\}$ we obtain similar results: the ball maximizes the total surface area of the smallest total perimeter partition. These results are illustrated in Figures \ref{fig:part-3D-equal} and \ref{fig:part-3D-large}.

\begin{figure}
	\centering 
	\includegraphics[width=0.3\textwidth]{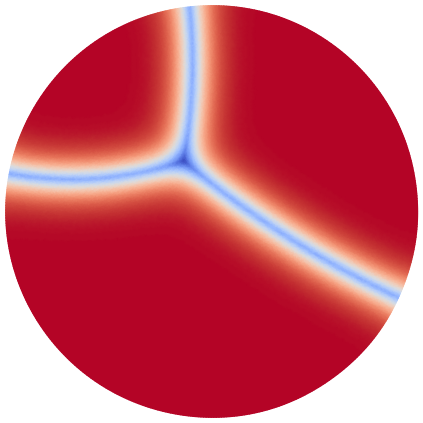}\quad
	\includegraphics[width=0.3\textwidth]{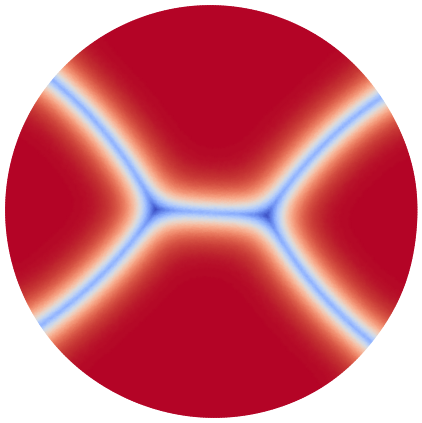}\quad
	\includegraphics[width=0.3\textwidth]{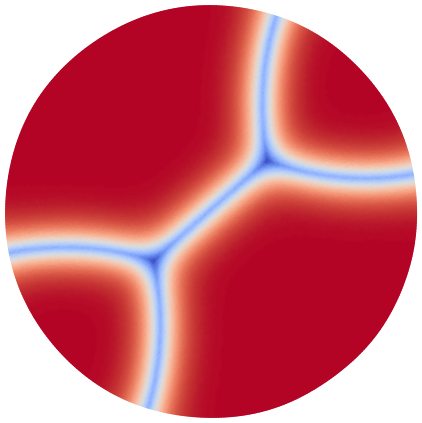}\quad
	\caption{Maximization of the length of the minimal perimeter partition into different areas: $n=3$, ratios $1:2:3$, $n=4$, ratios $1:1:2:2$, $n=4$, ratios $1:1:3:3$.}
	\label{fig:part-2D-non-equal}
\end{figure}

Note that simulations made in the case of one phase already show that when partitioning the domain into two regions with two non-equal areas, the maximizer of the minimal length partition is still the disk. This suggests that even in the case where the cells do not have the same prescribed area, the set $\Omega$ which maximizes the minimal perimeter of a partition is still the disk in 2D (the ball in 3D). Indeed, when considering more cells with different areas, the numerical result is the same: the disk seems to be the maximizer (see Figure \ref{fig:part-2D-non-equal} for some examples). As already underlined in \cite{cox-nonequal}, the study of partitions in cells with prescribed but different areas is more complex, since in this case there are even more local minima. 

The numerical simulations give rise to the following conjecture, which generalizes the results of \cite{longest_shortest_fence}.
\begin{conj}
	
   1. Given $c \in (0,1)$, the set $\Omega$ maximizing $I(\Omega,c)$ under the constraint $|\Omega|=v_d$ (i.e. solving \eqref{eq:longest_shortest_fence}) is the ball.
   
   2. Given $n>1$ and $\bo c = (c_i)_{i=1}^n \in \Bbb{R}^n_+$ with $\sum_{i=1}^n c_i=1$, the set $\Omega$ maximizing $PI(\Omega,\bo c)$ under the constraint $|\Omega|=v_d$ (i.e. solving \eqref{eq:longest_shortest_partition}) is the ball.
\end{conj}

\begin{rem}
	The same type of results seem to hold when maximizing the minimal geodesic perimeter for closed surfaces in $\Omega$ in $\Bbb{R}^3$ which are boundaries of convex sets with a constraint on the $\mathcal H^2(\partial \Omega)$. The techniques used in this case are those presented in \cite{bogosel_oudet_manifolds} and the theoretical and numerical framework is similar to what was done in dimension three. In this case the sphere seems to be the maximizing set, which is in accord with the conjecture stated above. 
\end{rem}

\section{Remarks on optimality conditions}
\label{sec:optim-cond}

As discussed in Section \ref{sec:sh-deriv}, existence of shape derivatives for $I_\varepsilon $ and $PI_\varepsilon $ depends on the uniqueness of the minimizers for these functionals. Therefore, it is not straightforward to obtain classical optimality conditions. It is possible, however, to obtain some qualitative information about sets maximizing the minimal values of $I_\varepsilon $ and $PI_\varepsilon $ under volume constraint. Recall that the optimizer of shape differentiable functional $J$ under volume constraint will verify an optimality condition of the form
\begin{equation}
J'(\Omega)(\theta)+\ell |\Omega|'(\theta) = 0,
\label{eq:lagrange-volume}
\end{equation}
where $\ell \in \Bbb{R}$ is a Lagrange multiplier associated to the volume constraint. Recall that the shape derivative of the volume functional is $|\Omega|'(\theta) = \int_{\partial \Omega} \theta.n$. 

{\bf Non-uniqueness of the minimal relative perimeter set/partition at the optimum.} Results in Section \ref{sec:sh-deriv} indicate that the shape derivatives of $I_\varepsilon $ and $PI_\varepsilon $ exist when they correspond to unique minimizers of the Modica-Mortola type functionals. In this case, the corresponding shape derivatives are boundary integrals of non-constant functions multiplied by the normal perturbation $\theta.n$. Therefore, it is straightforward to see that a relation of the type \eqref{eq:lagrange-volume} cannot hold. This allows us to conclude that of $\Omega^*$ is a minimizer of $I_\varepsilon (\Omega,c)$ the optimal minimal relative perimeter set is not unique. The same happens in the case of partitions: if $\Omega^*$ minimizes $PI_\varepsilon (\Omega,\bo c)$ and the minimal length partition of $\Omega$ with constraints $\bo c$ is not unique. 

Suppose now that $\Omega$ is a domain with fixed volume $|\Omega| = v_d$ such that $SF(\Omega,c)$ is unique. Then for $\varepsilon>0$ small enough the minimizer of $I_\varepsilon (\Omega,c)$ will also be unique. Thus $I_\varepsilon (\Omega,c)$ admits a shape derivative. Also the corresponding optimal density $u_\Omega$ is not constant on the boundary, and therefore the optimality relation \eqref{eq:lagrange-volume} cannot hold. This shows that such a domain $\Omega$ is not a solution of problem \eqref{eq:longest_shortest_fence}. The same argument can be applied for problem \eqref{eq:longest_shortest_partition}.

As a conclusion, a domain $\Omega$ that solves \eqref{eq:longest_shortest_fence} must have multiple minimal relative perimeter sets (or multiple minimal length partitions for problem \eqref{eq:longest_shortest_partition}).

\section{Conclusions}
The theoretical considerations and numerical simulations presented in this paper suggest that the results of \cite{berry-et-al}, \cite{wang-wang}, \cite{longest_shortest_fence} are valid in more general settings: in dimensions two and three, under volume and convexity constraints the ball is the set $\Omega$ who maximizes
\begin{itemize}[topsep=0pt,noitemsep]
	\item the minimal relative perimeter of a subset $\omega\subset \Omega$ with volume constraint $|\omega| = c|\Omega|$ for all $c \in (0,1)$.
	\item the minimal relative perimeter of a partition of $\Omega$ into sets $(\omega_i)_{i=1}^n$ with volume constraints $|\omega_i|=c_i|\Omega|$ given $c_i\in (0,1)$ with $\sum_{i=1}^n c_i=1$. The result seems to hold even in the case where the sets $|\omega_i|$ do not have the same volume constraints. 
\end{itemize}

The numerical maximization algorithm consists in solving at each iteration an optimization problem which approximates the least perimeter set or partition under the given constraint. Then a perturbation of the set which does not decrease the minimal perimeter is found and the set is modified. In all cases, the numerical result was close to the disk/ball. 

The initialization phase for the computation of the optimal partitions is made using Voronoi diagrams with prescribed capacity. We provide a new way of generating such Voronoi diagrams using the gradients of the areas with respect to the Voronoi points. The gradient of the perimeter of the Voronoi cells is also computed.

\begin{center}
	\textsc{Acknowledgments}
\end{center}

The authors were partially supported by the project ANR-18-CE40-0013 SHAPO financed by the French Agence Nationale de la Recherche (ANR). The authors thank Frank Morgan for indicating more references to previous works related to the convex isoperimetric problem. The authors thank the reviewers for the careful reading of the manuscript and for their remarks which helped improve the quality of the article.

\bibliography{./max_rel_per}

\begin{thebibliography}{10}

\bibitem{gammaconvalberti}
G.~Alberti.
\newblock Variational models for phase transitions, an approach via
  {$\Gamma$}-convergence.
\newblock In {\em Calculus of variations and partial differential equations
  ({P}isa, 1996)}, pages 95--114. Springer, Berlin, 2000.

\bibitem{ambrosiofuscopallara}
L.~Ambrosio, N.~Fusco, and D.~Pallara.
\newblock {\em Functions of bounded variation and free discontinuity problems}.
\newblock Oxford Mathematical Monographs. The Clarendon Press, Oxford
  University Press, New York, 2000.

\bibitem{capacity-constrained-continuous}
M.~{Balzer}.
\newblock Capacity-constrained voronoi diagrams in continuous spaces.
\newblock In {\em 2009 Sixth International Symposium on Voronoi Diagrams},
  pages 79--88, 2009.

\bibitem{capacity-constrained-discrete}
M.~Balzer, T.~Schl\"{o}mer, and O.~Deussen.
\newblock Capacity-constrained point distributions: A variant of lloyd's
  method.
\newblock In {\em ACM SIGGRAPH 2009 Papers}, SIGGRAPH '09, New York, NY, USA,
  2009. Association for Computing Machinery.

\bibitem{berry-et-al}
J.~Berry, E.~Bongiovanni, W.~Boyer, B.~Brown, P.~Gallagher, D.~Hu, A.~Loving,
  Z.~Martin, M.~Miller, B.~Perpetua, and S.~Tammen.
\newblock The convex body isoperimetric conjecture in the plane.
\newblock {\em Rose-Hulman Undergraduate Mathematics Journal}, 18(2), 2017.
\newblock \nolinkurl{https://scholar.rose-hulman.edu/rhumj/vol18/iss2/2}.

\bibitem{bogosel_oudet_manifolds}
B.~Bogosel and E.~Oudet.
\newblock Partitions of minimal length on manifolds.
\newblock {\em Exp. Math.}, 26(4):496--508, 2017.

\bibitem{bonnans}
J.~F. Bonnans and A.~Shapiro.
\newblock {\em Perturbation analysis of optimization problems}.
\newblock Springer Series in Operations Research. Springer-Verlag, New York,
  2000.

\bibitem{braides2}
A.~Braides.
\newblock {\em Approximation of {F}ree-{D}iscontinuity Problems}.
\newblock Springer, 1998.

\bibitem{many_cells}
E.~Bretin, R.~Denis, J.-O. Lachaud, and E.~Oudet.
\newblock Phase-field modelling and computing for a large number of phases.
\newblock {\em ESAIM Math. Model. Numer. Anal.}, 53(3):805--832, 2019.

\bibitem{extension_theorems}
V.~Burenkov.
\newblock Extension theorems for {S}obolev spaces.
\newblock In {\em The {M}az'ya anniversary collection, {V}ol. 1 ({R}ostock,
  1998)}, volume 109 of {\em Oper. Theory Adv. Appl.}, pages 187--200.
  Birkh\"{a}user, Basel, 1999.

\bibitem{buttazzogconv}
G.~Buttazzo.
\newblock Gamma-convergence and its {A}pplications to {S}ome {P}roblems in the
  {C}alculus of {V}ariations.
\newblock {\em School on {H}omogenization {ICTP}, {T}rieste, {S}eptember 6-17},
  1993.

\bibitem{sh-deriv-gauss}
A.~Chicco-Ruiz, P.~Morin, and M.~S. Pauletti.
\newblock The shape derivative of the {G}auss curvature.
\newblock {\em Rev. Un. Mat. Argentina}, 59(2):311--337, 2018.

\bibitem{cianchi}
A.~Cianchi.
\newblock On relative isoperimetric inequalities in the plane.
\newblock {\em Boll. Un. Mat. Ital. B (7)}, 3(2):289--325, 1989.

\bibitem{cox-partitions}
S.~J. Cox and E.~Flikkema.
\newblock The minimal perimeter for {$N$} confined deformable bubbles of equal
  area.
\newblock {\em Electron. J. Combin.}, 17(1):Research Paper 45, 23, 2010.

\bibitem{MMG}
C.~Dapogny, C.~Dobrzynski, and P.~Frey.
\newblock Three-dimensional adaptive domain remeshing, implicit domain meshing,
  and applications to free and moving boundary problems.
\newblock {\em J. Comput. Phys.}, 262:358--378, 2014.

\bibitem{delfour-zolesio}
M.~C. Delfour and J.-P. Zol\'{e}sio.
\newblock {\em Shapes and geometries}, volume~22 of {\em Advances in Design and
  Control}.
\newblock Society for Industrial and Applied Mathematics (SIAM), Philadelphia,
  PA, second edition, 2011.
\newblock Metrics, analysis, differential calculus, and optimization.

\bibitem{Du1999}
Q.~Du, V.~Faber, and M.~Gunzburger.
\newblock Centroidal voronoi tessellations: Applications and algorithms.
\newblock {\em {SIAM} Review}, 41(4):637--676, Jan. 1999.

\bibitem{longest_shortest_fence}
L.~Esposito, V.~Ferone, B.~Kawohl, C.~Nitsch, and C.~Trombetti.
\newblock The longest shortest fence and sharp {P}oincar\'{e}-{S}obolev
  inequalities.
\newblock {\em Arch. Ration. Mech. Anal.}, 206(3):821--851, 2012.

\bibitem{fiacco}
A.~V. Fiacco.
\newblock {\em Introduction to sensitivity and stability analysis in nonlinear
  programming}, volume 165 of {\em Mathematics in Science and Engineering}.
\newblock Academic Press, Inc., Orlando, FL, 1983.

\bibitem{Firey1965}
W.~J. Firey.
\newblock Lower bounds for volumes of convex bodies.
\newblock {\em Archiv der Mathematik}, 16(1):69--74, Dec. 1965.

\bibitem{gilbarg-trudinger}
D.~Gilbarg and N.~S. Trudinger.
\newblock {\em Elliptic partial differential equations of second order}.
\newblock Classics in Mathematics. Springer-Verlag, Berlin, 2001.
\newblock Reprint of the 1998 edition.

\bibitem{grisvard}
P.~Grisvard.
\newblock {\em Elliptic problems in nonsmooth domains}, volume~24 of {\em
  Monographs and Studies in Mathematics}.
\newblock Pitman (Advanced Publishing Program), Boston, MA, 1985.

\bibitem{gurtin-matano}
M.~E. Gurtin and H.~Matano.
\newblock On the structure of equilibrium phase transitions within the gradient
  theory of fluids.
\newblock {\em Quart. Appl. Math.}, 46(2):301--317, 1988.

\bibitem{cox-nonequal}
F.~Headley and S.~Cox.
\newblock Least perimeter partition of the disc into n bubbles of two different
  areas.
\newblock {\em The European Physical Journal E}, 42(7), July 2019.

\bibitem{freefem}
F.~Hecht.
\newblock {N}ew development in {F}ree{F}em++.
\newblock {\em J. Numer. Math.}, 20(3-4):251--265, 2012.

\bibitem{henrot-pierre-english}
A.~Henrot and M.~Pierre.
\newblock {\em Shape variation and optimization}, volume~28 of {\em EMS Tracts
  in Mathematics}.
\newblock European Mathematical Society (EMS), Z\"{u}rich, 2018.
\newblock A geometrical analysis, English version of the French publication [
  MR2512810] with additions and updates.

\bibitem{Iri}
M.~Iri, K.~Murota, and T.~Ohya.
\newblock A fast voronoi-diagram algorithm with applications to geographical
  optimization problems.
\newblock In {\em System Modelling and Optimization}, pages 273--288.
  Springer-Verlag.

\bibitem{nlopt}
S.~G. Johnson.
\newblock The nlopt nonlinear-optimization package.
\newblock \nolinkurl{http://github.com/stevengj/nlopt}.

\bibitem{image_des_maths}
B.~Kloeckner.
\newblock Dans quelle forme la plus petite paroi enfermant un volume donn\'{e}
  est-elle la plus grande ?
\newblock - Images des mathématiques :
  \nolinkurl{http://images.math.cnrs.fr/Dans-quelle-forme-la-plus-petite-paroi-enfermant-un-volume-donne-est-elle-la.html},
  12 2019.

\bibitem{loomis_whitney}
A.~Koldobsky, C.~Saroglou, and A.~Zvavitch.
\newblock Estimating volume and surface area of a convex body via its
  projections or sections.
\newblock {\em Studia Math.}, 244(3):245--264, 2019.

\bibitem{luckhaus-modica}
S.~Luckhaus and L.~Modica.
\newblock The {G}ibbs-{T}hompson relation within the gradient theory of phase
  transitions.
\newblock {\em Arch. Rational Mech. Anal.}, 107(1):71--83, 1989.

\bibitem{maggi-book}
F.~Maggi.
\newblock {\em Sets of finite perimeter and geometric variational problems},
  volume 135 of {\em Cambridge Studies in Advanced Mathematics}.
\newblock Cambridge University Press, Cambridge, 2012.
\newblock An introduction to geometric measure theory.

\bibitem{modica_gradient_theory}
L.~Modica.
\newblock The gradient theory of phase transitions and the minimal interface
  criterion.
\newblock {\em Arch. Rational Mech. Anal.}, 98(2):123--142, 1987.

\bibitem{modica-mortola}
L.~Modica and S.~Mortola.
\newblock Un esempio di {$\Gamma ^{-}$}-convergenza.
\newblock {\em Boll. Un. Mat. Ital. B (5)}, 14(1):285--299, 1977.

\bibitem{morgan-bubbles}
F.~Morgan.
\newblock Soap bubbles in {${\bf R}^2$} and in surfaces.
\newblock {\em Pacific J. Math.}, 165(2):347--361, 1994.

\bibitem{morgan_blog}
F.~Morgan.
\newblock Convex body isoperimetric conjecture.
\newblock
  \nolinkurl{https://sites.williams.edu/Morgan/2010/07/03/convex-body-isoperimetric-conjecture/},
  2010.

\bibitem{oudet}
{\'E}.~Oudet.
\newblock Approximation of partitions of least perimeter by
  {$\Gamma$}-convergence: around {K}elvin's conjecture.
\newblock {\em Exp. Math.}, 20(3):260--270, 2011.

\bibitem{polya}
G.~Polya.
\newblock Aufgabe 283.
\newblock {\em Elem. d. Math.}, 13:40--41, 1958.

\bibitem{mayavi}
P.~{Ramachandran} and G.~{Varoquaux}.
\newblock Mayavi: 3d visualization of scientific data.
\newblock {\em Computing in Science Engineering}, 13(2):40--51, 2011.

\bibitem{isoperimetric_convex}
M.~Ritor\'{e} and E.~Vernadakis.
\newblock Isoperimetric inequalities in {E}uclidean convex bodies.
\newblock {\em Trans. Amer. Math. Soc.}, 367(7):4983--5014, 2015.

\bibitem{sokolowski-zolesio}
J.~Soko\l~owski and J.-P. Zol\'{e}sio.
\newblock {\em Introduction to shape optimization}, volume~16 of {\em Springer
  Series in Computational Mathematics}.
\newblock Springer-Verlag, Berlin, 1992.
\newblock Shape sensitivity analysis.

\bibitem{sternberg}
P.~Sternberg and K.~Zumbrun.
\newblock Connectivity of phase boundaries in strictly convex domains.
\newblock {\em Arch. Rational Mech. Anal.}, 141(4):375--400, 1998.

\bibitem{Svanberg2002}
K.~Svanberg.
\newblock A class of globally convergent optimization methods based on
  conservative convex separable approximations.
\newblock {\em {SIAM} Journal on Optimization}, 12(2):555--573, Jan. 2002.

\bibitem{wang-wang}
B.-H. Wang and Y.-K. Wang.
\newblock A note on the convex body isoperimetric conjecture in the plane,
  2021.

\bibitem{wichiramala}
W.~Wichiramala.
\newblock Efficient cut for a subset of prescribed area.
\newblock {\em Thai J. Math.}, 5(3, Special issue):95--100, 2007.

\bibitem{xin2016}
S.-Q. Xin, B.~L\'{e}vy, Z.~Chen, L.~Chu, Y.~Yu, C.~Tu, and W.~Wang.
\newblock Centroidal power diagrams with capacity constraints: Computation,
  applications, and extension.
\newblock {\em ACM Trans. Graph.}, 35(6), Nov. 2016.

\end{thebibliography}
\bibliographystyle{abbrv}

\end{document}